\magnification=1200
\input amstex
\UseAMSsymbols

\noindent
\space

\centerline{\bf  ESSENTIAL COMMUTANTS}
\centerline{\bf  ON STRONGLY PSEUDO-CONVEX DOMAINS}

\bigskip
\bigskip

\centerline{\bf  Yi Wang${}^{\bold 1}$ and Jingbo Xia}

{\input amsppt.sty
\topmatter
\thanks
\vskip0pt
${}^1$  Supported in part by National Science Foundation grant DMS-1900076.
\endthanks
\endtopmatter
}

\noindent
{\bf Abstract.}    
 Consider a bounded strongly pseudo-convex domain  $\Omega $  with a smooth boundary in  ${\bold C}^n$.  
Let  ${\Cal T}$ be the Toeplitz algebra on the Bergman space  $L^2_a(\Omega )$.   That is,   ${\Cal T}$  
is the  $C^\ast $-algebra generated by the Toeplitz operators  $\{T_f : f \in L^\infty (\Omega )\}$.  
Extending the work  [27,28] in the special case of the unit ball,  we show that on any such  $\Omega $,  
${\Cal T}$   and   $\{T_f : f \in  {\text{VO}}_{\text{bdd}}\} + {\Cal K}$  are essential commutants of each other.   
On a general  $\Omega $ considered in this paper,  the proofs require many new ideas and techniques.
These same techniques also enable us to show that for  $A \in {\Cal T}$,  if  $\langle Ak_z,k_z\rangle \rightarrow 0$
as  $z \rightarrow \partial \Omega $,  then  $A$  is a compact operator.   

\bigskip
\centerline{\bf  1. Introduction}
\medskip
An enduring question in the study of Toeplitz operators is their essential commutativity.   
In this paper we consider this question on strongly pseudo-convex domains.  
It will be beneficial to start the paper with a recollection of necessary definitions and background.
\medskip
Suppose that   ${\Cal Z}$   is a   collection of bounded operators on a Hilbert space  
${\Cal H}$.   Then its essential commutant   is defined to be
$$
\text{EssCom}({\Cal Z})  =  \{A \in {\Cal B}({\Cal H}) : [A,T] \ \ \text{is compact for every} \ \  T \in {\Cal Z}\}.
$$
\noindent
The study of essential commutants began with the classic papers  of Johnson-Parrott [13],  Voiculescu [23]  and Popa  [19].  
Ever since, essential-commutant problems have become a mainstay of operator theory and operator algebras.  
As it turns out,   many of the most interesting examples in the study of essential commutants are associated with Toeplitz operators, 
of various kinds [4,6-9,11,24,25,27,28].
Perhaps, one reason why essential-commutant problems attract attention is that they are generally not easy.
\medskip
In this paper we consider an arbitrary  bounded, strongly pseudo-convex domain  $\Omega $  with a smooth boundary in  ${\bold C}^n$.  
Recall that the Bergman space   $L^2_a(\Omega )$   is the collection of analytic functions  $h$  on   $\Omega $  satisfying the condition 
$$
\int _\Omega |h|^2dv < \infty ,
$$
\noindent
where  $dv$  is the volume measure on   $\Omega $.   Let  $P : L^2(\Omega ) \rightarrow L^2_a(\Omega )$  be the orthogonal projection.  
For each  $f \in L^\infty (\Omega )$,  we have the Toeplitz operator   $T_f$  defined by the formula
$$
T_fh  =  P(fh),   \quad  h \in L^2_a(\Omega ).
$$
\noindent
Let  ${\Cal T}$  denote the  $C^\ast $-algebra generated by   $\{T_f : f \in  L^\infty (\Omega )\}$.   Then  ${\Cal T}$  is called the 
{\it Toeplitz algebra}  on  the Bergman space   $L^2_a(\Omega )$.   It is well known that   ${\Cal T}$  contains   ${\Cal K}$,  the collection of 
compact operators on  $L^2_a(\Omega )$  [22,Theorem 4.1.25].  Obviously,  this is a convenient fact for the study of essential commutants.
\medskip
In the case of the unit ball  {\bf B}  in  ${\bold C}^n$,    
the essential commutant problems related to   ${\Cal T}$   were solved in [27,28],  with [26,Theorem 1.3] playing a pivotal role.  
Specifically, in the case of the unit ball,   it was shown that EssCom$({\Cal T}) = \{T_f : f \in  {\text{VO}}_{\text{bdd}}\} + {\Cal K}$   in [27]  
and that  EssCom$\{T_f : f \in  {\text{VO}}_{\text{bdd}}\} = {\Cal T}$  in  [28].   Once one knows that, a question naturally presents itself:  what happens if 
one replaces the unit ball   {\bf B}  by a general strongly pseudo-convex domain  $\Omega $?  Equally naturally, one would expect that the same 
results hold on a general  $\Omega $.  But here one immediately runs into two difficulties:
\medskip
(1)  The works in [27,28], particularly in [26],  rely heavily on the explicit formula for the Bergman metric  $\beta $  on  {\bf B}.  Without such 
an explicit formula, it is not clear how to redo many of the estimates in [27,28].  By contrast,  in the case of a general strongly 
pseudo-convex domain,  we only know the  {\it asymptotics}  of the Bergman metric [10,20],  but we do not have a formula for it that is explicit enough.
In other words, on a general   $\Omega $,  we do not have good enough a handle on the Bergman metric to do 
many of the  necessary estimates.  The same is true if one considers the Kobayashi metric instead of the Bergman metric.
\medskip
(2)  The techniques in [26,27,28] depend heavily on the M\"obius transforms on  {\bf B}.   But on a general 
strongly pseudo-convex domain  $\Omega $,  there is no such thing as M\"obius transform.   In other words,  compared 
with the unit ball,  a general    $\Omega $  totally lacks global symmetry.   Compared with (1), this difficulty is more substantive,  but it 
also makes an exciting challenge:  can we prove the results in [27,28] on a domain without symmetry?
\medskip
We are pleased to report that we have managed to overcome these difficulties.  
The way we deal with difficulty (1) is to simply introduce a new metric that serves our purpose.   Since  $\Omega $  
is a strongly pseudo-convex domain,  it has a defining function  $r$,  i.e.,  $\Omega = \{z \in {\bold C}^n : r(z) < 0\}$.   
Then the formula
$$
b_{ij}(z)  =    \partial _i\bar \partial _j\log {1\over -r(z)},   \quad  1 \leq i,j \leq n,  
$$
\noindent
for  $z$  near  $\partial \Omega $  gives us  the infinitesimal generator of a metric  $d$  on  $\Omega $.   One might 
call this  $d$  a poor man's  imitation of the Bergman metric,  but the above formula is explicit enough 
to allow us to do all the necessary analysis.     We will have more to say about this point below,  
and the precise definition of  $d$  will be given at the beginning of Section 2.
\medskip
Difficulty (2) simply requires new approaches.  
Examining the involvements of M\"obius transforms in [27,28] one by one,  
we have managed to find a new idea or new technique as a replacement in each case.  Thus the results 
about essential commutants mentioned above can indeed be proved without symmetry.
\medskip
To state our results, we need the notion of   {\it vanishing oscillation},  which was first introduced  
in [5,3] for functions on bounded symmetric domains with respect to the Bergman metric.  In this paper 
we need to define functions of vanishing oscillation with respect to the metric  $d$  on  $\Omega $.  
Let  $f$  be a continuous function  on  $\Omega $.  Then   $f$  is said to have vanishing oscillation  if  
$$
\lim _{z\rightarrow \partial \Omega }\sup \{|f(z) - f(w)| : d(z,w) \leq 1\}   =   0.
$$
\noindent
Let  VO  denote the collection of  functions of vanishing oscillation  on  $\Omega $.  Further,  define
$$
{\text{VO}}_{\text{bdd}}  =  \text{VO}\cap L^\infty (\Omega ).
$$
\noindent
Our main results are the two theorems below:
\medskip
\noindent
{\bf Theorem 1.1.}     {\it On any  bounded},  {\it strongly pseudo-convex domain}   $\Omega $  
{\it with a smooth boundary in}  ${\bold C}^n$,  {\it the following hold true}:

\noindent
(i)   {\it The Toeplitz algebra}   ${\Cal T}$  {\it is the essential commutant of}  
$\{T_f : f \in {\text{VO}}_{\text{bdd}}\}$.

\noindent
(ii)  {\it The  essential commutant of}  ${\Cal T}$  {\it equals}   $\{T_f : f \in {\text{VO}}_{\text{bdd}}\} + {\Cal K}$.

\medskip
Let    ${\Cal Q}$   denote the Calkin algebra    ${\Cal B}(L^2_a(\Omega ))/{\Cal K}$,   and let    
$$
\pi : {\Cal B}(L^2_a(\Omega )) \rightarrow {\Cal Q}
$$   
\noindent
be the quotient homomorphism.    Then    $\pi (\text{EssCom}({\Cal Z})) = \{\pi ({\Cal Z})\}'$
for every subset   ${\Cal Z} \subset {\Cal B}(L^2_a(\Omega ))$.   
Obviously,  a subset   ${\Cal A}$   of   ${\Cal Q}$  satisfies the double-commutant relation  
${\Cal A} = {\Cal A}''$   if and only if   ${\Cal A} =  {\Cal G}'$   for some  ${\Cal G} \subset {\Cal Q}$.
Thus Theorem 1.1(i)  implies  that   $\pi ({\Cal T})$  satisfies the double-commutant relation in  ${\Cal Q}$.
\medskip 
As it turns out, the techniques that allow us to prove Theorem 1.1(i),  also  
give us a classic compactness criterion for  $A \in {\Cal T}$  in terms of its Berezin transform on  $\Omega $.   
Let us  write   $k_z$,  $z \in \Omega $,    for the normalized reproducing kernel for the 
Bergman space  $L^2_a(\Omega )$.   
\medskip
\noindent
{\bf Theorem 1.2.}     {\it Consider any bounded},  {\it strongly pseudo-convex domain}   $\Omega $  
{\it with a smooth boundary in}  ${\bold C}^n$.   {\it Let}  $A \in {\Cal T}$.  {\it If}
$$
\lim _{z\rightarrow \partial \Omega }\langle Ak_z,k_z\rangle = 0,
\tag 1.1
$$
\noindent
{\it then}   $A$     {\it is a compact operator on}   $L^2_a(\Omega )$.

\medskip
At this point, it is appropriate to briefly recall the long history of this line of investigations.   
The first result of this genre was due to Axler and Zheng [1],  where the domain was 
the unit disc in   {\bf C}   and   $A$  was a finite algebraic combination of Toeplitz operators.   
Later in [21], Su\'arez showed that this compactness criterion holds for all   $A \in {\Cal T}$  on the unit  
ball  {\bf B}  in    ${\bold C}^n$.   The fact that Su\'arez was able to do this for {\it arbitrary}  $A \in {\Cal T}$  on the ball,  
rather than just for finite algebraic combinations of Toeplitz operators,  was considered to be a major breakthrough.   
Consequently,  [21]  inspired many generalizations [2,12,29], including generalizations on the Fock space.   But all these papers depend on   
the M\"obius transforms on the domain in question.   In this regard,   Theorem 1.2 is the  
first to remove any and all involvement of M\"obius transforms, since in general there aren't any on  $\Omega $.
\medskip
The rest of the paper is taken up by the  proofs of these results.  Because we have to start from scratch,  
there are numerous steps involved.  We conclude the introduction by an outline of our plan.  
\medskip
First of all, in Section 2 we precisely define the metric  $d$  mentioned above.   In addition to  $d$,  another important 
quantity for the paper is the  ``gauge"
$$
\rho (z,w) = |z - w|^2 + |\langle z - w,(\bar \partial r)(z)\rangle |
$$ 
\noindent
on  $\Omega $.  Section 2 contains several fundamental estimates involving  $d$  and  $\rho (z,w)$.   Section 3  brings in 
another important ingredient for our analysis,  the function  
$$
F(z,w) = |r(z)| + |r(w)| + \rho (z,w),
$$
\noindent
which is a familiar fixture on strongly pseudo-convex domains.   The main result of the section 
is  Lemma 3.8,  which is a version of the Forelli-Rudin estimates  for  $\Omega $
in which  $d$  and  $F$  are quantitatively involved.
\medskip
Sections 4 and 5 are devoted to operators that are discrete sums constructed from the Bergman kernel  $K(z,w)$
over  $d$-lattices.   The main goal for these two sections is Corollary 5.3,   which provides the norm-continuity 
of such discrete sums under small perturbation of the lattice.
\medskip
In Section 6, we introduce   LOC$(A)$,  the class of ``localized versions of  $A$"  for any bounded 
operator  $A$  on   $L^2_a(\Omega )$.   Using Lemma 3.8 and Corollary 5.3 mentioned above and doing 
quite a bit of additional work,  we show in Section 6  that   LOC$(A) \subset {\Cal T}$  for every  
$A \in {\Cal B}(L^2_a(\Omega ))$.  This is a major step in the proof of Theorem 1.1(i).
\medskip
Section 7 is devoted to matters related to functions of vanishing oscillation.   
In particular,  we consider the scalar quantity
$$
\text{diff}(f) = \sup \{|f(z) - f(w)| : d(z,w) < 1\},
$$
\noindent
which is another essential ingredient in the proof of Theorem 1.1(i).    We show that 
every operator in  EssCom$\{T_f : f \in {\text{VO}}_{\text{bdd}}\}$  satisfies an  ``$\epsilon $-$\delta $"  
condition  involving  ``diff".
\medskip
In Section 8 we construct approximate partitions of the unity on  $\Omega $   that satisfy two competing 
requirements:  (1)  The ``diff" for the partition functions must be small.  (2)  There is a  fixed, finite cap on 
the overlaps of the sets involved.   This construction is based on a suitable   
analogue of   ``radial-spherical decomposition"  for  $\Omega $.  
As it turns out,   the gauge    $\rho (z,w)$    plays the role of  ``spherical coordinates"  in our decomposition,  whereas the defining 
function  $r$  gives us a convenient   ``radial coordinate".
\medskip
With all the above preparation,  we prove Theorem 1.1(i) in Section 9.  The gist of the proof is that the 
``$\epsilon $-$\delta $"   condition mentioned above characterizes the membership   $X \in {\Cal T}$.  
The same work also shows that for    $A \in {\Cal T}$,  if  
$$
\lim _{z \rightarrow \partial \Omega }\sup \{|\langle Ak_w,k_z\rangle |  :  d(z,w) < R\}  =  0
\tag 1.2
$$
\noindent
for every given  $0 < R < \infty $,   then $A$  is a compact operator.   This is a major step in the proof of Theorem 1.2.   
In fact, what remains for the proof of Theorem 1.2 is to show that (1.1) implies (1.2).
\medskip
Then in Section 10, we turn to the proof of part (ii) in Theorem 1.1.   With the work in Section 9, this is now relatively easy.  
First of all, Theorem 1.1(i) tells us that    $\text{EssCom}({\Cal T})$   coincides with the essential center of   ${\Cal T}$.   That is,  
$\text{EssCom}({\Cal T}) \subset {\Cal T}$.   Then we show that the membership   $A \in \text{EssCom}({\Cal T})$
implies that the  Berezin transform   $\tilde A$  of  $A$  is in  ${\text{VO}}_{\text{bdd}}$.   Since  $A  - T_{\tilde A} \in {\Cal T}$,   
the membership  $\tilde A \in {\text{VO}}_{\text{bdd}}$   and the work in Section 9  
lead to an easy proof of the fact that  $A  - T_{\tilde A} \in {\Cal K}$,  which proves Theorem 1.1(ii).
\medskip
Finally, in Section 11 we show that (1.1) indeed implies (1.2).      
For all previous works involving this step, this was easy, because one could use M\"obius transforms.   
But in our case of a general strongly pseudo-convex domain,  this becomes a non-trivial step.  
Material from Sections 2-4 will be needed for this step. 

\bigskip
\centerline{{\bf 2. A metric on} $\Omega $  {\bf and related facts}}

\medskip
First of all,  we cite [15,20] as general references for strongly pseudo-convex domains.  
Throughout the paper,   $\Omega $ denotes a bounded, connected,
strongly pseudo-convex domain in  ${\bold C}^n$  with smooth boundary.    
More precisely,  we always assume that    $\Omega $  is bounded and connected,  
and that there is a real-valued  $C^\infty $  function   $r$  defined in an open neighborhood of  the closure of   
$\Omega $  such that the following three conditions are satisfied:

(1)  $\Omega = \{z \in {\bold C}^n : r(z) < 0\}$.

(2)  $|(\nabla r)(z)|  \neq 0$   for every   $z \in \partial \Omega $.

(3)  There is a    $c > 0$  such that

$$
\sum _{i,j=1}^n(\partial _i\bar \partial _jr)(z)\xi _i\bar \xi _j \geq  c(|\xi _1|^2 + \cdots + |\xi _n|^2)
\tag 2.1
$$

for all   $z \in \partial \Omega $   and  $\xi _1,\dots ,\xi _n \in {\bold C}$.

\noindent
Such an  $r$  is called a defining function for the domain,  and will be fixed along with  $\Omega $.

\medskip
It will be convenient to adopt the following convention:
We will consider  ${\bold C}^n$  as a column space whenever an  $n\times n$  matrix acts on it.
When there is no matrix involved, we will consider  ${\bold C}^n$  either as a column space or as a row space,  
whichever is more appropriate.
\medskip
Let  $A(z)$  be the  $n\times n$  matrix whose entry in the intersection of $i$-th column and  $j$-row is 
$(\partial _i\bar \partial _jr)(z)$,  $i, j = 1, \dots, n$.    By (2) and (3),   there is a  $\theta > 0$  such that
if   $w \in \Omega $  and  $r(w) > -3\theta $,    then   $|(\nabla r)(w)| \neq 0$  and  
$$
\langle A(w)\xi ,\xi \rangle \geq  (c/2)|\xi |^2
\tag 2.2
$$
\noindent
for all    $\xi \in {\bold C}^n$.   Let   $\psi  : {\bold R} \rightarrow  [0,1]$
be a  $C^\infty $  function  such that  $\psi = 1$  on   $[-\theta ,\infty )$    and  $\psi = 0$   on   $(-\infty ,-2\theta ]$.  
Write  $\delta _{ij}$  for  Kronecker's delta.  We then define
$$
b_{ij}(z)  =  \psi (r(z))\bigg({1\over -r(z)}(\partial _i\bar \partial _jr)(z)  +  {1\over r^2(z)}(\partial _ir)(z)(\bar \partial _jr)(z)\bigg)  +  (1 - \psi (r(z)))\delta _{ij}
\tag 2.3
$$
\noindent
for  $i,j \in \{1,\dots ,n\}$   and   $z \in \Omega $.    Let  ${\Cal B}(z)$  be the  $n\times n$  matrix whose 
entry in the intersection of $i$-th column and  $j$-row is  $b_{ij}(z)$,  $i, j = 1, \dots, n$.   
From  (2.2) and the definition of   $\psi $  we see that the    ${\Cal B}(z)$  is invertible for every   $z \in \Omega $.   
Thus the local Hermitian form
$$
H_z(\xi ,\eta )   =  \langle {\Cal B}(z)\xi ,\eta \rangle,    \quad  \xi ,\eta \in  T_z\Omega  = {\bold C}^n,
$$
\noindent
generates a non-degenerate metric   $d$   on  $\Omega $.   That is,  for  $z, w \in \Omega $,
$$
d(z,w)  =  \inf \int _0^1\sqrt{\langle {\Cal B}(g(t))g'(t),g'(t)\rangle }dt,
\tag 2.4
$$
\noindent
where the infimum is taken over all  $C^1$   maps  $g : [0,1] \rightarrow \Omega $   satisfying the conditions  
$g(0)$  $=$  $z$   and   $g(1) = w$.   The definition of   $\psi $  ensures that for  $i,j \in \{1,\dots ,n\}$,
$$
b_{ij}(z)  =    \partial _i\bar \partial _j\log {1\over -r(z)}   \quad  \text{whenever}  \ \  -\theta \leq r(z) < 0.
$$
\noindent
Denote   $\bar \partial = (\bar \partial _1,\dots ,\bar \partial _n)$,  which will play a prominent role throughout the paper.  Then  
$$
\langle {\Cal B}(z)\xi ,\xi \rangle =  {\langle A(z)\xi ,\xi \rangle \over -r(z)} +  \bigg({|\langle \xi ,(\bar \partial r)(z)\rangle |\over -r(z)}\bigg)^2
\quad  \text{whenever}  \ \  -\theta \leq r(z) < 0,
\tag 2.5
$$
\noindent
$\xi \in {\bold C}^n$.   These identities make  $d$   an {\it imitation} of the Bergman metric on  $\Omega $.   Compared with the real  
Bergman metric  $\beta $,  our imitation  $d$  has the advantage that 
the explicit formulas above will greatly simplify many of the estimates below.

\medskip
\noindent
{\bf Lemma 2.1.}   {\it There is a}  $c_{2.1} > 0$  {\it  such that}   $-r(w) \geq - c_{2.1}2^{-4d(z,w)}r(z)$  {\it  for all}  
$z, w \in \Omega $.
\medskip
\noindent
{\it Proof}.  Consider any $z, w \in \Omega $  such that   $-r(z) \leq \theta $  and
$$
-r(w) \leq -(1/2)r(z).
$$
\noindent   
Let    $g : [a,b] \rightarrow \Omega $  be a  $C^1$  map such that  $g(a) = z$   and  $g(b) = w$.    
Let  $a'$  be the largest number in  $[a,b]$  such that    $r(g(a')) = r(z)$.   Then
$-r(g(t))  \leq  -r(g(a'))$   for    $t \in [a',b]$  and
$$
\int _{a'}^b{d\over dt}r(g(t))dt  =  r(g(b)) - r(g(a'))  \geq  -(1/2)r(z).
$$
\noindent
Note that
$$
{d\over dt}r(g(t)) =  2\text{Re}\langle g'(t),(\bar \partial r)(g(t))\rangle .
$$
\noindent
Therefore the above implies
$$
2\int _{a'}^b|\langle g'(t),(\bar \partial r)(g(t))\rangle |dt  \geq  \int _{a'}^b{d\over dt}r(g(t))dt  \geq  -(1/2)r(z).
$$
\noindent
Since   $-r(g(t)) \leq - r(z) \leq \theta $  for every  $t \in [a',b]$,    by (2.5)  we have
$$
\int _{a'}^b\sqrt{\langle {\Cal B}(g(t))g'(t),g'(t)\rangle }dt  \geq
\int _{a'}^b{|\langle g'(t),(\bar \partial r)(g(t))\rangle |\over -r(g(t))}dt   
\geq  \int _{a'}^b{|\langle g'(t),(\bar \partial r)(g(t))\rangle |\over -r(z)}dt   \geq  {1\over 4}.
$$
\noindent
From the above we see that  if  $z, w \in \Omega $   satisfy the conditions   $-r(z) \leq \theta $  and
$$
-r(w) \leq -2^{-m}r(z)  \quad  \text{for some} \ \  m \in {\bold N},
$$
\noindent
then for any   $C^1$   map  $g : [0,1] \rightarrow \Omega $ with the properties   $g(0) = z$   and  $g(1) = w$  we have
$$
\int _0^1\sqrt{\langle {\Cal B}(g(t))g'(t),g'(t)\rangle }dt    
\geq  {m\over 4}.
$$
\noindent
Combining this with   (2.4),  we have  $m \leq 4d(z,w)$.    This implies the inequality 
$$
-r(w) > -(2^{-4d(z,w)-1})r(z)
$$
\noindent
for all  $w \in \Omega $  whenever  $-r(z) \leq \theta $.
\medskip
Suppose that  $-r(z) > \theta $.  If  $-r(w) \geq \theta $,  then the case is trivial,  as the function  $-r$  has a maximum on  $\Omega $.  
Suppose that  $-r(w) < \theta $.     There is a   $C^1$  map  
$g : [0,1] \rightarrow \Omega $ such that  $g(0) = z$,  $g(1) = w$   and
$$
\int _0^1\sqrt{\langle {\Cal B}(g(t))g'(t),g'(t)\rangle }dt  \leq  d(z,w) + (1/4). 
$$  
We have  $-r(g(0)) > \theta $   and  $-r(g(1)) < \theta $.   Thus there is an  $a \in  [0,1]$  such that  
$-r(g(a)) = \theta $.     Define   $z' = g(a)$.   By what we proved above,  $-r(w) > 2^{-4d(z',w)-1}\theta $.  Since   $a \in  [0,1]$,  $z' = g(a)$    
and  $g(1) = w$,    we have    $d(z',w) \leq d(z,w) + (1/4)$.   Hence
$$
-r(w)\geq  2^{-4\{d(z,w)+(1/4)\}-1}\theta = -\{\theta /(-r(z))\}(2^{-4d(z,w)-2})r(z).
$$  
Since   $-r$  has a maximum  on  $\Omega $,  the lemma also holds in the case  $-r(z) > \theta $.     $\square $
\medskip
\noindent
{\bf Lemma 2.2.}   {\it  There is a constant}   $0 < C_{2.2} < \infty $   {\it  such that}
$$
|z - w|^2 + |\langle z - w,(\bar \partial r)(z)\rangle |  \leq  C_{2.2}\{d(z,w) +d^2(z,w)\}2^{12d(z,w)}(-r(z))
$$
\noindent
{\it  for all}  $z, w \in \Omega $.
\medskip
\noindent
{\it Proof}.  We first show that there is a  $C$  such that
$$
|z - w|  \leq Cd(z,w)2^{4d(z,w)}\sqrt{-r(z)}
\tag 2.6
$$
\noindent
for all  $z, w \in \Omega $.   By (2.4),  for  any given    $z, w \in \Omega $,  there is a  $C^1$  map  $g : [0,1] \rightarrow \Omega $
such that   $g(0) = z$,  $g(1) = w$,  and  
$$
\int _0^1\sqrt{\langle {\Cal B}(g(t))g'(t),g'(t)\rangle }dt  \leq  2d(z,w).
$$
\noindent
There is a  $t_0 \in [0,1]$  such that  $-r(g(t_0)) \geq -r(g(t))$   for every  $t \in [0,1]$.  By (2.5) and (2.3),  
there is a  $c_1 > 0$  such that
$$
{c_1\over \sqrt{-r(g(t_0))}}\int _0^1|g'(t)|dt  \leq \int _0^1\sqrt{\langle {\Cal B}(g(t))g'(t),g'(t)\rangle }dt  \leq  2d(z,w).
$$
\noindent
Set  $C_1 = 2/c_1$.   Since    $g(0) = z$  and  $g(1) = w$,  the above implies
$$
|z - w|  \leq  C_1d(z,w) \sqrt{-r(g(t_0))}.
\tag 2.7
$$
\noindent
If we write  $\zeta = g(t_0)$,   then  
$$
d(z,\zeta )  \leq  \int _0^{t_0}\sqrt{\langle {\Cal B}(g(t))g'(t),g'(t)\rangle }dt  \leq  2d(z,w).
$$
\noindent
By Lemma 2.1,  we have  $-r(\zeta ) \leq c_{2.1}^{-1}2^{4d(z,\zeta )}(-r(z)) \leq  c_{2.1}^{-1}2^{8d(z,w)}(-r(z))$.  Combining this with  
(2.7),  (2.6)  follows.
\medskip
The same argument also shows that  $d(z,g(t)) \leq 2d(z,w)$  for every
$t \in [0,1]$.  Therefore 
$$
|z - g(t)|  \leq 2Cd(z,w)2^{8d(z,w)}\sqrt{-r(z)},
$$
\noindent
$t \in [0,1]$.     Using Lemma 2.1 and the obvious Lipschitz condition for  $\bar \partial r$,  we have
$$
\align
|\langle z -w,&(\bar \partial r)(z)\rangle |  =     |\langle g(1) - g(0),(\bar \partial r)(z)\rangle |   
\leq  \int _0^1|\langle g'(t),(\bar \partial r)(z)\rangle |dt  \\
&\leq  \int _0^1|\langle g'(t),(\bar \partial r)(g(t))\rangle |dt  +  \int _0^1|g'(t)||(\bar \partial r)(z) - (\bar \partial r)(g(t))|dt  \\
&\leq  {2^{8d(z,w)}\over c_{2.1}}(-r(z))\int _0^1{|\langle g'(t),(\bar \partial r)(g(t))\rangle |\over -r(g(t))}dt  \\
&\quad  \quad + 2Cd(z,w)2^{8d(z,w)}\cdot C_2 \bigg({2^{8d(z,w)}\over c_{2.1}}\bigg)^{1/2}(-r(z))\int _0^1\sqrt{\langle {\Cal B}(g(t))g'(t),g'(t)\rangle }dt  \\
&\leq  C_3(-r(z))\{2^{8d(z,w)}d(z,w) + 2^{12d(z,w)}d^2(z,w)\}.
\endalign
$$
\noindent
Combining this with (2.6),  the lemma is proved.   $\square $
\medskip
For  $z \in \Omega $   and   $a > 0$,  define the imitation Bergman metric ball
$$
D(z,a) = \{w \in \Omega : d(z,w) <  a\}.
$$
\noindent
{\bf Definition  2.3.}   For    $\eta \in {\bold C}^n\backslash \{0\}$,   $a > 0$ and  $b > 0$,  we let  
${\Cal P}(\eta ;a,b)$   be the collection of  vectors   $u + v$  satisfying the following three conditions:

(1)    $u, v \in {\bold C}^n$   with   $|u| < a$   and  $|v| < b$.

(2)    $u \perp \eta $.

(3)    $v \in \{\xi \eta : \xi \in {\bold C}\}$.

\medskip
\noindent
{\bf Proposition 2.4.}   {\it  Given any}  $0 < a < \infty $,   {\it  there are}   $0 < c \leq C < \infty $   {\it  such that}
$$
z + {\Cal P}((\bar \partial r)(z);c\sqrt{-r(z)},-cr(z))  \subset  D(z,a)  \subset  z + {\Cal P}((\bar \partial r)(z);C\sqrt{-r(z)},-Cr(z)) 
$$
\noindent
{\it  for every}   $z \in \Omega $  {\it  satisfying the condition}  $-r(z) < \theta $.

\medskip
\noindent
{\it Proof}.  Let  $0 < a < \infty $  be given and consider a sufficiently small   $c > 0$.  
Let  $u, v \in {\bold C}^n$  satisfy the conditions  
$u \perp (\bar \partial r)(z)$,    $|u| < c\sqrt{-r(z)}$,  $v \in \{\xi (\bar \partial r)(z) : \xi \in {\bold C}\}$,   and
$|v| < -cr(z)$.   We want to show that   $d(z,z+u+v) < a$.    To prove this,  consider the path
$$
g(t)  = z + t(u+v),
$$
\noindent
$t \in [0,1]$.    Then   $g'(t) = u + v$.    By the Taylor expansion for  $\bar \partial r$,   we have
$$
\align
\langle g'(t),(\bar \partial r)(g(t))\rangle   &= \langle u+v,(\bar \partial r)(z + t(u+v))\rangle  \\
&=   \langle u+v,(\bar \partial r)(z) + tX(z)(u+v) + o(|u+v|)\rangle   \\
&=  \langle v,(\bar \partial r)(z)\rangle + \langle u + v,tX(z)(u+v) + o(|u+v|)\rangle ,
\endalign
$$
\noindent
where  $X(z)$  is the derivative of   $\bar \partial r$  at   $z$,  which is a linear map from  ${\bold C}^n$  to  ${\bold C}^n$.  
Since  $|v| < -cr(z)$  and  $|u + v| \leq c(-r(z) + \sqrt{-r(z)})$,  we see that
$$
|\langle g'(t),(\bar \partial r)(g(t))\rangle |  \leq  cM(-r(z))
$$
\noindent  
for   $t \in [0,1]$.    Since  $r$  is real-valued,  Taylor expansion gives us
$$
\align
r(g(t)) = r(z + t(u+v))  &=  r(z)  +  2t\text{Re}\langle u+v,(\bar \partial r)(z)\rangle +  O(|u+v|^2)  \\
&=  r(z)  +  2t\text{Re}\langle v,(\bar \partial r)(z)\rangle +  O(-c^2r(z)).
\endalign
$$
Since    $|v| < -cr(z)$  and since  $c$  is small,  we obtain
$$
{|\langle g'(t),(\bar \partial r)(g(t))\rangle | \over -r(g(t)) }  \leq  cM'.
$$
\noindent
Similarly,  we have
$$
{\langle A(g(t))g'(t),g'(t)\rangle \over -r(g(t))}  =  {\langle A(g(t))(u+v),u+v\rangle \over -r(g(t))}   \leq  c^2M''.
$$
\noindent
Combining these two inequalities with (2.5) and (2.4),  we see that the smallness of  $c$  ensures  $d(z,z+u+v) < a$.   
This proves the first inclusion in the proposition.   

\medskip
It is easy to see that the second inclusion,   $D(z,a) \subset \cdots $,  is simply a consequence of Lemma 2.2.
This completes the proof.     $\square $

\medskip
\noindent
{\bf Proposition 2.5.}   {\it  There is a}    $0 < C_{2.5} < \infty $   {\it  such that if} $0 < a < 1/2$,  {\it  then}
$$
D(z,a)  \subset  z + {\Cal P}((\bar \partial r)(z);C_{2.5}a\sqrt{-r(z)},-C_{2.5}ar(z)) 
$$
\noindent
{\it  for every}   $z \in \Omega $   {\it  satisfying the condition}   $-r(z) < \theta $.

\medskip
\noindent
{\it Proof}.    Suppose that   $0 < a < 1/2$   and that    $z, w \in \Omega $  satisfy the condition  $d(z,w) < a$.
Then there is a  $C^1$  map   $g : [0,1] \rightarrow \Omega $   with   $g(0) = z$   and  $g(1) = w$  such that
$$
\int _0^1\sqrt{\langle {\Cal B}(g(t))g'(t),g'(t)\rangle }dt <  2a.
$$
\noindent
Thus   $d(z,g(t)) <  2a < 1$   for every  $t \in [0,1]$.   Suppose that    $w = z + u + v$   with   $u \perp (\bar \partial r)(z)$   and  
$v \in \{\xi (\bar \partial r)(z) : \xi \in {\bold C}\}$.     To estimate  $|u|$,  we again apply
Lemma 2.1,  which gives us  $-r(z) \geq -c(1)r(g(t))$  for every   $t \in [0,1]$.  Since  $|u| \leq |z - w|$,  we have  
$$
|u| \leq \int _0^1|g'(t)|dt  \leq  C_3\sqrt{{-r(z)\over c(1)}}\int _0^1\sqrt{\langle {\Cal B}(g(t))g'(t),g'(t)\rangle }dt 
<  C_3\sqrt{{-r(z)\over c(1)}}\cdot 2a.
\tag 2.8
$$
\noindent
To estimate  $|v|$,  we  apply Lemma 2.2.  Since  $0 < a < 1/2$,   Lemma 2.2 gives us
$$
|\langle v,(\bar \partial r)(z)\rangle |  \leq  C_{2.2}d(z,w)(3/2)2^6(-r(z)).
\tag 2.9
$$
\noindent
Recall that  $\theta $  was chosen so that   $(\bar \partial r)(\zeta ) \neq 0$   whenever    $0 < -r(\zeta ) < 3\theta $.  
Hence the proposition follows from (2.8) and (2.9).     $\square $

\medskip
On the domain  $\Omega $  we define the measure  
$$
d\mu (z)  =   {dv(z)\over (-r(z))^{n+1}}.
\tag 2.10
$$ 
\medskip
\noindent
{\bf  Proposition 2.6.}     {\it  For each}  $a \in (0,\infty )$,   {\it  there are}   $0 < c(a) \leq C(a) < \infty $  {\it  such that}
$$
c(a)  \leq  \mu (D(z,a))  \leq  C(a)
$$
\noindent
{\it  for every}   $z \in \Omega $.

\medskip
\noindent
{\it Proof}.  Since  $v({\Cal P}(\eta ;x,y)) = C_nx^{2n-2}y^2$,  this follows immediately from Proposition 2.4 and Lemma 2.1.  $\square $

\medskip
For each   $0 \leq  \rho < \theta $,  define the surface
$$
S_\rho  = \{ z \in {\bold C}^n :  -r(z) = \rho \}.
$$
\noindent
In particular,  we have  $S_0 = \partial \Omega $,  the boundary of the domain   $\Omega $.
\medskip
\noindent
{\bf Proposition 2.7.}   {\it  There exist a finite  open cover}  $U_1, \dots ,U_m$  of  
$$
H  =  \{z \in {\bold C}^n :  0 \leq  -r(z) \leq \theta /2\}
$$
\noindent
{\it  in}  ${\bold C}^n$   {\it  and a}  $1 \leq  C < \infty $   {\it  such that the following holds true}:   
{\it  Suppose that}  $0 < \rho \leq \theta /2$  {\it  and that}  $z, w \in S_\rho \cap U_i$  {\it  for some}  $i \in \{1,\dots ,m\}$. 
{\it  Furthermore},  {\it  suppose that there is an}   $R \geq 1$  {\it  such that}  $|z - w| \leq R\sqrt{\rho }$   {\it  and}     
$|\langle z - w,(\bar \partial r)(z)\rangle | \leq  R^2\rho $.   {\it  Then}  $d(z,w) \leq CR^2$.
\medskip
\noindent
{\it Proof}.    For  $\zeta \in {\bold C}^n$  and  $a > 0$,  denote  $B(\zeta ,a) = \{\xi \in {\bold C}^n : |\zeta - \xi | < a\}$   as usual.  
Note that by assumption,  $H$  is a compact set on which  $|\nabla r|$  does not vanish.  
By the usual open covering argument,  there is a   $\tau  > 0$  such that if  $z_0 \in  H$,  
then the conclusion of the standard implicit function theorem holds on   $B(z_0,\tau )$   for the equation  $r = r(z_0)$.   
See, e.g.,  [17,page 74].    
Since  $H$  is compact,  there are  $z_1, \dots, z_m \in H$   such that   $\cup _{i=1}^mB(z_i,\tau /2) \supset H$.    
We define  $U_i = B(z_i,\tau /2)$,  $i = 1, \dots ,m$.

\medskip
Now let   $0 < \rho \leq \theta /2$,   and let  $z, w \in S_\rho \cap U_i$  satisfy the conditions   
$|z - w| \leq R\sqrt{\rho }$   and  $|\langle z - w,(\bar \partial r)(z)\rangle | \leq  R^2\rho $   for some   $R \geq 1$.   
Then, of course,  $|z - w| < \tau $.    
By the discussion in the first paragraph,  every point in  $S_\rho \cap B(z,\tau )$   can be expressed in the form   
$$
z + x + f_z(x),
$$
\noindent
where   $x \in {\bold C}^n$   satisfies the conditions   $\text{Re}\langle x,(\bar \partial r)(z)\rangle = 0$   and   $|x| < \tau '$,
and where  $f_z$   satisfies the condition  $|f_z(x)| \leq  C_1|x|^2$.   
Since  the implicit function theorem provides bounds that are independent of the points in  $H$,
reducing the value of  $\tau $ if necessary,  we may assume that
$|f_z(x)| \leq (1/2)|x|$   when   $|x + f_z(x)| < \tau $.  Let   $x_0 \in B(0,\tau ')$  be such that    $w =$  $z + x_0 + f_z(x_0)$.  
Then   $(1/2)|x_0| \leq  |x_0 + f_z(x_0)| =  |z - w|$.   Hence  $|x_0| \leq  2R\sqrt{\rho }$.  We have
$$
\align
|\langle x_0,(\bar \partial r)(z)\rangle | &\leq   |\langle w - z,(\bar \partial r)(z)\rangle | + |\langle f_z(x_0),(\bar \partial r)(z)\rangle |   
\leq  R^2\rho + C_2|f_z(x_0)|    \\
&\leq    R^2\rho + C_3|x_0|^2   \leq  R^2\rho + 4C_3R^2\rho  =  C_4R^2\rho .
\endalign
$$
\noindent
Now define the map   $g : [0,1] \rightarrow S_\rho $  by the formula  
$$
g(t)  =  z + tx_0  + f_z(tx_0),
$$
\noindent
$t \in [0,1]$.   We have   $g'(t) = x_0 + (Df_z)(tx_0)x_0$,   where   $Df_z$  is the derivative of  $f_z$.  
Recalling (2.5),  we have
$$
{1\over \rho}\langle A(g(t))g'(t),g'(t)\rangle \leq  {C_5\over \rho }|x_0|^2 \leq {C_5\over \rho }4R^2\rho = C_6R^2.
$$
\noindent
Therefore
$$
\int _0^1\bigg({\langle A(g(t))g'(t),g'(t)\rangle\over - r(g(t))}\bigg)^{1/2}dt  \leq  \sqrt{C_6}R = C_7R.
\tag 2.11
$$
\noindent
Note that the condition  $|f_z(x)| \leq  C_1|x|^2$  implies that   $(Df_z)(0) = 0$.   Hence
$$
g'(t)  =    x_0 + \{(Df_z)(tx_0) - (Df_z)(0)\}x_0  =  x_0 + h(t)
$$
\noindent
with   $|h(t)|  \leq  C_8|x_0|^2  \leq  4C_8R^2\rho =  C_9R^2\rho $.    Consequently    
$$
\align
|\langle g'(t),(\bar \partial r)(g(t))\rangle |   &\leq  |\langle g'(t),(\bar \partial r)(z)\rangle |  +  |\langle g'(t),(\bar \partial r)(g(t)) - (\bar \partial r)(z)\rangle |  \\
&\leq  |\langle x_0,(\bar \partial r)(z)\rangle | + |\langle h(t),(\bar \partial r)(z)\rangle |  +  C_{11}|g'(t)||g(t) -z|  \\
&\leq  C_4R^2\rho  + C_{10}C_9R^2\rho  + C_{12}|x_0|^2  \leq   C_{13}R^2\rho .
\endalign
$$
\noindent
Thus we have
$$
\int _0^1{|\langle g'(t),(\bar \partial r)(g(t))\rangle | \over - r(g(t))}dt  \leq   C_{13}R^2.
$$
\noindent
Recalling (2.5), (2.4)  and combining the above with (2.11),  we find that  $d(z,w) \leq  C_7R + C_{13}R^2$.   Since we assume  
$R \geq 1$,  it follows that    $d(z,w) \leq  (C_7 + C_{11})R^2$.   $\square $
\medskip
For each  $0 \leq \rho \leq \theta /2$,   we write  $d\sigma _\rho $ for the natural surface measure 
on    $S_\rho $.   For every  triple of    $0 \leq \rho \leq \theta /2$,    $\zeta \in S_\rho $  and    $t > 0$,  we define 
$$
Q_\rho (\zeta ,t)  =  \{\xi \in S_\rho :  |\zeta - \xi |^2 + |\langle \zeta - \xi ,(\bar \partial r)(\zeta )\rangle | < t\}.
$$
\medskip
\noindent
{\bf Proposition 2.8.}   {\it There are constants}  $0 < \tau  \leq \theta /2$  {\it and}  $0 < c_{2.8} \leq  C_{2.8} < \infty $  {\it such that}
$$
c_{2.8}t^n  \leq  \sigma _\rho (Q_\rho (\zeta ,t)) \leq  C_{2.8}t^n
\tag 2.12
$$
\noindent
{\it for all}    $0 \leq \rho \leq  \tau $,     $\zeta \in  S_\rho $  {\it and}   $0 < t \leq T_0$,  {\it where}  
$T_0 = \sup \{|u - v|^2 +  |\langle u - v,(\bar \partial r)(u)\rangle | : u, v \in \Omega \}$.
\medskip
\noindent
{\it Proof}.   First, we remark that the restriction  $t \leq T_0$  is only necessary to guarantee the {\it lower bound} in (2.12).  
Second,  adjusting the constants  $c_{2.8}$  and  $C_{2.8}$  if necessary,   it suffices to find an  $a \in (0,T_0]$  such that  
(2.12) holds for all  $0 < t  <  a$.
\medskip
For each   $\zeta \in S_\rho $,  $0 \leq \rho \leq \theta /2$,   denote    $T_\zeta = \{x \in {\bold C}^n : \text{Re}\langle x,(\bar \partial r)(\zeta )\rangle = 0\}$,  
which is the real tangent space to  $S_\rho $   at  $\zeta $.  
For each pair of  $\zeta \in S_\rho $   and   $s > 0$,  define   
$E_\zeta (s) = \{x \in T_\zeta : |x|^2 + |\langle x,(\bar \partial r)(\zeta )\rangle | < s\}$.  
Each  $x \in T_\zeta $  has the decomposition  $x = y + z$,  where  $\langle y,(\bar \partial r)(\zeta )\rangle = 0$   
and   $z \in \{w(\bar \partial r)(\zeta ) : w \in {\bold C}\}$.   Since   $\text{Re}\langle x,(\bar \partial r)(\zeta )\rangle = 0$,  we have  
$z = ih(\bar \partial r)(\zeta )$   for some  $h \in $  {\bf R}.   Let  $v_{2n-1}$  denote the real $(2n-1)$-dimensional  volume measure on  
$T_\zeta $.   Using this  $x = y +z$   decomposition,  it is elementary that there are   $0 < c_1 \leq C_1 < \infty $  such that
$$
c_1s^n  \leq  v_{2n-1}(E_\zeta (s)) \leq C_1s^n
\tag 2.13
$$
\noindent
for all  $\zeta \in S_\rho $,  $0 \leq \rho \leq \theta /2$,  and  $s > 0$.
\medskip
As in the proof of Proposition 2.7, we apply the standard implicit function theorem.   
There is a    $0 <  \tau  \leq \theta /2$  such that the conclusion of the implicit function theorem
holds on  $\{z \in {\bold C}^n : 0 \leq -r(z) \leq  \tau \}$   for  $r$  with uniform bounds.   Namely,  
there are constants  $b > 0$,  $0 < c \leq 1$  and  $0 < C_2 < \infty $  such that if  $0 \leq \rho \leq  \tau $  and  
$\zeta \in S_\rho $,  then every element in   $B(\zeta ,b)\cap S_\rho $  can be expressed in the form 
$$
\Phi _{\zeta }(x)  = \zeta + x + f_\zeta (x) 
$$
\noindent
for some  $x \in T_\zeta\cap B(0,c)$,   where   
$f_\zeta $  satisfies the conditions    $|f_\zeta (x)| \leq C_2|x|^2$    and   $|f_\zeta (x)| \leq (1/2)|x|$  when  $|x| < c$.  
Furthermore,  there are constants $0 < c_3 \leq C_3 < \infty $  such that the matrix inequality   $c_3 \leq (D\Phi _\zeta (x))^\ast D\Phi _\zeta (x) \leq C_3$  
holds whenever   $0 \leq - r(\zeta ) \leq \tau $  and    $x \in T_\zeta\cap B(0,c)$,  where    $D\Phi _\zeta $   is the derivative of 
$\Phi _\zeta $,  which is a  $2n\times (2n-1)$  real matrix.   
\medskip
Now take  $a_1 =  b^2$.   Let    $0 \leq \rho \leq  \tau $  and  $0 < t < a_1$.   
Then for any pair of  $\zeta \in S_\rho $  and  $\xi \in Q_\rho (\zeta ,t)$,  we have  $\xi \in B(\zeta ,b)$.  Therefore we can write
$$
\xi = \zeta + x + f_\zeta (x)
$$ 
\noindent
for some   $x \in T_\zeta\cap B(0,c)$.  We have  $|\xi - \zeta | = |x + f_\zeta (x)| \geq (1/2)|x|$   and  $\langle x,(\bar \partial r)(\zeta )\rangle $  
$=$  $\langle \xi - \zeta ,(\bar \partial r)(\zeta )\rangle - \langle f_\zeta (x),(\bar \partial r)(\zeta )\rangle $  with  $|f_\zeta (x)| \leq C_2|x|^2$.  
Hence
$$
|x|^2 + |\langle x,(\bar \partial r)(\zeta )\rangle |  \leq  C_4\{|\xi - \zeta |^2 +  |\langle \xi - \zeta ,(\bar \partial r)(\zeta )\rangle |\}
\tag 2.14
$$
\noindent
for some constant  $1 \leq C_4 < \infty $.  Similarly,  $|\xi - \zeta | = |x + f_\zeta (x)| \leq  (1+C_3)|x|$  and  $\langle \xi - \zeta ,(\bar \partial r)(\zeta )\rangle $  
$=$   $\langle x,(\bar \partial r)(\zeta )\rangle +  \langle f_\zeta (x),(\bar \partial r)(\zeta )\rangle $  with  $|f_\zeta (x)| \leq C_2|x|^2$.
Consequently
$$
|\xi - \zeta |^2 +  |\langle \xi - \zeta ,(\bar \partial r)(\zeta )\rangle |  \leq  C_5\{|x|^2 + |\langle x,(\bar \partial r)(\zeta )\rangle |\}
\tag 2.15
$$
\noindent
for some constant  $1 \leq C_5 < \infty $.  Set  $a_2 = c^2/C_4$  and  $a = \min \{a_1,a_2,T_0\}$.  
If  $0 < t < a$,   then   $E_\zeta (C_4t) \subset T_\zeta\cap B(0,c)$.  
Thus   (2.14)  implies that  for  $0 < t < a$,  we have
$$
\Phi _\zeta (E_\zeta (C_4t))  \supset Q_\rho (\zeta ,t).
$$
\noindent
Combining  the smoothness of   $\Phi _\zeta $  on  
$T_\zeta\cap B(0,c)$   with the upper bound in (2.13),  we obtain
$$
\sigma _\rho (Q_\rho (\zeta ,t))  \leq   \sigma _\rho (\Phi _\zeta (E_\zeta (C_4t)) )
\leq C_6v_{2n-2}(E_\zeta (C_4t))  \leq   C_6C_1(C_4t)^n  =  C_7t^n,
$$
\noindent
which gives us the upper bound in (2.12).
Similarly, since  $a < c^2$  and  $C_5 \geq 1$,  for   $0 < t < a$   we have   $E_\zeta (t/C_5)  \subset T_\zeta\cap B(0,c)$.    Therefore for  $0 < t < a$,
(2.15)  implies  
$$
\Phi _\zeta (E_\zeta (t/C_5))  \subset  Q_\rho (\zeta ,t).
$$
\noindent
From  the non-singularity of   $\Phi _\zeta $  on   $T_\zeta\cap B(0,c)$  and  the lower bound in (2.13)  we obtain
$$
c_1(t/C_5)^n \leq  v_{2n-1}(E_\zeta (t/C_5))  \leq  C_8\sigma _\rho (\Phi _\zeta (E_\zeta (t/C_5)))  \leq  C_8\sigma _\rho (Q_\rho (\zeta ,t)),
$$
\noindent
proving the lower bound in (2.12).   This completes the proof.   $\square $
\medskip
\noindent
{\bf Proposition 2.9.}   {\it  There is a constant}   $0 < C_{2.9} < \infty $    {\it  such that the following holds true}:   {\it  Let}  
$z \in \Omega $,   $k \in {\bold Z}$   {\it  and}   $j \in {\bold Z}_+$.   {\it  Then the volume of the set}
$$
\align
W_{z;k,j}  =   \{w \in \Omega : \  &2^{k-1}(-r(z))   < - r(w) \leq 2^k(-r(z)) \\
&\text{and} \ \  |z - w|^2 + |\langle z - w,(\bar \partial r)(z)\rangle | \leq 2^{k+j}(-r(z))\} 
\endalign
$$
\noindent
{\it  does not exceed}   $C_{2.9}2^{nj}(-2^kr(z))^{n+1}$.
\medskip
\noindent
{\it Proof}.   First of all,  there is a   $C_1$  such that   $|(\bar \partial r)(\zeta ) - (\bar \partial r)(\zeta ')|  \leq  C_1|\zeta - \zeta '|$    
for all  $\zeta ,\zeta ' \in \Omega $.     Suppose that  $2^k(-r(z)) \leq \tau $,  where   $\tau $  is the same as in Proposition 2.8.
For any value  $2^{k-1}(-r(z))   < \rho  \leq 2^k(-r(z))$,  denote
$$
\Sigma (\rho ;k,j) =  \{w \in S_\rho :   |z - w|^2 + |\langle z - w,(\bar \partial r)(z)\rangle | \leq 2^{k+j}(-r(z))\}.
$$
\noindent
Suppose that   $\Sigma (\rho ;k,j)  \neq \emptyset $,  and  pick a   $w_\rho \in  \Sigma (\rho ;k,j)$.  Then elementary estimates yield
$$
\Sigma (\rho ;k,j)  \subset  Q_\rho (w_\rho ,C_22^{k+j}(-r(z))).
$$
\noindent
Applying Proposition 2.8,   we obtain 
$$
\sigma _\rho (\Sigma (\rho ;k,j))  \leq  \sigma _\rho (Q_\rho (w_\rho ,C_22^{k+j}(-r(z))))  \leq  C_3(-2^{k+j}r(z))^n
\tag 2.16
$$
\noindent
for every  $2^{k-1}(-r(z))   < \rho  \leq 2^k(-r(z))$   under the condition    $2^k(-r(z)) \leq \tau $.
\medskip
Now consider  $\Omega $  as a domain in the real space   ${\bold R}^{2n}$  under the usual identification.  We know that  $(\nabla r)(x) \neq 0$  for every  
$x \in \partial \Omega $.   For each  $a \in \partial \Omega $,  there is a  $j = j(a) \in \{1,2,\dots ,2n-1,2n\}$  such that the map  
$$
F_a(x_1,x_2,\dots ,x_{2n-1},x_{2n})   =  (x_1,\dots ,x_{j-1},-r(x_1,x_2,\dots ,x_{2n-1},x_{2n}),x_{j+1},\dots ,x_{2n})
$$ 
\noindent
from  $\Omega $  to    ${\bold R}^{2n}$  has the property that the  derivative   $(DF_a)(a)$   is invertible.  Thus
there is an open neighborhood   $U_a$  of  $a$  in    ${\bold R}^{2n}$  such that the inverse mapping theorem
holds on   $U_a$   for   $F_a$.   Shrinking  $U_a$  slightly if necessary,
we may assume that     $DF_a^{-1}$  is bounded on  $F_aU_a$.
By the  compactness of $\partial \Omega $,  there is a finite subset  $A \subset \partial \Omega $  and a  $0 < \tau _1 \leq \tau $   such that  
$$
\bigcup _{a\in A}U_a  \supset  \{w \in \Omega :  0 < -r(w) \leq  \tau _1\}.
\tag 2.17
$$
\noindent
To complete the proof of the proposition, consider the following two cases.

\medskip
(1)   Suppose that  $-2^kr(z) \geq \tau _1$.  Then  the conclusion of the proposition is trivial.   (2)  Suppose that  
 $-2^kr(z) < \tau _1$.     In this case (2.17)  gives us  
$$
W_{z;k,j}  =  \bigcup _{a\in A}\{U_a\cap W_{z;k,j}\}.
$$
\noindent
Since   $A$  is an {\it a priori}  determined finite set,  it suffices to estimate the volume of  $U_a\cap W_{z;k,j}$   
for each   $a \in A$.   Given any  $a \in A$,  there is a  $j = j(a)$  such that  for  $0 < \rho < \tau _1$,
$$
S_\rho \cap U_a  =  \{F_a^{-1}(x_1,\dots ,x_{j-1},\rho ,x_{j+1},\dots ,x_{2n})  :  (x_1,\dots ,x_{j-1},\rho ,x_{j+1},\dots ,x_{2n}) \in F_aU_a\}.
$$
\noindent
By (2.16),   for each  $0 < \rho < \tau _1$,   the  real $(2n-1)$-dimensional volume of the set
$$
\{(x_1,\dots ,x_{j-1},x_{j+1},\dots ,x_{2n})  :   (x_1,\dots ,x_{j-1},\rho ,x_{j+1},\dots ,x_{2n}) \in F_a(U_a\cap \Sigma (\rho ;k,j))\}
$$
\noindent
does not exceed   $C_9(-2^{k+j}r(z))^n$.   Hence we have the (real) $2n$-dimensional  volume estimate
$$
v(F_a(U_a\cap W_{z;k,j}))   \leq C_9(-2^{k+j}r(z))^n\cdot 2^k(-r(z)).
$$
\noindent
Applying  $F_a^{-1}$,  we find that
$$
v(U_a\cap W_{z;k,j})   \leq  C_{10}2^{nj}(-2^kr(z))^{n+1}.
$$
\noindent
Since  card$(A) < \infty $,  this completes the proof.   $\square $
\medskip
\noindent
{\bf Definition  2.10.}     (i)     Let  $a$   be a positive number.   
A subset   $\Gamma $  of  $\Omega $  is said to be  $a$-separated if
$D(z,a)\cap D(w,a) = \emptyset $  for all distinct elements   $z$,  $w$  in  $\Gamma $.

\noindent
(ii)   A subset   $\Gamma $  of    $\Omega $  is simply said to be separated if it is $a$-separated  
for some  $a > 0$.
\medskip
\noindent
{\bf Lemma 2.11.}    (1)  {\it  For any pair of}    $0 < a < \infty $  {\it  and}   $0 < R < \infty $,   {\it  there is a natural number}  $N = N(a,R)$    
{\it  such that for every}  $a$-{\it  separated set}   $\Gamma $  {\it  in}   $\Omega $   {\it  and every}  $z \in \Omega $,   {\it  we have}  
$$
\text{card}\{u \in \Gamma :  d(u,z) \leq  R\}  \leq  N.
$$
\noindent
(2)   {\it  For any pair of}    $0 < a \leq  R < \infty $,   {\it  there is a natural number}  $m = m(a,R)$  {\it  such that every}      
$a$-{\it  separated set}   $\Gamma $  in  $\Omega $  {\it  admits a partition}  $\Gamma = \Gamma _1\cup \cdots \cup \Gamma _m$  
{\it  with the property that for every}  $j \in \{1,\dots ,m\}$,  {\it  the set}   $\Gamma _j$   {\it  is}  $R$-{\it  separated}.
\medskip
\noindent
{\it Proof}.   By Proposition 2.6,  any integer   $N \geq C(R+a)/c(a)$  will do for (1).  Then,  by(1),  for any   $0 < a \leq  R < \infty $,  
there is an  $m \in $  {\bf N}  such that  if   $\Gamma $  is any  $a$-separated set in   $\Omega $,  
then  card$\{u \in \Gamma : d(u,v) \leq 2R\} \leq m$  for every  $v \in \Gamma $.  By a standard maximality argument,   $\Gamma $  
admits a partition  $\Gamma = \Gamma _1\cup \cdots \cup \Gamma _m$  such that for every  $j \in \{1,\dots ,m\}$,  the conditions  
$u, v \in \Gamma _j$   and  $u \neq v$  imply  $d(u,v) > 2R$.  Thus each  $\Gamma _j$  is  $R$-separated,
proving (2).   $\square $

\bigskip
\centerline{{\bf 3. Forelli-Rudin estimates on}  $\Omega $}
\medskip
We will need the familiar functions
$$
X(z,w) = -r(w)  - \sum_{j=1}^n\frac{\partial r(w)}{\partial w_j}(z_j-w_j)   
- {1\over 2}\sum _{j,k=1}^n \frac{\partial^2r(w)}{\partial w_j\partial w_k}(z_j-w_j)(z_k-w_k),
\tag 3.1
$$
$$
\rho(z,w)  =  |z-w|^2 +   |\langle z - w,(\bar \partial r)(z)\rangle |
$$
\noindent
and
$$
F(z,w)   =  |r(z)|  +   |r(w)|  +  \rho (z,w)
$$
\noindent
associated with  $\Omega $  and  $r$,  which are standard fixtures on strongly pseudo-convex domains.
\medskip
\noindent
{\bf Lemma 3.1.}  [18,20]    {\it  There is a}  $\delta > 0$  {\it  such that}
$$
|X(z,w)| \approx |r(z)| + |r(w)| + |\text{Im}X(z,w)| + |z - w|^2  \approx F(z,w)
$$
{\it  in the region}   ${\Cal R}_{\delta}  = \{(z,w) \in \Omega \times \Omega  : |r(z)|+ |r(w)| + |z-w| < \delta\}$. 
\medskip
Below is what one usually refers to as the Forelli-Rudin estimates:
\medskip
\noindent
{\bf Lemma 3.2.}   [18,20]    {\it  Let}   $a \in $  {\bf R}   {\it   and}  $\kappa  > - 1$.  {\it  Then for} $z\in\Omega$,
$$
\int_{\Omega}\frac{|r(w)|^{\kappa}}{F(z,w)^{n+1+\kappa+a}}dv(w)\approx
\left\{
\matrix
1   &\text{\it if}  \ a <  0   \\
\log \big\{|r(z)|^{-1}\big\}    &\text{\it if}  \  a  =  0\\
|r(z)|^{-a}  &\text{\it if}  \  a > 0
\endmatrix   
\right.  .
$$
\medskip
Recall  that  for any  $z \in \overline{\Omega }$ with   $0 \leq -r(z) \leq  \theta $,    
we have  $(\bar \partial r)(z) \neq 0$  as a vector in  ${\bold C}^n$.
\medskip
\noindent
{\bf  Definition 3.3.}   For  $z \in \Omega $  satisfying the condition  $0 < -r(z) < \theta $,  
let  $u_z$   denote the unit vector  $(\bar \partial r)(z)/|(\bar \partial r)(z)|$  in  ${\bold C}^n$.
\medskip
\noindent
{\bf Lemma 3.4.}    {\it  There exist constants}   $\delta _0 > 0$   {\it  and}   $0 < C_{3.4} < \infty $   {\it  such that if}    
$z \in \Omega $  {\it  and}   $i \in {\bold Z}_+$   {\it  satisfy the condition}  $- 2^{i+1}r(z) < \delta _0$,  
{\it  then  for every}   $x \in [1,2]$  {\it  we have}  
$$
d(z + 2^ir(z)u_z,z + x2^ir(z)u_z)  \leq  C_{3.4}.
\tag 3.2
$$
\noindent
{\it  Moreover},   {\it  if}    $-r(z) < \delta _0$,  {\it  then}   $d(z,z + su_z)  \leq  C_{3.4}$  {\it  for every}   $s \in [r(z),0]$.
\medskip
\noindent
{\it Proof}.
Let   $z \in \Omega $  be such that   $0 < -r(z) < \theta $.   Then 
$$
{d\over dt}r(z + tu_z)  = 2\text{Re}\langle u_z,(\bar \partial r)(z + tu_z)\rangle   =   2|(\bar \partial r)(z)| + O(|t|).
$$
Thus there is a  $0 < \delta _0 <  \theta $    such that if  $-r(z) < \delta _0$,  then
$$
\text{the function} \ \  t \mapsto r(z + tu_z)  \ \  \text{is increasing on}  \ \  [- \delta _0, \delta _0].
$$
\noindent
Now let   $z \in \Omega $  and   $i \in {\bold Z}_+$ be  such that  $-2^{i+1}r(z) <  \delta _0$.   Let   $x \in [1,2]$.   Then for any
$s \leq s'$  in  the interval  $[x2^ir(z),2^ir(z)]$,   the above monotonicity guarantees
$r(z +su_z)$  $\leq $  $r(z+s'u_z)$,  i.e.,  $- r(z +s'u_z) \leq -r(z+su_z)$.   
For such a pair of  $s$  and  $s'$,  it follows from (2.3),  (2.4)  and the above monotonicity that
$$
d(z +s'u_z,z +su_z)  \leq  C{|(z +s'u_z) - (z +su_z)|\over -r(z +s'u_z)}  \leq  C{|-2^{i+1}r(z) - 2^i(-r(z))|\over -2^ir(z)}  =  C,
$$
\noindent
which proves (3.2).   Similarly,  if  $-r(z) < \delta _0$,     then for every   
$s \in [r(z),0]$  we have   $r(z+su_z) \leq r(z)$,  i.e.,  $-r(z+su_z) \geq -r(z)$.  Hence the same argument shows that  
$d(z,z+su_z) \leq C|z - (z+su_z)|/(-r(z)) \leq C$.      This proves the lemma.   $\square $
\medskip
\noindent
{\bf Lemma 3.5.}    {\it  There exist constants}    $0 < c_{3.5}  \leq 1$   {\it  and}    $0 < \delta _1 \leq \delta _0$,  
{\it  where}  $\delta _0$  {\it  was given in Lemma} 3.4,  {\it  such that if}  
$z \in \Omega $   {\it  satisfies the condition}  $-r(z) < \delta _1$   {\it  and if}  $-\delta _1 \leq t \leq 0$,     {\it  then}
$$
-r(z + tu_z) +  r(z)  \geq c_{3.5}|t|.
$$
\medskip
\noindent
{\it Proof}.   Taylor expansion gives us
$$
r(z + tu_z) = r(z) + 2t\text{Re}\langle u_z,(\bar \partial r)(z)\rangle + O(t^2)  =  r(z) + 2t|(\bar \partial r)(z)| + O(t^2). 
$$
\noindent
In other words,  $r(z + tu_z)  - r(z) = \{2|(\bar \partial r)(z)| + O(t)\}t$.   From this the desired conclusion becomes obvious.  
$\square $
\medskip
\noindent
{\bf Proposition 3.6.}   {\it  There is a constant}   $0 < C_{3.6} < \infty $   {\it  such that if}  
$z, w \in \Omega $  {\it  satisfy the conditions}   $r(z) = r(w)$   {\it  and}  
$|z - w|^2 + |\langle z-w,(\bar \partial r)(z)\rangle | \leq -2^jr(z)$,  $j \in {\bold Z}_+$,   {\it  then}  $d(z,w) \leq C_{3.6}(1+j)$.
\medskip
\noindent
{\it Proof}.    Recall the sets  $U_1,\dots ,U_m$  from  Proposition 2.7,  which are an open cover of   
$H = \{\zeta \in {\bold C}^n :  0 \leq -r(\zeta ) \leq \theta /2\}$  in  ${\bold C}^n$.  By general topology, there is an  $a_1 > 0$  
such that for any pair of  $z, w \in H$,  if  $|z - w| < a_1$,  then there is an  $i(z,w) \in \{1,\dots ,m\}$  such that  
$z, w \in U_{i(z,w)}$.   Another elementary exercise gives us a pair of constants  $0 < \theta _0 < \theta /2$  
and  $0 < a < \min \{1,a_1/4\}$ which have the following property:  Suppose that  $z, w \in  \{\zeta \in \Omega : -r(\zeta ) < \theta _0\}$
and that  $z', w' \in \Omega $.  If the inequalities  $|z - w| <a$,  $|z - z'| < a$  and  $|w - w'| < a$  hold,  then there is an  $i^\ast  \in \{1,\dots ,m\}$  
such that    $z', w' \in U_{i^\ast }$.
\medskip
Define  $\delta = \min \{\theta _0,a^2,\delta _1/2\}$,  where   $\delta _1$  is given in Lemma 3.5.  We divide the rest of the proof  
into two cases.   
\medskip
(1)   Suppose that   $-2^jr(z) < \delta c_{3.5}$, where  $c_{3.5}$  is also from Lemma 3.5.     
Then    $|z - w|^2 < \delta $,  which implies   $|z - w| < a$.
We set   $s = 2^jr(z)/c_{3.5}$.  Since  $-\delta _1 < s < 0$,  by Lemma 3.5,
$$
- r(z+su_z)  \geq  c_{3.5}|s| =  -2^jr(z).
$$ 
\noindent
Since    $j \geq 0$,  there is an  $s(z) \in [s,0]$   such that    $-r(z + s(z)u_z) = -2^jr(z)$.  We set  
$z' = z + s(z)u_z$.   Then  $r(z')  = 2^jr(z)$   and   $|z - z'| = |s(z)| \leq |s| < \delta < a$.   
Since   $|s| < \delta \leq  \delta _0/2$,  it follows from Lemma 3.4 that   $d(z,z') \leq  C_1(1+j)$.  
\medskip
Similarly,  since  $r(w) = r(z)$,  
there is an  $s(w) \in [s,0]$  such that if we set  $w' = w + s(w)u_w$,   then  $r(w') = 2^jr(z)$,
$|w - w'|  \leq |s|  < a$  and   $d(w,w') \leq  C_1(1+j)$.  
\medskip
Since   $\bar \partial r$  satisfies a Lipschitz condition on  $\Omega $,  from the conditions 
$$
|z - w|^2 + |\langle z-w,(\bar \partial r)(z)\rangle | \leq -2^jr(z),
$$
\noindent
$|z - z'| \leq  -2^jr(z)/c_{3.5}$,   $|w - w'| \leq -2^jr(z)/c_{3.5}$  and  $-2^jr(z) < \delta c_{3.5}$   it is easy to deduce
$$
|z' - w'|^2 + |\langle z'-w',(\bar \partial r)(z')\rangle | \leq -C_22^jr(z) =  C_2|r(z')|.
$$
\noindent
Since   $|z - z'| < a$,    $|w - w'| < a$   and   $|z - w| < a$,  by the first paragraph,  there is an  $i^\ast $  $\in $  
$\{1,\dots ,m\}$  such  that   $z', w' \in U_{i^\ast }$.   Hence it follows from the above bound and Proposition 2.7 that    
$d(z',w') \leq C_3$.  Combining this with the last two paragraphs,  we obtain  
$$
d(z,w)  \leq  d(z,z') + d(z',w') + d(w,w')  \leq  2C_1(1+j) + C_3 \leq C_4(1+j).
$$
\noindent
This proves the proposition under the condition   $-2^jr(z) < \delta c_{3.5}$.
\medskip
(2)   Suppose that   $-2^jr(z) \geq \delta c_{3.5}$.   (2.a)   Further, suppose that   $-r(z) \geq \delta c_{3.5}/2$.   In this case the  
conclusion is trivial,  for  $\{\zeta \in {\bold C}^n  : - r(\zeta ) \geq \delta c_{3.5}/2\}$  is a compact subset of   $\Omega $.  
(2.b)  Suppose that  $-r(z) < \delta c_{3.5}/2$,  i.e.,    $-2r(z) < \delta c_{3.5}$.   Let  $j_0$  be the largest natural number such that  
$-2^{j_0}r(z) < \delta c_{3.5}$.   Then obviously  $j_0 < j$.   By the work in case (1) we know that there are  $\tilde z, \tilde w \in \Omega $   such that  
$r(\tilde z) = 2^{j_0}r(z) = r(\tilde w)$,  $d(z,\tilde z) \leq C_1(1+j_0)$   and  $d(w,\tilde w) \leq C_1(1+j_0)$.  The choice of  $j_0$   ensures that
$-2^{j_0+1}r(z) \geq \delta c_{3.5}$,  which means  $-r(\tilde z) \geq \delta c_{3.5}/2$  and   $-r(\tilde w) \geq \delta c_{3.5}/2$.  
Hence  $d(\tilde z,\tilde w) \leq C_5$,  and consequently
$$
d(z,w)  \leq  d(z,\tilde z) + d(\tilde z,\tilde w) + d(w,\tilde w)  \leq  2C_1(1+j_0) + C_5 \leq C_6(1+j)
$$
\noindent
in this subcase.   This completes the proof.   $\square $
\medskip
\noindent
{\bf Lemma 3.7.}    {\it  There is a}  $1 \leq  C_{3.7} < \infty $  {\it  such that if}  
$z, w \in \Omega $    {\it  satisfy the conditions}  $2^{k-1}(-r(z))   \leq  - r(w) \leq 2^k(-r(z))$     {\it  and} 
$|z - w|^2 + |\langle z - w,(\bar \partial r)(z)\rangle | < 2^{k+j}(-r(z))$,  {\it  where} $k \in $  {\bf Z}    {\it  and}   
$j \in {\bold Z}_+$,   {\it  then}   $d(z,w) <  C_{3.7}(1+|k|+j)$.
\medskip
\noindent
{\it Proof}.     (1)   First, let us consider the case where  $k \geq 1$.   (1.a)  Further, suppose that
$2^k(-r(z)) \geq c_{3.5}\delta _1/4$,  where  $c_{3.5}$  and  $\delta _1$  are the constants in Lemma 3.5.  Then the condition  
$2^{k-1}(-r(z))   \leq  - r(w) \leq 2^k(-r(z))$   implies  $-r(w) \geq c_{3.5}\delta _1/8$.  If we also have  
$-r(z) \geq c_{3.5}\delta _1/4$,  then of course,  $d(z,w) \leq C_1$,  regardless of other conditions.  
Suppose that    $-r(z) < c_{3.5}\delta _1/4$.   Let  $k'$   be the largest integer such that   $-2^{k'+1}r(z) < c_{3.5}\delta _1$.  
Set
$$
z' = z + 2^{k'+1}r(z)u_z.
$$
\noindent
Since    $-2^kr(z) \geq c_{3.5}\delta _1/4$,  we have  $k' +1 <  k +2$,  i.e.,  $k' \leq k$.     It follows from 
Lemma 3.4  that   $d(z,z') \leq C_{3.4}(k'+2) \leq C_{3.4}(k+2)$.  By Lemma 3.5,  we have
$c_{3.5}2^{k'+1}|r(z)|   \leq  -r(z')$. 
The choice of  $k'$ ensures that    $-2^{k'+2}r(z) \geq c_{3.5}\delta _1$.  Hence the above implies
$$
c_{3.5}^2\delta _1 /2  \leq  -r(z').
$$
\noindent
Thus  $d(z',w) \leq C_2$,  and consequently  $d(z,w) \leq C_{3.4}(k+2) + C_2 \leq C_3k$ in this subcase.
\medskip
(1.b)  Suppose that   $2^k(-r(z)) < c_{3.5}\delta _1/4$.   Then  $-r(w) \leq 2^k(-r(z)) < c_{3.5}\delta _1/4$.  
By Lemma 3.5,   we have
$$
c_{3.5}\delta _1 \leq   -r(z - \delta _1u_z).
$$
\noindent
Hence   $-r(z - \delta _1u_z)  > -r(w)$.   Since  $-r(z) \leq 2^{k-1}(-r(z))   \leq  -r(w)$,  there is an  $s \in [-\delta _1,0]$  such that  
$r(z + su_z) = r(w)$.   Also, Lemma 3.5 tells us that  
$$
c_{3.5}|s|  \leq   -r(z + su_z) = - r(w)  \leq  -2^kr(z).
\tag 3.3
$$
\noindent
Thus     $|s| \leq -c_{3.5}^{-1}2^kr(z)$.  Now the condition  $2^k(-r(z)) < c_{3.5}\delta _1/4$  implies  $-c_{3.5}^{-1}2^kr(z)  \leq \delta _1/4$.   
Therefore it follows from Lemma 3.4 and the inequality  $|s| \leq -c_{3.5}^{-1}2^kr(z)$  that
$$
d(z,z + su_z)  \leq  C_{3.4}\{1 + C_3\log (2^k/c_{3.5})\}   \leq  C_4k.
$$
\noindent
Thus what remains for this subcase is to show that
$$
d(z + su_z,w)  \leq  C_5(1+j).
$$
\noindent
For convenience, let us denote   $\zeta = z + su_z$.   Since  $r(\zeta ) = r(w)$,
to prove the above inequality,   by Proposition 3.6,  it suffices to show that
$$
|\zeta - w|^2 + |\langle \zeta - w, (\bar \partial r)(\zeta )\rangle |  \leq  C_62^j(-r(\zeta )).
\tag 3.4
$$
\noindent
By (3.3),   
$|\zeta - z|  =  |s|  \leq  c_{3.5}^{-1}2^k(-r(z))  \leq  (2/c_{3.5})(-r(\zeta ))$.
Since   $\Omega $  is bounded,   we have  $|\zeta - z|^2 \leq C_7|\zeta -z|$.   Therefore
$$
|\zeta - w|^2  \leq  2|\zeta - z|^2 + 2 |z - w|^2  \leq C_8(-r(\zeta )) + 2^{k+j+1}(-r(z))  \leq  C_92^j(-r(\zeta )).
$$
\noindent
We have  $|\langle z - w,(\bar \partial r)(z)\rangle | < 2^{k+j}(-r(z))$  by assumption.  Also,
$$
|\langle \zeta - w, (\bar \partial r)(\zeta )\rangle  - \langle z- w, (\bar \partial r)(z)\rangle |  \leq  C_{10}|\zeta - z| \leq  C_{11}(-r(\zeta )).
$$
\noindent
Since  $2^{k-1}(-r(z))   \leq  - r(\zeta )$,  these inequalities prove (3.4).  Thus the case   $k \geq 1$  is proved.
\medskip
(2)  Now suppose that   $k \leq 0$.   Note that the condition  $2^{k-1}(-r(z))   \leq  - r(w) \leq 2^k(-r(z))$  
implies  that  $2^{-k}(-r(w)) \leq -r(z) \leq 2^{-k+1}(-r(w))$.   Also, the condition  $|z - w|^2 + |\langle z - w,(\bar \partial r)(z)\rangle | < 2^{k+j}(-r(z))$ 
can be rewritten as
$$
|z - w|^2 + |\langle z - w,(\bar \partial r)(z)\rangle | < 2^{1+j}(-r(w)).
$$
\noindent
Since  $|\langle z - w,(\bar \partial r)(z)\rangle  -  \langle z - w,(\bar \partial r)(w)\rangle | \leq  C_{12}|z-w|^2$, 
we now have  
$$
|z - w|^2 + |\langle z - w,(\bar \partial r)(w)\rangle | < C_{13}2^j(-r(w)).
$$
\noindent
Thus case (2)  follows from case (1) by reversing the roles of  $z$  and  $w$.   $\square $
\medskip
We  need the following  ``vanishing" version of Lemma 3.2.
\medskip
\noindent
{\bf Lemma 3.8.}    {\it  Given any}  $a > 0$   {\it  and}  $\kappa  > - 1$,    {\it  there are}   $0 < C < \infty $  
{\it  and}  $s > 0$   {\it  such that}  
$$
\int_{\Omega \backslash D(z,R)}{|r(w)|^\kappa |r(z)|^a\over F(z,w)^{n+1+\kappa+a}}dv(w)
\leq  C2^{-sR}
$$
\noindent
{\it  for all}   $z \in \Omega $   {\it  and}  $R \geq 3C_{3.7}$,  {\it  where}   $C_{3.7}$  {\it  is the constant in Lemma}  3.7.
\medskip
\noindent
{\it Proof}.    For  $z \in \Omega $  and  $k \in {\bold Z}$   we define the sets
$$
\align
Z_{z;k,0}  =  \ &\{w \in \Omega :   2^{k-1}(-r(z))   \leq  - r(w) < 2^k(-r(z)) \\
&\quad \  \text{and} \ \  |z - w|^2 + |\langle z - w,(\bar \partial r)(z)\rangle | < 2^{k}(-r(z))\}   \quad \text{and}  \\
Z_{z;k,j}  =  \ &\{w \in \Omega :   2^{k-1}(-r(z))   \leq  - r(w) < 2^k(-r(z)) \\
&\quad \ \text{and} \ \  2^{k+j-1}(-r(z)) \leq |z - w|^2 + |\langle z - w,(\bar \partial r)(z)\rangle | < 2^{k+j}(-r(z))\},  \quad  j \geq 1. 
\endalign
$$
\noindent
By the definition of  $F(z,w)$,      for all   $k \geq  0$  and  $j \geq 0$,  if  $w \in Z_{z;k,j}$,  then
$$
{|r(w)|^\kappa |r(z)|^a\over F(z,w)^{n+1+\kappa+a}}  \leq  C_1{(2^k|r(z)|)^\kappa  |r(z)|^a\over (2^{k+j}|r(z)|)^{n+1+\kappa + a}}
=  {C_1\over 2^{(n+1+a)k}2^{(n+1+\kappa +a)j}|r(z)|^{n+1}}.
$$
\noindent
In the case  $k <  0$,  $j \geq 0$   and  $w \in Z_{z;k,j}$,   we have
$$
{|r(w)|^\kappa |r(z)|^a\over F(z,w)^{n+1+\kappa+a}}  \leq  {C_2(2^k|r(z)|)^\kappa  |r(z)|^a\over (|r(z)| + 2^{k+j}|r(z)|)^{n+1+\kappa + a}}
=  {C_22^{\kappa k}\over (1+2^{k+j})^{n+1+\kappa +a}|r(z)|^{n+1}}.
$$
\noindent
By Proposition 2.9,  $v(Z_{z;k,j})  \leq C_32^{nj}(-2^kr(z))^{n+1}$.   Thus if  $k \geq  0$  and  $j \geq 0$,   then
$$
\int_{Z_{z;k,j}}{|r(w)|^\kappa |r(z)|^a\over F(z,w)^{n+1+\kappa+a}}dv(w)
\leq   {C_1C_32^{nj}(-2^kr(z))^{n+1}\over 2^{(n+1+a)k}2^{(n+1+\kappa +a)j}|r(z)|^{n+1}}     
=   {C_4\over 2^{ak}2^{(1+\kappa + a)j}}.
\tag 3.5
$$
\noindent
Similarly,  in the case  $k <  0$ and  $j \geq 0$,  we have
$$
\int_{Z_{z;k,j}}{|r(w)|^\kappa |r(z)|^a\over F(z,w)^{n+1+\kappa+a}}dv(w)
\leq  {C_22^{\kappa k}C_32^{nj}(-2^kr(z))^{n+1}\over (1+2^{k+j})^{n+1+\kappa +a}|r(z)|^{n+1}}   
=  {C_52^{n(k+j)}2^{(1+\kappa )k}\over (1+2^{k+j})^{n+1+\kappa +a}}
\tag 3.6
$$
\noindent
and   $2^{(1+\kappa )k}$ $=$  $2^{(1+\kappa )k/2}\cdot 2^{(1+\kappa )(k+j)/2}\cdot 2^{-(1+\kappa )j/2}$. 
Let  $R \geq 3C_{3.7}$.   By Lemma 3.7,  the condition   $Z_{z;k,j}\backslash D(z,R) \neq \emptyset $  
implies either   $|k| \geq (2C_{3.7})^{-1}R$  or  $j \geq (2C_{3.7})^{-1}R$.   Therefore
$$
\int_{\Omega \backslash D(z,R)}{|r(w)|^\kappa |r(z)|^a\over F(z,w)^{n+1+\kappa+a}}dv(w)
\leq  \sum _{(k,j)\in E(R)}\int_{Z_{z;k,j}}{|r(w)|^\kappa |r(z)|^a\over F(z,w)^{n+1+\kappa+a}}dv(w),
$$
\noindent
where  $E(R) = \{(k,j) \in {\bold Z}\times {\bold Z}_+$  :  $\text{either}$  $|k| \geq (2C_{3.7})^{-1}R$  $\text{or}$  $j \geq (2C_{3.7})^{-1}R\}$.   
Using (3.5)  and  (3.6),  it is now  elementary to verify that the lemma holds for every  
$0 < s < (2C_{3.7})^{-1}\min \{a,(1+\kappa )/2\}$.    $\square $
\medskip
\noindent
{\bf Lemma 3.9.}    {\it  Given any}  $a > 0$   {\it  and}  $\kappa  > - 1$,    {\it  there is a}  $0 < C < \infty $    {\it  such that}   
$$
\int_\Omega d(z,w){|r(w)|^\kappa |r(z)|^a\over F(z,w)^{n+1+\kappa+a}}dv(w)
\leq  C
$$
\noindent
{\it  for every}    $z \in \Omega $.

\medskip
\noindent
{\it Proof}.    Given any  $z \in \Omega $,  define   $E_0 = D(z,3C_{3.7})$   and  
$$
E_i =  D(z,(3+i)C_{3.7})\backslash D(z,(3+i-1)C_{3.7})
$$
\noindent
for  $i \geq 1$.   For each   $i \in {\bold Z}_+$,  if  $w \in E_i$,  then   $d(z,w) < (3+i)C_{3.7}$.   Hence
$$
\align
\int_\Omega d(z,w){|r(w)|^\kappa |r(z)|^a\over F(z,w)^{n+1+\kappa+a}}dv(w) 
&= \sum _{i=0}^\infty \int_{E_i}d(z,w){|r(w)|^\kappa |r(z)|^a\over F(z,w)^{n+1+\kappa+a}}dv(w)   \\
&\leq  \sum _{i=0}^\infty  (3+i)C_{3.7}\int_{E_i}{|r(w)|^\kappa |r(z)|^a\over F(z,w)^{n+1+\kappa+a}}dv(w) .
\endalign
$$ 
\noindent
We now apply  Lemma 3.2 to the term where  $i = 0$  and  Lemma 3.8 to the terms where  $i \geq 1$.  
The result of  this is
$$
\int_\Omega d(z,w){|r(w)|^\kappa |r(z)|^a\over F(z,w)^{n+1+\kappa+a}}dv(w) 
\leq   3C_{3.7}C_0 + \sum _{i=1}^\infty  (3+i)C_{3.7}C2^{-s(3+i-1)C_{3.7}}. 
$$ 
\noindent
Since Lemma 3.8  guarantees that  $s > 0$,  the right-hand side is finite.   $\square $
\medskip
\noindent
{\bf Lemma 3.10.}   {\it  There exist constants}   $0 < a_0 <  1/2$   {\it  and}   $0 < C_{3.10} < \infty $  {\it  such that  for any}  
$z, z', w, w' \in \Omega $  {\it  satisfying the conditions}  $d(z,z') < a_0$   {\it  and}  $d(w,w') < a_0$,  {\it  we have}
$$
F(z,w)  \leq  C_{3.10}F(z',w').
$$   
\medskip
\noindent
{\it Proof}.    By Lemma 2.1,  it suffices to consider the case where  $-r(\zeta ) < \theta $  for every  $\zeta \in \{z,z',w,w'\}$.
Since    $|(\bar \partial r)(z) - (\bar \partial r)(w)| \leq  C_1|z - w|$    for all   $z, w \in \Omega $,   there is a  $C_2$   such that
$$
F(w,z)  \leq  C_2F(z,w)
$$
\noindent
for all   $z, w \in \Omega $.   Therefore it suffices to find   $0 < a_0  < 1/2$  and  $C$   such that
$$
F(z,w)  \leq  CF(z',w)
\tag 3.7
$$
\noindent
for all  $z, z', w \in \Omega $    satisfying the condition  $d(z,z') < a_0$.   
Let  $z, z' \in \Omega $  be such that   $d(z,z') < a$  for some  $0 < a  < 1/2$. 
By Lemma 2.1,   $-r(z) \leq  (4/c_{2.1})(-r(z'))$.   Hence,  to prove (3.7),  it suffices to consider the case  
where    
$$
\rho (z,w) = |z - w|^2 + |\langle z - w,(\bar \partial r)(z)\rangle | \geq -r(z).
\tag 3.8
$$
\noindent
Proposition 2.5 tells us that    $z' = z + u + v$   with  
$u \perp (\bar \partial r)(z)$  and  $v \in \{\xi (\bar \partial r)(z) : \xi \in {\bold C}\}$  satisfying the conditions  
$|u|$  $\leq $   $C_{2.5}a|r(z)|^{1/2}$  and  $|v|  \leq  C_{2.5}a|r(z)|$.   Therefore
$$
|z - z'|^2  =  |u|^2 + |v|^2  \leq  C_3a^2|r(z)|.
$$
\noindent
By a simple completion of square,  we find that
$$
|z' - w|^2  \geq  |z' - z|^2   - 2|z'-z||z-w|  +   |z - w|^2  \geq  (1/2) |z - w|^2 -  |z' - z|^2.
$$
\noindent
Also,    $|\langle z - z',(\bar \partial r)(z)\rangle | =   |\langle v,(\bar \partial r)(z)\rangle |  \leq C_4a|r(z)|$.    Consequently
$$
\align
|\langle z' - w,&(\bar \partial r)(z')\rangle | \geq  |\langle z' - w,(\bar \partial r)(z)\rangle | - C_1|z' -w||z - z'| \\
&\geq   |\langle z - w,(\bar \partial r)(z)\rangle | - C_4a|r(z)| - C_1|z -w||z - z'| - C_1|z - z'|^2.
\endalign
$$
\noindent
Combining the above inequalities,  we see that there is a   $1 \leq C_5 < \infty $  such that
$$
\align
\rho &(z',w)  \geq {1\over 2}\rho (z,w) - C_5a\{\rho ^{1/2}(z,w)|r(z)|^{1/2}  +  |r(z)|\}  \\
&= {1\over 2}(1 - C_5a)\rho (z,w)  + {1\over 2}C_5a\big\{\rho (z,w) - 2\rho ^{1/2}(z,w)|r(z)|^{1/2} + |r(z)|\big\} - {3\over 2}C_5a|r(z)|. 
\endalign
$$
\noindent
Note that the  $\{\cdot \cdot \cdot \}$  above is a complete square.   Thus if  $a \leq (2C_5)^{-1}$,   then
$$
\rho (z',w)  \geq   {1\over 4}\rho (z,w) - {3\over 2}C_5a|r(z)|.
$$
\noindent
Recalling (3.8),  if we further require that  $a \leq (12C_5)^{-1}$,    then    
$$
\rho (z',w)  \geq (1/8)\rho (z,w).
$$
\noindent
Thus if we set  $a_0 = (12C_5)^{-1}$,    then there is a  $C$  such that  (3.7) holds for all  
$z, z', w \in \Omega $    satisfying the condition  $d(z,z') < a_0$.   $\square $ 

\bigskip
\centerline{\bf 4. Estimates related to the Bergman kernel}
\medskip
Let  $K(z,w)$  be the  Bergman kernel  for   $\Omega $.   By definition,    
it has the symmetry   $K(w,z)$  $=$   $\overline{K(z,w)}$.
The following well-known result of Fefferman gives us a good handle on  $K$:
\medskip
\noindent
{\bf Theorem 4.1.}  [10,Theorem 2]  {\it  The Bergman kernel has the form}
$$
K(z,w)  =  C|(\nabla r)(w)|^2\text{det}{\Cal L}(w)X^{-(n+1)}(z,w) + \tilde K(z,w)
$$
\noindent
{\it  on}  ${\Cal R}_\delta = \{(z,w) \in \Omega \times \Omega : |r(z)|+ |r(w)| + |z-w| < \delta\}$  {\it  for some}   $\delta > 0$,  
{\it  where}   ${\Cal L}$  {\it  is the Levi form for the domain}   $\Omega $,  $X$  {\it  is given by} (3.1),  {\it  and}  
$\tilde K$  {\it  is an admissible kernel of weight}  $\geq - n - (1/2)$.   {\it  That is},  {\it  there is a constant}   $C'$  {\it  such that}  
$|\tilde K(z,w)| \leq  C'F(z,w)^{-n-(1/2)}$.   
\medskip
For any  $\delta > 0$,  the Bergman kernel  $K$  is known to be bounded on
$(\Omega \times \Omega )\backslash {\Cal R}_\delta $  [14].    One obvious implication of  Theorem 4.1 is that  
$$
c|r(z)|^{-n-1}  \leq   |K(z,z)|  \leq  C|r(z)|^{-n-1},   \quad  z \in \Omega .
\tag 4.1
$$
\noindent
For each   $z \in \Omega $,  let us denote   $K_z(w) = K(w,z)$.  Then it has the reproducing property 
$$
h(z)  =  \langle h,K_z\rangle
$$
\noindent
for  $h \in L^2_a(\Omega )$.  We write  $k_z$  for the normalized reproducing kernel,  i.e.,  
$k_z = K_z/\|K_z\|$.  
\medskip
\noindent
{\bf Lemma 4.2.}      {\it  Given any}  $0 < \eta < 1/2$  and  $a > 0$,  {\it  there are constants}   
$s > 0$   {\it  and}  $0 < C_{4.2} < \infty $   {\it  such that}
$$
\sup _{z \in \Omega }|r(z)|^{-(n/2)-\eta }\sum  _{w\in \Gamma \backslash D(z,R)}
|r(w)|^{(n/2)+\eta }\bigg({|r(z)|^{1/2}|r(w)|^{1/2}\over F(z,w)}\bigg)^{n+1}  \leq  C_{4.2}2^{-sR}
$$
\noindent
{\it  for every}  $a$-{\it  separated set}  $\Gamma $  {\it  in}  $\Omega $  {\it  and every}  $R \geq 3C_{3.7} + 1$.
\medskip
\noindent
{\it Proof}.     Let     $0 < \eta < 1/2$  and  $a > 0$  be given.   Define  $\alpha = (1/3)\min \{a_0,a\}$,  where  
$a_0$  is the constants  in Lemma 3.10.   
Suppose that  $\Gamma $ is an $a$-separated set   in  $\Omega $.  Then 
$$
D(w,\alpha )\cap D(w',\alpha )  =  \emptyset   \quad  \text{for all}  \ \ w \neq w'  \ \  \text{in}  \ \  \Gamma .
$$ 
\noindent
Applying  Lemmas  2.1 and 3.10,   for  $\zeta \in D(w,\alpha )$  we have
$$
{|r(w)|^{n+(1/2)+\eta }|r(z)|^{(n+1)/2}\over F(z,w)^{n+1}}  \leq    C{|r(\zeta )|^{n+(1/2)+\eta }|r(z)|^{(n+1)/2}\over F(z,\zeta )^{n+1}}.
$$
\noindent
Thus  for  $z \in \Omega $  we have
$$
\align
\sum  _{w\in \Gamma \backslash D(z,R)}
|r(w)|&^{(n/2)+\eta }\bigg({|r(z)|^{1/2}|r(w)|^{1/2}\over F(z,w)}\bigg)^{n+1}   \\
&\leq \sum  _{w\in \Gamma \backslash D(z,R)} {C\over \mu (D(w,\alpha ))}
\int _{D(w,\alpha )}{|r(\zeta )|^{n+(1/2)+\eta }|r(z)|^{(n+1)/2}\over F(z,\zeta )^{n+1}}d\mu (\zeta )  \\
&\leq  {C\over c(\alpha )}\int _{\Omega \backslash D(z,R-\alpha )}{|r(\zeta )|^{-(1/2)+\eta }|r(z)|^{(n+1)/2}\over F(z,\zeta )^{n+1}}dv(\zeta ),
\endalign
$$
\noindent
where the second  $\leq $  is justified by Proposition 2.6.
Applying Lemma 3.8 to the last integral,  the desired conclusion follows.   $\square $
\medskip
\noindent
{\bf Lemma 4.3.}   {\it  There is a  constant}   $0 < C_{4.3} < \infty $    {\it  such that} 
$$
|f(z)|  \leq  C_{4.3}|r(z)|^{-(n+1)/2}\|f\chi _{D(z,1)}\|
$$
\noindent 
{\it  for all}   $f \in L^2_a(\Omega )$  {\it  and}  $z \in \Omega $,  {\it  where}  $\|f\chi _{D(z,1)}\|$  {\it  is the norm of}     
$f\chi _{D(z,1)}$   {\it  in}  $L^2(\Omega )$.
\medskip
\noindent
{\it Proof}.   It is easy to see that the conclusion is trivial if   $-r(z) \geq \theta $.  
Suppose that   $-r(z) < \theta $.  Then Proposition 2.4 provides a  $c > 0$  such that
$$
D(z,1)  \supset z + {\Cal P}((\bar \partial r)(z);c|r(z)|^{1/2},c|r(z)|)
$$
\noindent
for every such   $z$.   Averaging on the polyball,  for $f \in L^2_a(\Omega )$ we have
$$
|f(z)|  \leq  {1\over v({\Cal P}((\bar \partial r)(z);c|r(z)|^{1/2},c|r(z)|))}\int _{z + {\Cal P}((\bar \partial r)(z);c|r(z)|^{1/2},c|r(z)|)}|f|dv.
$$
\noindent
Applying the Cauchy-Schwarz inequality on the right,  the desired conclusion follows.      $\square $
\medskip
\noindent
{\bf Lemma 4.4.}    {\it  Given any complex dimension}  $m \in $  {\bf N},  {\it  there is a constant}    
$0 < C_{4.4}(m) < \infty $  {\it  such that the following bound holds}:   {\it  Let}   $0 < \rho  < \infty $  
{\it  and define}  $B(\rho ) = \{z \in {\bold C}^m : |z| < \rho \}$.   {\it  Then for every}   $u \in {\bold C}^m$  {\it  with}   
$|u| < \rho /2$   {\it  and every analytic function}    $f$  {\it  on}  $B(\rho )$,  {\it  we have}
$$
|f(u) - f(0)| \leq  {|u|\over \rho }\cdot {C_{4.4}(m)\over v_m(B(\rho ))}\int _{B(\rho )}|f|dv_m.
$$
\medskip
\noindent
{\it Proof}.     By  standard integration formulas on the ball,  there is a   $C = C(m)$  such that
$$
|(\partial _jg)(0)|  \leq  {C\over v_m(B(1))}\int _{B(1)}|g|dv_m
\tag 4.2
$$
\noindent
for every analytic function    $g$  on   $B(1) = \{z \in {\bold C}^m : |z| < 1\}$
and  every   $j \in \{1,\dots ,m\}$.   Suppose that  $u = (u_1,\dots ,u_m)$.   If  $f$  is analytic on   $B(\rho )$,  then    
$$
f(u) - f(0)  =  \int _0^1{d\over dt}f(tu)dt =   \int _0^1\sum _{j=1}^m(\partial _jf)(tu)u_jdt.
$$
\noindent
Since    $|u| < \rho /2$,  for every  $t \in [0,1]$  we have   $tu + B(\rho /2) \subset B(\rho )$.  From  (4.2) and the scaling 
properties of    $\partial _j$  and  $dv_m$   we deduce
$$
|(\partial _jf)(tu)|  \leq  {2\over \rho }\cdot {C\over v_m(B(\rho /2))}\int _{tu+B(\rho /2)}|f|dv_m.
$$
\noindent
Since  $v_m(B(\rho /2)) = 2^{-2m}v_m(B(\rho ))$,  we see that the constant  $C_{4.4}(m) = m2^{2m+1}C$  
will do for the lemma.   $\square $

\medskip
\noindent
{\bf Lemma 4.5.}   {\it  There exist constants}   $0 < C_{4.5} < \infty $  {\it  and}  $0 < c_{4.5} < 1$   {\it  such that}  
$$
|f(w) - f(z)|  \leq  C_{4.5}d(z,w)|r(z)|^{-(n+1)/2}\|f\chi _{D(z,1)}\|
$$
\noindent 
{\it  for every pair of}  $z, w\in \Omega $  {\it  satisfying the condition}   $d(z,w) < c_{4.5}$    
{\it  and every}   $f \in L^2_a(\Omega )$.
\medskip
\noindent
{\it Proof}.   By Lemma 2.1,  there is a  $0 < \theta _1 < \theta $  such that  if   $-r(z) \geq \theta $  and  $d(z,w) \leq 1$,  
then  $-r(w)  \geq \theta _1$.   Since    $\{\zeta \in \Omega : -r(\zeta ) \geq \theta _1\}$  
is a compact subset of   $\Omega $,  we see that the case    $-r(z) \geq \theta $  is trivial.
\medskip
Suppose that   $-r(z) < \theta $.  Then Proposition 2.4 provides a  $c > 0$  such that
$$
D(z,1)  \supset z + {\Cal P}((\bar \partial r)(z);c|r(z)|^{1/2},c|r(z)|)
\tag 4.3
$$
\noindent
for every such  $z$.   By Proposition 2.5,  there is an  $0 < \alpha < 1/2$  such that
$$
D(z,\alpha )  \subset z + {\Cal P}((\bar \partial r)(z);(c/2)|r(z)|^{1/2},(c/2)|r(z)|)
\tag 4.4
$$
\noindent
for every such  $z$.   Set  $c_{4.5} = \alpha $.   Let  $w \in \Omega $  be such that  
$d(z,w) < \alpha $.   Then we can write  $w = z + x + y$,  where  $x \perp (\bar \partial r)(z)$  
and  $y \in \{\eta (\bar \partial r)(z) : \eta  \in {\bold C}\}$.  By (4.4),  we have   $|x| < (c/2)|r(z)|^{1/2}$  
and  $|y| < (c/2)|r(z)|$.
\medskip
Let  $f \in L^2_a(\Omega )$  be given.  Define    $F(\xi ) = f(z + \xi + y)$   for   $\xi  \perp (\bar \partial r)(z)$  
with  $|\xi | < c|r(z)|^{1/2}$.   Applying Lemma 4.4 to the case where  $\rho = c|r(z)|^{1/2}$,  we have
$$
|f(w) - f(z+y)|  =  |F(x) - F(0)| \leq  {|x|\over c|r(z)|^{1/2}}\cdot {C_{4.4}(n-1)\over v_{n-1}(B)}\int _{B}|F|dv_{n-1},
\tag 4.5
$$
\noindent
where  $B = \{\xi \in {\bold C}^n : \langle \xi ,(\bar \partial r)(z)\rangle = 0$  and   $|\xi | < c|r(z)|^{1/2}\}$.  
On the other hand,  for every  $\xi \in B$  we have
$$
|F(\xi )| = |f(z + \xi + y)|  \leq  {1\over A(D((c/2)|r(z)|))}\int _{D((c/2)|r(z)|)}|f(z + \xi + y +  \eta u_z)|dA(\eta ),
\tag 4.6
$$
\noindent
where  $u_z = (\bar \partial r)(z)/| (\bar \partial r)(z)|$   and  $D((c/2)|r(z)|) = \{\eta \in {\bold C} : |\eta | < (c/2)|r(z)|\}$.   
Note that   $v_{n-1}(B)$  $=$  $c_{n-1}(c|r(z)|^{1/2})^{2n-2}$  $= a_1|r(z)|^{n-1}$.   Combining (4.5), (4.6) and (4.3),  
we obtain  
$$
|f(w) - f(z+y)|  \leq  {|x|\over |r(z)|^{1/2}}\cdot {C_1\over |r(z)|^{n+1}}\int _{D(z,1)}|f|dv.
$$
\noindent
Applying the Cauchy-Schwarz inequality and Proposition 2.6,  we have
$$
|f(w) - f(z+y)|  \leq  {|x|\over |r(z)|^{1/2}}\cdot {C_1\sqrt{v(D(z,1))}\over |r(z)|^{n+1}}\|f\chi _{D(z,1)}\|
 \leq  {|x|\over |r(z)|^{1/2}}\cdot {C_2\|f\chi _{D(z,1)}\|\over |r(z)|^{(n+1)/2}}.
$$
\noindent
Since  $d(z,w) < \alpha $  and  $\alpha < 1/2$,  (2.6) gives us that  $|x|/|r(z)|^{1/2} \leq C_3d(z,w)$.   Hence
$$
|f(w) - f(z+y)|  \leq   C_4d(z,w)|r(z)|^{-(n+1)/2}\|f\chi _{D(z,1)}\|.
\tag 4.7
$$
\noindent
Next we show that  
$$
|f(z+y) - f(z)|  \leq   C_5d(z,w)|r(z)|^{-(n+1)/2}\|f\chi _{D(z,1)}\|,
\tag 4.8
$$
\noindent
which together with (4.7) will complete the proof of the lemma.
\medskip
To prove (4.8),  we write   $y = \beta u_z$,  where  $\beta \in $  {\bf C}.   Define   $G(\eta ) = f(z + \eta u_z)$  for  
$\eta \in $  {\bf C}  with  $|\eta | < c|r(z)|$.     Now,   applying Lemma 4.4 to the case where  $\rho = c|r(z)|$,  we have
$$
|f(z+y) - f(z)|  =  |G(\beta ) - G(0)|  \leq  {|\beta |\over c|r(z)|}\cdot {C_{4.4}(1)\over A(D(c|r(z)|))}\int _{D(c|r(z)|)}|G|dA,
\tag 4.9
$$
\noindent
where     $D(c|r(z)|) = \{\eta \in {\bold C} : |\eta | < c|r(z)|\}$.      For each    $\eta \in D(c|r(z)|)$   we have
$$
|G(\eta )|  =  |f(z + \eta u_z)|  \leq  {1\over v_{n-1}(B)}\int _B |f(z + \xi + \eta u_z)|dv_{n-1}(\xi ),
\tag 4.10
$$
\noindent
where   $B = \{\xi \in {\bold C}^n : \langle \xi ,(\bar \partial r)(z)\rangle = 0$  and   $|\xi | < c|r(z)|^{1/2}\}$.   Note that  $|\beta | = |y|$.   
Thus (4.9), (4.10) and (4.3) together give us
$$
|f(z+y) - f(z)|  \leq  {|y|\over |r(z)|}\cdot {C_6\over |r(z)|^{n+1}}\int _{D(z,1)}|f|dv.
$$
\noindent  
Since  $d(z,w) < \alpha $  and  $\alpha < 1/2$,   Lemma 2.2  implies that   $|y|/|r(z)| \leq C_7d(z,w)$.   
Applying the Cauchy-Schwarz inequality and Proposition 2.6  on the right-hand side,  we obtain (4.8).  
This completes the proof.   $\square $
\medskip
\noindent
{\bf Proposition 4.6.}   {\it  There is a constant}   $C_{4.6}$  {\it  such that if}  $z, w \in \Omega $  {\it  satisfies the condition}  
$d(z,w) < c_{4.5}$,   {\it  where}  $c_{4.5}$  {\it  was given in Lemma} 4.5,   {\it  then}
$$
|\langle f,k_z - k_w\rangle |  \leq  C_{4.6}d(z,w)\|f\chi _{D(z,1)}\|
$$
\noindent
{\it  for every}   $f \in L^2_a(\Omega )$.   {\it  Consequently},  {\it  if}  $d(z,w) < c_{4.5}$,   {\it  then}   $\|k_z - k_w\| \leq  C_{4.6}d(z,w)$.
\medskip
\noindent
{\it Proof}.  Write  $K_z(\zeta ) = K(\zeta ,z)$,  the unnormalized reproducing kernel.   Note that  Lemma 4.5 implies that  
$\|K_z - K_w\| \leq C_{4.5}d(z,w)|r(z)|^{-(n+1)/2}$   if  $d(z,w) < c_{4.5}$.   Therefore
$$
|\|K_z\| - \|K_w\||  \leq  C_{4.5}d(z,w)|r(z)|^{-(n+1)/2}  \quad \text{if}  \ \   d(z,w) < c_{4.5}.
$$
\noindent
Combining this with (4.1), the condition  $d(z,w) < c_{4.5}$,  and  Lemma 2.1,  we obtain
$$
|\|K_z\|^{-1}  - \|K_w\|^{-1}|  =  {|\|K_z\| - \|K_w\|| \over \|K_z\|\|K_w\|} \leq  C_1d(z,w)|r(z)|^{(n+1)/2}
\tag 4.11
$$
\noindent
when  $d(z,w) < c_{4.5}$.   Let  $f \in L^2_a(\Omega )$.    Then
$$
\align
\langle f,k_z - k_w\rangle &= f(z)\|K_z\|^{-1} - f(w)\|K_w\|^{-1}    \\
&=   (f(z) - f(w))\|K_w\|^{-1} + f(z)(\|K_z\|^{-1} - \|K_w\|^{-1}).
\tag 4.12
\endalign
$$
\noindent
Applying Lemma 4.5,  we have
$$
\align
|f(z) - f(w)|\|K_w\|^{-1}  &\leq  C_{4.5}d(z,w)|r(z)|^{-(n+1)/2}\|f\chi _{D(z,1)}\|\|K_w\|^{-1}  \\
&\leq  C_2d(z,w)\|f\chi _{D(z,1)}\|,
\tag 4.13
\endalign
$$
\noindent
where the second  $\leq $  follows from (4.1),   the condition  $d(z,w) < c_{4.5}$,  and Lemma 2.1.   
On the other hand,    Lemma 4.3 tells us that
$$
|f(z)|  \leq   C_{4.3}|r(z)|^{-(n+1)/2}\|f\chi _{D(z,1)}\|.
$$
\noindent
Combining this with (4.11),  we obtain
$$
|f(z)||\|K_z\|^{-1}  - \|K_w\|^{-1}|   \leq  C_1C_{4.3}d(z,w)\|f\chi _{D(z,1)}\|.
\tag 4.14
$$
\noindent
Obviously, the lemma follows from (4.12), (4.13)  and (4.14).    $\square $
\medskip
\noindent
{\bf Lemma 4.7.}  {\it  There is a}  $c_{4.7} > 0$  {\it  such that for any pair of}    $z, w \in \Omega $,    
{\it  if}  $d(z,w) \leq c_{4.7}$,    {\it  then}  $|\langle k_z,k_w\rangle | \geq 1/2$.
\medskip
\noindent
{\it Proof}.   We have  $1 - \text{Re}\langle k_z,k_w\rangle  = 2^{-1}\|k_z - k_w\|^2$.   By Proposition 4.6,  there is 
a  $c_{4.7} > 0$  such that for any pair of    $z, w \in \Omega $,  if  $d(z,w) \leq c_{4.7}$, then  $\|k_z - k_w\| \leq 1$.  
Thus if  $d(z,w) \leq c_{4.7}$,  then   $1 - \text{Re}\langle k_z,k_w\rangle  \leq 1/2$,  which implies  
$|\langle k_z,k_w\rangle | \geq 1/2$.   $\square $

\bigskip
\centerline{\bf 5. Discrete sums}
\medskip
We now consider operators constructed from the Bergman kernel.
\medskip
\noindent
{\bf Lemma 5.1.}    {\it  There is a constant}  $0 < C_{5.1} < \infty $  {\it  such that the following estimate holds}:  
{\it  Let}  $\Gamma $  {\it  be any}  $1$-{\it  separated set in}  $\Omega $.   {\it  Suppose that}  $\{e_z : z \in \Gamma \}$    
{\it  is an orthonormal set and}  $\{c_z : z \in \Gamma \}$  {\it  is a bounded set of complex coefficients}.   {\it  Then}
$$
\bigg\|\sum _{z\in \Gamma }c_zk_z\otimes e_z\bigg\|  \leq  C_{5.1}\sup _{z\in \Gamma }|c_z|.
$$
\medskip
\noindent
{\it Proof}.   Given such  $\Gamma $,   $\{e_z : z \in \Gamma \}$  and  $\{c_z : z \in \Gamma \}$,  define the operator  
$$
A  =  \sum _{z\in \Gamma }c_zk_z\otimes e_z.
$$
\noindent
Then  for every  $f \in L^2_a(\Omega )$   we have
$$
A^\ast f =  \sum  _{z\in \Gamma }\bar c_z\|K_z\|^{-1}f(z)e_z.
$$
\noindent 
From  Lemma 4.3  and (4.1) we obtain 
$$
\|A^\ast f\|^2 \leq C_{4.3}^2\sum _{z\in \Gamma }|c_z|^2\|K_z\|^{-2}|r(z)|^{-n-1}\|f\chi _{D(z,1)}\|^2  
\leq  (C_{4.3}^2/c)\sup _{z\in \Gamma }|c_z|^2\|f\|^2.
$$
\noindent 
Since   $f \in L^2_a(\Omega )$  is arbitrary,  this means that  $\|A\| = \|A^\ast \| \leq c^{-1/2}C_{4.3}\sup _{z\in \Gamma }|c_z|$.   $\square $
\medskip
\noindent
{\bf Lemma 5.2.}   {\it  Let}  $\Gamma $  {\it  be a}   $1$-{\it  separated set in}  $\Omega $.    {\it  Suppose that for every}  $z \in \Gamma $,  
{\it  we have a}   $\zeta (z) \in \Omega $   {\it  with}   $d(z,\zeta (z))< c_{4.5}$.   {\it  Then for every orthonormal set}   $\{e_z : z \in \Gamma \}$  
{\it  and  for every bounded set of complex coefficients}  $\{c_z : z \in \Gamma \}$,  {\it  we have}
$$
\bigg\|\sum _{z\in \Gamma }c_zk_z\otimes e_z  -  \sum _{z\in \Gamma }c_zk_{\zeta (z)}\otimes e_z \bigg\|  
\leq  C_{4.6}\sup _{z\in \Gamma }|c_z|d(z,\zeta (z)).
$$
\medskip
\noindent
{\it Proof}.    Write
$$
D  =  \sum _{z\in \Gamma }c_zk_z\otimes e_z  -  \sum _{z\in \Gamma }c_zk_{\zeta (z)}\otimes e_z  
=  \sum _{z\in \Gamma }c_z(k_z-k_{\zeta (z)})\otimes e_z.
$$
\noindent
For any  $f \in L^2_a(\Omega )$,  we have  
$$
D^\ast f  =  \sum _{z\in \Gamma }\bar c_z\langle f,k_z - k_{\zeta (z)}\rangle e_z.
$$
\noindent
Applying Proposition 4.6,  if  $d(z,\zeta (z))< c_{4.5}$  for every  $z \in \Gamma $,  then
$$
\align
\|D^\ast f\|^2  &=  \sum _{z\in \Gamma }|\bar c_z\langle f,k_z - k_{\zeta (z)}\rangle |^2  
\leq  C_{4.6}^2\sum _{z\in \Gamma }|c_z|^2d^2(z,\zeta (z))\|f\chi _{D(z,1)}\|^2  \\
&\leq   C_{4.6}^2\sup _{z\in \Gamma }|c_z|^2d^2(z,\zeta (z))\|f\|^2.
\endalign
$$
\noindent
Since   $f \in L^2_a(\Omega )$  is arbitrary,  this implies  $\|D\| = \|D^\ast \| \leq C_{4.6}\sup _{z\in \Gamma }|c_z|d(z,\zeta (z))$.      
$\square $
\medskip
\noindent
{\bf Corollary 5.3.}   {\it  Given any}  $a > 0$,  $0 \leq C < \infty $  {\it  and}  $\epsilon > 0$,  {\it  there is a}   $\delta > 0$  
{\it  such that the following estimate holds}:   {\it  Let}  $\Gamma $   {\it  be an} $a$-{\it  separated set in}  $\Omega $.    
{\it  Suppose that}    $\varphi $,  $\varphi '$,   $\psi $  {\it  and}  $\psi '$  {\it  are maps  from}  $\Gamma $ 
{\it  into}  $\Omega $.    {\it  If the inequalities}    
$$
d(z,\varphi (z)) \leq C,  \quad d(z,\psi (z)) \leq C,  \quad  d(\varphi (z),\varphi '(z)) \leq \delta ,  \quad 
d(\psi (z),\psi '(z)) \leq \delta 
$$
{\it  hold for every}  $z \in \Gamma $,  {\it  then for any  bounded set of coefficients}   $\{c_z : z \in \Gamma \}$  {\it  we have}  
$$
\bigg\|\sum _{z\in \Gamma }c_zk_{\varphi (z)}\otimes k_{\psi (z)}  
-  \sum _{z\in \Gamma }c_zk_{\varphi '(z)}\otimes k_{\psi '(z)}\bigg\|  
\leq  \epsilon \sup _{z\in \Gamma }|c_z|.
$$
\medskip
\noindent
{\it Proof}.    For  $z, w \in \Gamma $,  if   $d(z,w) >  2C + 2$,  then   $d(\varphi (z),\varphi (w)) > 2$   and  
$d(\psi (z),\psi (w)) > 2$.   By Lemma 2.11,  there is an   $N \in $  {\bf N}   determined by  $a$  and  $C$   such that  
$\Gamma $  admits a partition  
$$
\Gamma =  \Gamma _1\cup \cdots \cup \Gamma _N
$$
\noindent
with the property that for each  $j \in \{1,\dots ,N\}$,  the sets  $\{\varphi (z) : z \in \Gamma _j\}$  and  
$\{\psi (z) : z \in \Gamma _j\}$  are  $1$-separated.  Pick an orthonormal set  $\{e_z : z \in \Gamma \}$.   Fixing a  
$j \in \{1,\dots ,N\}$  for the moment,  we have
$$
\sum _{z\in \Gamma _j}c_zk_{\varphi (z)}\otimes k_{\psi (z)}  -  \sum _{z\in \Gamma _j}c_zk_{\varphi '(z)}\otimes k_{\psi '(z)}
=  AB^\ast - A{'B'}^\ast ,
$$
\noindent
where
$$
\align
A    &=  \sum _{z\in \Gamma _j}c_zk_{\varphi (z)}\otimes e_z,    \quad       B   =  \sum _{z\in \Gamma _j}k_{\psi (z)}\otimes e_z,   \\
A'   &=  \sum _{z\in \Gamma _j}c_zk_{\varphi '(z)}\otimes e_z,   \quad        B'  =  \sum _{z\in \Gamma _j}k_{\psi '(z)}\otimes e_z.
\endalign
$$
\noindent
We have
$$
AB^\ast - A{'B'}^\ast =  (A - A')B^\ast + A'(B^\ast - {B'}^\ast ).
$$
\noindent
Since  $\{\varphi (z) : z \in \Gamma _j\}$  and   $\{\psi (z) : z \in \Gamma _j\}$  are  $1$-separated,  
if we apply  Lemma 5.2 to $A - A'$  and  $B - B'$  and Lemma 5.1 to   $B$  and  $A'$,  we see that
$$
\|AB^\ast - A{'B'}^\ast\|  \leq  {\epsilon \over N}\sup _{z\in \Gamma }|c_z|
$$
\noindent
when  $\delta $  is sufficiently small.   This completes the proof.  $\square $

\bigskip
\centerline{{\bf 6.  Operators in the Toeplitz algebra}  ${\Cal T}$}
\medskip
Define the measure  
$$
d\tilde \mu (w)  =  K(w,w)dv(w)
$$
\noindent
on  $\Omega $.   By  (4.1),  this is just a slightly different version of 
the measure  $d\mu $  defined by  (2.10).  Given an  $f \in L^\infty (\Omega)$,  we have the integral representation 
$$
T_f  =    \int f(w)k_w\otimes k_wd\tilde \mu (w)
\tag 6.1
$$
\noindent
for the Toeplitz operator  $T_f$.  This formula is obtained by direct verification.   
Starting from this representation, we will  show that the Toeplitz algebra  ${\Cal T}$   contains certain classes of operators.    
The two main steps in the section are Propositions 6.4 and 6.6  below.
\medskip
\noindent
{\bf Proposition 6.1.}     {\it Suppose that}  $\Gamma $  {\it is a separated set in}   $\Omega $    
{\it and that}  $\{c_z : z \in \Gamma \}$  {\it  is a bounded set of  complex coefficients}.  
{\it Then the operator}
$$
\sum _{z\in \Gamma }c_zk_z\otimes k_z
$$   
\noindent
{\it  belongs to the closure of}   $\{T_f : f \in L^\infty (\Omega )\}$  {\it with respect to the operator norm}.
\medskip
\noindent
{\it Proof}.    By Lemma 2.11,   we may assume that  $\Gamma $  is  $1$-separated.     Let  $\epsilon > 0$  be given.
Since  $\sup _{z\in \Gamma }|c_z| < \infty $,  it follows from  Corollary 5.3 that there is a  $\delta > 0$  such that
$$
\bigg\|\sum _{z\in \Gamma }c_zk_z\otimes k_z - \sum _{z\in \Gamma }c_zk_{\zeta (z)}\otimes k_{\zeta (z)}\bigg\|   
\leq  \epsilon 
\tag 6.2
$$
\noindent
if  $\zeta (z) \in D(z,\delta )$   for every   $z \in \Gamma $.  We may, of course, assume that  $\delta < 1$,  
consequently   $D(z,\delta )\cap D(w,\delta ) = \emptyset $  for all  $z \neq w $  in  $\Gamma $.
\medskip
Define the function
$$
\varphi   =  \sum _{z\in \Gamma }{c_z\over \tilde \mu (D(z,\delta ))}\chi _{D(z,\delta )}
\tag 6.3
$$
\noindent
on  $\Omega $.   By   (4.1)  and   Proposition 2.6,  there is an  $\beta  > 0$  
such that   $\tilde \mu (D(z,\delta )) \geq \beta $  for every  $z \in \Omega $.  
Hence   $\varphi \in L^\infty (\Omega )$.   We will show that
$$
\bigg\|\sum _{z\in \Gamma }c_zk_z\otimes k_z - T_\varphi \bigg\|   \leq  \epsilon .
\tag 6.4
$$
\noindent
To prove  this, we define the measure   $d\nu _z = \{\tilde \mu (D(z,\delta ))\}^{-1}\chi _{D(z,\delta )}d\tilde \mu $   for every  $z \in \Gamma $.  
Then it follows from (6.1) and (6.3)  that 
$$
T_\varphi = \sum _{z\in \Gamma }c_z\int k_w\otimes k_wd\nu _z(w).
$$
\noindent
Note that each  $d\nu _z$  is a probability measure concentrated on  $D(z,\delta )$.   Hence   $d\nu _z $  is in the weak-* closure of  
convex combinations of unit point masses on    $D(z,\delta )$.   Therefore  $T_\varphi $  is the limit in weak operator topology of operators 
of the form
$$
T  =  {1\over k}\sum  _{j=1}^k\sum _{z\in \Gamma }c_zk_{\zeta (z;j)}\otimes k_{\zeta (z;j)},
$$
\noindent
where  $k \in $  {\bf N}  and   $\zeta (z;j) \in D(z,\delta )$  for all  $z \in \Gamma $  and  $j \in \{1,\dots ,k\}$.  By  (6.2),
$$
\bigg\|\sum _{z\in \Gamma }c_zk_z\otimes k_z - T\bigg\|     \leq  
{1\over k}\sum  _{j=1}^k\bigg\|\sum _{z\in \Gamma }c_zk_z\otimes k_z - \sum _{z\in \Gamma }c_zk_{\zeta (z;j)}\otimes k_{\zeta (z;j)}\bigg\|
\leq \epsilon .
$$
\noindent
Since this holds for every such  $T$  and since  $T_\varphi $  is the weak limit of such  $T$'s,  (6.4) follows.  
Since  $\epsilon > 0$  is arbitrary,  this completes the proof.   $\square $
\medskip
Next we remind the reader of a well-known fact:
\medskip
\noindent
{\bf Proposition 6.2.}   [22,Theorem 4.1.25]   {\it The  Toeplitz algebra}  ${\Cal T}$ {\it  contains}  ${\Cal K}$,  
{\it the collection of compact operators on the Bergman space}  $L^2_a(\Omega )$.  
\medskip
\noindent
{\bf Definition 6.3.}     (a)   Let  ${\Cal D}_0$   denote the collection of operators of the form
$$
\sum _{z\in \Gamma }c_zk_z\otimes k_{\gamma (z)},
$$  
\noindent
where   $\Gamma $  is any separated set in   $\Omega $,    $\{c_z : z \in \Gamma \}$  is any bounded set of complex 
coefficients,  and  $\gamma :  \Gamma \rightarrow \Omega $   is any map for which there is a   $0 \leq C < \infty $   
such that
$$
d(z,\gamma (z))  \leq   C
\tag 6.5
$$
\noindent
for every    $z \in \Gamma $.   

\noindent
(b)   Let   ${\Cal D}$   denote the closure of the linear span of   ${\Cal D}_0$ with respect to the operator norm.

\noindent
(c)  For any  $A \in {\Cal B}(L^2_a(\Omega ))$,    ${\Cal D}_0(A)$  denotes  the collection of operators of the form  
$$
\sum _{z\in \Gamma }c_z\langle Ak_{\psi (z)},k_{\varphi (z)}\rangle k_{\varphi (z)}\otimes k_{\psi (z)},
$$  
\noindent
where  $\Gamma $ is a separated set in  $\Omega $,   $\{c_z : z \in \Gamma \}$  is a bounded set of coefficients,   
and  $\varphi , \psi  :  \Gamma \rightarrow \Omega $  are maps for which there is a  $0 \leq C < \infty $  
such that    $d(z,\varphi (z)) \leq C$  and    $d(z,\psi (z)) \leq C$  for every   $z \in \Gamma $.    

\noindent
(d)  For any  $A \in {\Cal B}(L^2_a(\Omega ))$,     ${\Cal D}(A)$  denotes    
the closure of the linear span of   ${\Cal D}_0(A)$ with respect to the operator norm.
\medskip
\noindent
{\bf Proposition 6.4.}   {\it We have the inclusion}  ${\Cal D}_0 \subset {\Cal T}$.  {\it Consequently},  
${\Cal D} \subset {\Cal T}$.
\medskip
\noindent
{\it Proof}.     Let  $\Gamma $, $\{c_z : z \in \Gamma \}$,  $\gamma $  and  $C$   be as described in Definition 6.3(a),  and consider   
$$
T  =  \sum _{z\in \Gamma }c_zk_z\otimes k_{\gamma (z)}.
$$
\noindent
To show that  $T \in {\Cal T}$,   by Lemma 2.11, we may assume that  
$$
d(z,w) > 4C + 2  \quad \text{for all} \ \  z \neq w  \ \  \text{in}  \ \   \Gamma .   
\tag 6.6
$$
\noindent
For each  $z \in \Gamma $,  since  $d(z,\gamma (z)) \leq C$,  by (2.4)  there is a 
$C^1$  map  $g_z : [0,1] \rightarrow \Omega $  such that   $g_z(0) = z$,  $g_z(1) = \gamma (z)$,   and such that the number
$$
\ell _z  =  \int _0^1\sqrt{\langle {\Cal B}(g_z(t))g_z'(t),g_z'(t)\rangle }dt
$$
\noindent
satisfies the condition  $\ell _z \leq 2C$.    Pick a  $k \in $  {\bf N}  such that  $2C/k < \min \{a_0,c_{4.7}\}$,  
where  $a_0$  and  $c_{4.7}$  are the constants  in Lemmas 3.10 and  4.7 respectively.   For each  $z \in \Gamma $,  there are  
$$
0 = x(z,0) \leq x(z,1) \leq \cdots \leq x(z,k-1) \leq x(z,k) = 1
$$  
\noindent
such that  
$$
\int _0^{x(z,j)}\sqrt{\langle {\Cal B}(g_z(t))g_z'(t),g_z'(t)\rangle }dt  =  {j\over k}\ell _z
$$
\noindent
for  $j= 0,1,\dots ,k$.  For each pair of  $z \in \Gamma $  and  $j \in \{0,1,\dots ,k\}$,  we now define
$$
\gamma _j(z)  =  g_z(x(z,j)).
$$
\noindent
We have  $\gamma _0(z) = z$   and  $\gamma _k(z) = \gamma (z)$,   $z \in \Gamma $.   
Since   $\ell _z \leq  2C$,    for  all  $0 \leq j < k$  and  $z \in \Gamma $,  
$$
d(\gamma _j(z),\gamma _{j+1}(z))  \leq  \int _{x(z,j)}^{x(z,j+1)}\sqrt{\langle {\Cal B}(g_z(t))g_z'(t),g_z'(t)\rangle }dt  
=  \ell _z/k  <  \min \{a_0,c_{4.7}\}.
$$
\noindent
By Lemma 4.7,  this ensures that
$$
|\langle k_{\gamma _j(z)},k_{\gamma _{j+1}(z)}\rangle |  \geq 1/2  
\tag 6.7
$$
\noindent
for  all  $0 \leq j < k$  and  $z \in \Gamma $.
\medskip
To prove that   $T \in {\Cal T}$,   it suffices to show that for every  $j \in \{0,1,\dots ,k\}$   
and every subset  $E$   of  $\Gamma $,  we have 
$$
\sum _{z\in E}c_zk_z\otimes k_{\gamma _j(z)}  \in  {\Cal T}.
\tag 6.8  
$$
\noindent
We will accomplish this by an induction on  $j$.   Since   $\gamma _0(z) = z$  for every  $z \in \Gamma $,   the case  
$j = 0$  follows from Proposition 6.1.   Suppose now that  $0 \leq j < k$   and that (6.8) holds for this  $j$  and for every  
$E \subset \Gamma $.  To simplify notation,  for every   $S \subset \Gamma $,    let us denote 
$$
X_S  =  \sum _{z\in S}c_zk_z\otimes k_{\gamma _j(z)}   \quad  \text{and}  \quad  
Y_S  =  \sum _{z\in S}{1\over \langle k_{\gamma _{j+1}(z)},k_{\gamma _j(z)}\rangle }k_{\gamma _{j+1}(z)}\otimes k_{\gamma _{j+1}(z)}.
$$ 
\noindent
By the induction hypothesis,  we have  $X_S \in {\Cal T}$.  By (6.7)  and Proposition 6.1,  we also have  $Y_S \in {\Cal T}$.  Therefore  
$X_SY_S \in {\Cal T}$    for every  $S \subset \Gamma $.  To complete the induction,  it suffices to show that 
given  $E \subset \Gamma $  and  $\epsilon > 0$,  there is a finite partition   $E = S_1\cup \cdots \cup S_N$   such that
$$
\bigg\|X_{S_1}Y_{S_1} + \cdots +  X_{S_N}Y_{S_N} - \sum _{z\in E}c_zk_z\otimes k_{\gamma _{j+1}(z)}\bigg\|   \leq  \epsilon .
\tag 6.9
$$ 
\noindent
To see how this is done, first note that for any partition   $E = S_1\cup \cdots \cup S_N$,
$$
\align
X_{S_1}Y_{S_1} &+ \cdots +  X_{S_N}Y_{S_N} - \sum _{z\in E}c_zk_z\otimes k_{\gamma _{j+1}(z)}\\
&=  \sum _{\nu =1}^N\sum \Sb z,w\in S_\nu \\ z\neq w \endSb 
c_z{ \langle k_{\gamma _{j+1}(w)},k_{\gamma _j(z)}\rangle \over  \langle k_{\gamma _{j+1}(w)},k_{\gamma _j(w)}\rangle }
k_z\otimes k_{\gamma _{j+1}(w)}  =  UWV^\ast  ,
\endalign
$$
\noindent
where  
$$
\align
U  &=  \sum _{z\in E}c_zk_z\otimes e_z,   \quad  V  =  \sum _{z\in E}k_{\gamma _{j+1}(z)}\otimes e_z
\quad  \text{and}   \\  
W  &=  \sum _{\nu =1}^N\sum \Sb z,w\in S_\nu \\ z\neq w \endSb 
{\langle k_{\gamma _{j+1}(w)},k_{\gamma _j(z)}\rangle \over  \langle k_{\gamma _{j+1}(w)},k_{\gamma _j(w)}\rangle }
e_z\otimes e_w,
\endalign
$$
\noindent
where  $\{e_z : z \in E\}$  is an orthonormal set.
\medskip
By (6.6),  $\{\gamma _{j+1}(z) : z \in E\}$  is a  $1$-separated set.  Thus Lemma 5.1 provides the bound    
$\|V\| \leq  C_{5.1}$.  Similarly,  $\|U\| \leq C_{5.1}c$,  where  $c = \sup _{z\in \Gamma }|c_z|$.   Consequently  
$$
\bigg\|X_{S_1}Y_{S_1} + \cdots +  X_{S_N}Y_{S_N} - \sum _{z\in E}c_zk_z\otimes k_{\gamma _{j+1}(z)}\bigg\|   \leq C_{5.1}^2c\|W\|.
\tag 6.10
$$ 
\noindent
Thus we need to find a partition  $E = S_1\cup \cdots \cup S_N$  such that  $\|W\|$  is small.   To do this,  consider an  $R > 3C_{3.7} + 1$,
whose value will be determined below.  By Lemma 2.11,  there is a partition  
$E = S_1\cup \cdots \cup S_N$  such that for every  $\nu \in \{1,\dots ,N\}$,  the conditions  $z, w \in S_\nu $  and  
$z \neq w$  imply  $d(z,w) > R$.   With    $S_1, \dots ,S_N$  so chosen,  we define
$$
{\Cal F} = \bigcup _{\nu =1}^N\{(z,w) \in S_\nu \times S_\nu :  z \neq w\}.
$$
\noindent
We can rewrite  $W$  in the form
$$
W  =  \sum _{(z,w)\in E\times E}a(z,w)e_z\otimes e_w,  
$$
\noindent
where
$$
a(z,w)  =  
\left\{
\matrix
{\langle k_{\gamma _{j+1}(w)},k_{\gamma _j(z)}\rangle \over  \langle k_{\gamma _{j+1}(w)},k_{\gamma _j(w)}\rangle }   &\text{if}  &(z,w) \in {\Cal F}  \\
\ \  \\
0   &\text{if}  &(z,w) \notin {\Cal F}
\endmatrix
\right.  .
$$
\noindent
Recall that  $d(\gamma _p(z),\gamma _{p+1}(z)) < a_0$  for all $z \in \Gamma $  and  $0 \leq p < k$.   
Recalling (6.7)  and applying  Theorem 4.1  and  Lemmas 3.1, 2.1,  and Lemma 3.10 multiple times,    we obtain
$$
\align
|a(z,w)|  \leq  2{|K(\gamma _j(z),\gamma _{j+1}(w))|\over \|K_{\gamma _{j+1}(w)}\|\|K_{\gamma _j(z)}\|}      
&\leq C_1\bigg({|r(\gamma _j(z))|^{1/2}|r(\gamma _{j+1}(w))|^{1/2}\over  F(\gamma _j(z),\gamma _{j+1}(w))}\bigg)^{n+1}   \\
&\leq C_2\bigg({|r(z)|^{1/2}|r(w)|^{1/2}\over  F(z,w)}\bigg)^{n+1} 
\endalign
$$
\noindent
for  $(z,w) \in {\Cal F}$.   Pick an  $\eta \in (0,1/2)$  and define   $h(w) = |r(w)|^{(n/2)+\eta }$,    $w \in \Gamma$.   
If  $(z,w) \in {\Cal F}$,  then  $d(z,w) > R$  by design.
Since  $E$  is  $1$-separated,  it follows from Lemma 4.2 that   
$$
\sum _{w\in E}|a(z,w)|h(w) \leq  C_2\sum _{w\in E\backslash D(z,R)} 
|r(w)|^{(n/2)+\eta }\bigg({|r(z)|^{1/2}|r(w)|^{1/2}\over  F(z,w)}\bigg)^{n+1}   
\leq  {C_2C_{4.2}\over 2^{sR}}h(z)
$$
\noindent 
for every  $z \in E$.  A similar inequality holds for  $\sum _{z\in E}|a(z,w)|h(z)$,  $w \in E$.   By the standard Schur test,  we conclude that  
$\|W\| \leq  C_2C_{4.2}2^{-sR}$.   Recalling (6.10),  we see that (6.9) holds if we pick  $R > 3C_{3.7} + 1$  such that  
$C_{5.1}^2cC_2C_{4.2}2^{-sR} \leq \epsilon $.  This completes the proof.   $\square $
\medskip
Following the ideas in  [28],  we will now generalize the notion of  {\it localized operators}  to strongly pseudo-convex domains.
\medskip
\noindent
{\bf Definition 6.5.}    Let  $A$   be a bounded operator on the Bergman space  $L_a^2(\Omega )$.   Then   LOC$(A)$   
denotes the collection of operators of the form   
$$
T  =  \sum _{z\in \Gamma }T_{f_z}AT_{f_z},   
\tag 6.11
$$
\noindent
where  $\Gamma $ is any separated set in   $\Omega $ and 
$\{f_z : z \in \Gamma \}$   is any family of continuous functions on   $\Omega $  satisfying the following three conditions:

\noindent
(1)    There is a   $0 < \rho  < \infty $   such that    $f_z = 0$   on  $\Omega \backslash D(z,\rho )$   for every   $z \in \Gamma $.  

\noindent
(2)  The inequality  $0 \leq f_z \leq 1$      holds on   $\Omega $  for every   $z \in \Gamma $.

\noindent
(3)    The family  $\{f_z : z \in \Gamma \}$  satisfies a uniform Lipschitz condition on  $\Omega $  with  respect to the metric  $d$.
That is, there is a  $0 < C < \infty $   such that     $|f_z(\zeta ) - f_z(\xi )| \leq Cd(\zeta , \xi )$   
for all  $z \in \Gamma $  and $\zeta , \xi \in \Omega $.     
\medskip
\noindent
{\bf Proposition 6.6.}   {\it For every bounded operator}   $A$  {\it on}   $L_a^2(\Omega )$,  {\it we have} LOC$(A)  \subset {\Cal D}(A)$.
\medskip
\noindent
{\it Proof}.    Let  $A$   be a bounded operator  on  $L_a^2(\Omega )$,  and consider a  $T$  given by (6.11).  
To prove that  $T \in {\Cal D}(A)$,  by Lemma 2.11,  we may assume that  $\Gamma $  is  $1$-separated.  
For convenience,  let us define the product measure   $\nu = \tilde \mu \times \tilde \mu $  on  $\Omega \times \Omega $.
By (6.1),   for each  $z \in \Gamma $  we have
$$
T_{f_z}AT_{f_z}  =  \iint h_z(u,v)k_u\otimes k_vd\nu (u,v),
\tag 6.12
$$
\noindent
where
$$
h_z(u,v)  =  f_z(u)f_z(v)\langle Ak_v,k_u\rangle .  
\tag 6.13
$$
\noindent
By condition (1) above,  $h_z$  vanishes on the complement of  $D(z,\rho )\times D(z,\rho )$.  
It follows from  Proposition 4.6 and condition (3) above that for any   $a > 0$,  there is a  $b  > 0$  
such that
$$
\sup _{z\in \Gamma }|h_z(u,v) - h_z(u',v')|  \leq a \quad  \text{if} \ \  
d(u,u') \leq b  \ \ \text{and}  \ \  d(v,v') \leq b.
\tag 6.14
$$
\noindent
The rest of the proof  is divided into two steps.
\medskip
Step I.  We first show that for any   $\epsilon > 0$,  there is a  $0 < \delta \leq \rho $  such that the following holds true:  
Suppose that  $\Lambda $ is a subset of  $\Gamma $.  For each  $z \in  \Lambda $,  let   $\varphi (z),  \psi (z) \in  D(z,\rho )$.   
For each  $z \in \Lambda $,   suppose that we have a Borel set   $E_z = F_z\times G_z$   with    
$F_z$  $\subset $  $D(\varphi (z),\delta )$,  $G_z \subset  D(\psi (z),\delta )$ and  $\nu (E_z) > 0$.   
Finally,  for  each  $z \in \Lambda $,  let $a_z \in [0,2]$.  Then
$$
\bigg\|\sum _{z\in \Lambda }{a_z\over \nu (E_z)}\iint _{E_z}h_z(u,v)k_u\otimes k_vd\nu (u,v) - 
\sum _{z\in \Lambda  }a_zh_z(\varphi (z),\psi (z))k_{\varphi (z)}\otimes k_{\psi (z)}\bigg\|  \leq  \epsilon .
\tag 6.15
$$
\noindent
To prove this,  denote
$$
\align
W  &=  \sum _{z\in \Lambda }{a_z\over \nu (E_z)}\iint _{E_z}h_z(u,v)k_u\otimes k_vd\nu (u,v)   \quad  \text{and}   \\
Z   &=  \sum _{z\in \Lambda }a_zh_z(\varphi (z),\psi (z))k_{\varphi (z)}\otimes k_{\psi (z)}.
\endalign
$$
\noindent
Note that for each   $z \in \Lambda $,  $\{\chi _{E_z}/\nu (E_z)\}d\nu $  is a probability measure concentrated on  $E_z$.   
Thus it is in the weak-* closure of convex combinations of unit point masses on  $E_z$.   Consequently  $W$  is in the 
closure in weak operator topology of operators of the form
$$
W'  =  {1\over k}\sum _{j=1}^k\sum _{z\in \Lambda }a_zh_z(u(z;j),v(z;j))k_{u(z;j)}\otimes k_{v(z;j)},
$$
\noindent
where   $k \in $  {\bf N}  and,  for each  $1 \leq j \leq k$,  we have   $(u(z;j),v(z;j)) \in E_z$,  i.e.,  $u(z;j) \in F_z$  and  
$v(z;j) \in G_z$,  $z \in \Lambda $.  It is easy to see that
$$
W' - Z   =  {1\over k}\sum _{j=1}^k(X_j + Y_j),
$$
\noindent
where
$$
\align 
X_j  &=  \sum _{z\in \Lambda }a_z\{h_z(u(z;j),v(z;j)) - h_z(\varphi (z),\psi (z))\}k_{u(z;j)}\otimes k_{v(z;j)}  \quad  \text{and}   \\
Y_j  &=  \sum _{z\in \Lambda }a_zh_z(\varphi (z),\psi (z))k_{u(z;j)}\otimes k_{v(z;j)} 
- \sum _{z\in \Lambda }a_zh_z(\varphi (z),\psi (z))k_{\varphi (z)}\otimes k_{\psi (z)}.
\endalign
$$
\noindent
From (6.14) and Lemmas 2.11 and 5.1 we see that there is a  $\delta _ 1 > 0$  such that   $\|X_j\| \leq \epsilon /2$  
for every  $1 \leq j \leq k$   if   $\delta \leq \delta _1$.   By Corollary 5.3,  there is a  
$\delta _ 2 > 0$  such that   $\|Y_j\| \leq \epsilon /2$   for every  $1 \leq j \leq k$  if    $\delta \leq \delta _2$.   
Hence for any   $0 < \delta \leq \min \{\delta _1,\delta _2,\rho \}$,  we have   $\|W' - Z\| \leq \epsilon $.   
Since    $W - Z$  is the weak limit of operators of the form  
$W' - Z$,  we have  $\|W - Z\| \leq \epsilon $   for any choice of  $0 < \delta \leq \min \{\delta _1,\delta _2,\rho \}$.  
This proves (6.15)  and completes Step I.
\medskip
Step II.  Recall that  $\nu = \tilde \mu \times \tilde \mu $.   By  (4.1)  and Proposition 2.6,  there is an  
$N \in $  {\bf N}  such that   $N \geq \nu (D(w,2\rho )\times D(w,2\rho ))$  for every   $w \in \Omega $.   Let  $\epsilon > 0$  
be given.  We will now find a   $B \in \text{span}({\Cal D}_0(A))$  such that
$$
\|T - B\|  \leq  N\epsilon .
\tag 6.16
$$
\noindent
Since   $\epsilon > 0$  is arbitrary,  this will imply the membership  $T \in {\Cal D}(A)$.  To find such a  
$B \in \text{span}({\Cal D}_0(A))$,   let  $\delta $  be the number provided for this  $\epsilon $  in Step I.   
For each  $z \in \Gamma $,  there is a subset    $S_z$  in   $D(z,\rho )$  that is maximal with respect to the property
$$
D(x,\delta /2)\cap D(y,\delta /2)  =  \emptyset  \quad  \text{for all}  \ \  x \neq y \ \  \text{in} \ \  S_z.
$$
\noindent
By Proposition 2.6  and the fact that  $\mu (D(z,2\rho )) < \infty $,  we see that   $S_z$  is a finite set,  and consequently we can 
represent it in the form   $S_z = \{\varphi (z;1),\dots ,\varphi (z;m(z))\}$  with some  $m(z) \in $  {\bf N}.  The maximality of  
$S_z$  implies that   $\cup _{j=1}^{m(z)}D(\varphi (z;j),\delta ) \supset D(z,\rho )$.   Thus for each  $z \in \Gamma $,  
a standard set-theoretical argument gives us Borel sets  
$$
F(z;1), \dots , F(z;m(z))
$$
\noindent
with the following properties:

(i)  $D(\varphi (z;j),\delta /2) \subset F(z;j)  \subset D(\varphi (z;j),\delta )$  for each  $j \in \{1,\dots ,m(z)\}$.

(ii)  $F(z;i)\cap F(z;j) = \emptyset $  for all  $i \neq j$  in  $\{1,\dots ,m(z)\}$.

(iii)   $D(z,\rho )  \subset \cup _{j=1}^{m(z)}F(z;j) \subset D(z,2\rho )$.

\noindent
We now define  $E_{z;i,j} = F(z;i)\times F(z;j)$  for  $z \in \Gamma $  and  $i, j \in \{1,\dots ,m(z)\}$.

\medskip
Let   $z \in \Gamma $.  Since  $h_z$  vanishes on   $(\Omega \times \Omega )\backslash (D(z,\rho )\times D(z,\rho ))$,  (iii) and (ii) imply
$$
T_{f_z}AT_{f_z}  =  \sum _{i=1}^{m(z)} \sum _{j=1}^{m(z)} \iint _{E_{z;i,j}}h_z(u,v)k_u\otimes k_vd\nu (u,v).
\tag 6.17
$$
\noindent
By (i) and Proposition 2.6,  there is a  $k \in $  {\bf N}   such that   $1/k  < \nu (E_{z;i,j})$   for all  $z \in \Gamma $  and  $i, j \in \{1,\dots ,m(z)\}$.  
For such a triple of  $z, i, j$,  we let   $p(z;i,j)$  be the largest natural number satisfying the condition   $p(z;i,j)/k \leq \nu (E_{z;i,j})$.    Define 
$$
a(z;i,j) = {k\over p(z;i,j)}\nu (E_{z;i,j})
$$
\noindent
for   $z \in \Gamma $  and  $i, j \in \{1,\dots ,m(z)\}$.   Then   $0 <  a(z;i,j)  \leq 2$,  because the definition of  $p(z;i,j)$  ensures that  
$\{p(z;i,j)+1\}/k > \nu (E_{z;i,j})$.    We can now rewrite  (6.17)  in the form
$$
T_{f_z}AT_{f_z}  =  {1\over k}\sum _{i=1}^{m(z)} \sum _{j=1}^{m(z)}p(z;i,j){a(z;i,j)\over \nu (E_{z;i,j})}
\iint _{E_{z;i,j}}h_z(u,v)k_u\otimes k_vd\nu (u,v).
\tag 6.18
$$
\noindent
On the other hand,  for every   $z \in \Gamma $,  we have
$$
\align
\sum _{i=1}^{m(z)} \sum _{j=1}^{m(z)}p(z;i,j)  &=  k\sum _{i=1}^{m(z)} \sum _{j=1}^{m(z)}{p(z;i,j)\over k}  
\leq  k\sum _{i=1}^{m(z)} \sum _{j=1}^{m(z)}\nu (E_{z;i,j})   \\
&\leq k\nu (D(z,2\rho )\times D(z,2\rho ))  \leq  kN.
\endalign
$$
\noindent
We can regard  $p(z;i,j)$  as the  ``multiplicity"  with which the triple  $(z,i,j)$  appears in  (6.18).  
The above estimate shows that for a fixed  $z \in \Gamma $,  
all the multiplicities add up to something less than or equal to  $kN$.
Thus there are subsets  $\Gamma _1, \Gamma _2 , \dots , \Gamma _{kN}$  of  $\Gamma $  such that
$$
\sum _{z\in \Gamma }T_{f_z}AT_{f_z}  
=  {1\over k}\sum _{\ell =1}^{kN}\sum _{z\in \Gamma _\ell }{a(z;i(z,\ell ),j(z,\ell ))\over \nu (E_{z;i(z,\ell ),j(z,\ell )})}
\iint _{E_{z;i(z,\ell ),j(z,\ell )}}h_z(u,v)k_u\otimes k_vd\nu (u,v),
\tag 6.19
$$
\noindent
where for each pair of  $\ell \in \{1,\dots ,kN\}$   and   $z \in \Gamma _\ell $  we have  $i(z,\ell ), j(z,\ell ) \in \{1,\dots ,m(z)\}$.   
For each  $\ell \in \{1,\dots ,kN\}$,   define
$$
B_\ell  = \sum _{z\in \Gamma _\ell }a(z;i(z,\ell ),j(z,\ell ))h_z(\varphi (z;i(z,\ell )),\varphi (z;j(z,\ell )))k_{\varphi (z;i(z,\ell ))}\otimes k_{\varphi (z;j(z,\ell ))}.
$$
\noindent
Since  $\varphi (z;i(z,\ell )), \varphi (z;j(z,\ell )) \in D(z,\rho )$  for every  $z \in \Gamma _\ell $,   recalling  Definition 6.3(c)  and  
(6.13),   we have  $B_\ell \in {\Cal D}_0(A)$.  Therefore the operator  
$$
B  =  {1\over k}\sum _{\ell =1}^{kN}B_\ell 
\tag 6.20
$$
\noindent
belongs to the linear span of  ${\Cal D}_0(A)$.    By the choice of  $\delta $  and  Step I,  we have
$$
\bigg\|\sum _{z\in \Gamma _\ell }{a(z;i(z,\ell ),j(z,\ell ))\over \nu (E_{z;i(z,\ell ),j(z,\ell )})}
\iint _{E_{z;i(z,\ell ),j(z,\ell )}}h_z(u,v)k_u\otimes k_vd\nu (u,v) - B_\ell \bigg\|  \leq  \epsilon ,
$$
\noindent
$1 \leq \ell \leq  kN$.  Combining this with (6.19) and (6.20),  we see that  $\|T - B\|$  does not exceed
$$
 {1\over k}\sum _{\ell =1}^{kN}\bigg\|\sum _{z\in \Gamma _\ell }{a(z;i(z,\ell ),j(z,\ell ))\over \nu (E_{z;i(z,\ell ),j(z,\ell )})}
\iint _{E_{z;i(z,\ell ),j(z,\ell )}}h_z(u,v)k_u\otimes k_vd\nu (u,v) - B_\ell \bigg\|  \leq  N\epsilon .
$$
\noindent
This proves (6.16)  and completes the proof of the proposition.   $\square $
\medskip
It follows from Lemma 2.11 that   ${\Cal D}(A) \subset {\Cal D}$  for every  $A \in {\Cal B}(L^2_a(\Omega ))$  
(cf. Definition 6.3).     Thus from  Propositions 6.6 and 6.4   we immediately obtain
\medskip
\noindent
{\bf Corollary 6.7.}   {\it For every bounded operator}   $A$  {\it on}   $L_a^2(\Omega )$,  {\it we have} LOC$(A)  \subset {\Cal T}$.
\medskip
To conclude this section, we recall  
\medskip
\noindent
{\bf  Lemma 6.8.}    {\it Let}   $\{f_1, \dots ,f_\ell \}$  {\it be a finite set of functions in}  $L^\infty (\Omega )$  {\it with the property that}
$f_jf_k = 0$  {\it for all}   $j \neq k$  {\it in}  $\{1,\dots ,\ell \}$.   {\it Let}  $A$  {\it be any bounded operator on the Bergman space}  
 $L_a^2(\Omega )$.   {\it Then there exist complex numbers}   $\{\gamma _1,\dots ,\gamma _\ell \}$   {\it with}   $|\gamma _k| = 1$   
 {\it for every}  $k \in \{1,\dots ,\ell \}$  {\it and a subset}  $E$  of  $\{1,\dots ,\ell \}$   {\it such that if we define}
 $$
 F = \sum _{k\in E}f_k,  \quad   G = \sum _{k\in \{1,\dots ,\ell \}\backslash E}f_k,  \quad 
 F' = \sum _{k\in E}\gamma _kf_k   \quad \text{\it and}  \ \  G' = \sum _{k\in \{1,\dots ,\ell \}\backslash E}\gamma _kf_k, 
 $$
 \noindent
 {\it then}
 $$
 \bigg\|\sum _{j\neq k}T_{f_j}AT_{f_k}\bigg\|   \leq  4(\|T_{F'}AT_G\| + \|T_{G'}AT_F\|).
 $$
\medskip
This lemma was proved in the case of the unit ball as Lemma 5.1 in [28].  But the proof in the case of a general  $\Omega $ 
is exactly the same.  The only property of Toeplitz operators that was used in the proof of [28,Lemma 5.1] was that a Toeplitz operator 
is the compression to a subspace of a multiplication operator on an  $L^2$.   Thus not only does Lemma 6.8 hold, its analogue also holds,  
for example, in the setting of Hardy spaces.  For that reason we will not repeat the proof of Lemma 6.8 here.

\bigskip
\centerline{\bf 7.  Oscillation and compactness}  
\medskip
For a continuous function  $f$  on  $\Omega $,  we define
$$
\text{diff}(f)  =  \sup \{|f(z) - f(w)| :  d(z,w) \leq 1\}.
$$
\medskip
\noindent
{\bf Lemma 7.1.}   {\it  For any continuous function}  $f$  {\it  on}  $\Omega $  {\it  and any}   $k \in $  {\bf N},  {\it  we have}
$$
|f(z) - f(w)| \leq  (k+1)\text{diff}(f) 
\tag 7.1
$$
\noindent
{\it  for any pair of}  $z, w \in \Omega $  {\it  satisfying the condition}  $d(z,w) \leq k$.
\medskip
\noindent
{\it Proof}.    Let   $z, w \in \Omega $  be such that  $d(z,w) \leq k$.   By (2.4),  there is a  $C^1$  map  
$\gamma : [0,1] \rightarrow \Omega $   such that  $\gamma (0) = z$,  $\gamma (1) = w$   and
$$
\int _0^1\sqrt{\langle {\Cal B}(\gamma (t))\gamma '(t),\gamma '(t)\rangle }dt   \leq  k+1.
\tag 7.2
$$
\noindent
There are   $0 = x_0 \leq x_1 \leq  \cdots \leq x_k \leq x_{k+1} = 1$  such that
$$
\int _{x_j}^{x_{j+1}}\sqrt{\langle {\Cal B}(\gamma (t))\gamma '(t),\gamma '(t)\rangle }dt   =
{1\over k+1}\int _0^1\sqrt{\langle {\Cal B}(\gamma (t))\gamma '(t),\gamma '(t)\rangle }dt
\tag 7.3
$$
\noindent
for  $0 \leq j \leq k$.   Define  $z_j = \gamma (x_j)$,  $j = 0, 1, \dots ,k+1$.   Then  $z_0 = z$  and  $z_{k+1} = w$.  
It follows from (7.2), (7.3) and (2.4) that   $d(z_j,z_{j+1}) \leq 1$,  consequently  
$$
|f(z_j) - f(z_{j+1})|  \leq  \text{diff}(f)
$$
\noindent
for  $j = 0,1, \dots ,k$.  With this inequality, (7.1) follows from a standard telescoping sum.   $\square $
\medskip
Recall that  $P$  denotes the orthogonal projection from  $L^2(\Omega )$  onto  $L^2_a(\Omega )$.
\medskip
\noindent
{\bf Lemma 7.2.}  {\it  There is a constant}  $0 < C_{7.2} < \infty $  {\it  such that}   $\|[M_f,P]\| \leq  C_{7.2}\text{diff}(f)$  
{\it  for every bounded continuous function}   $f$   {\it  on}   $\Omega $.
\medskip
\noindent
{\it Proof}.     Let  $T$  be the integral operator on  $L^2(\Omega )$  with the function 
$$
\{d(z,w) + 2\}|K(z,w)|
$$
\noindent
as its integral kernel.  We know that  $|K(z,w)| = |K(w,z)|$.    Recall that   the Bergman kernel  $K$  is  bounded on
$(\Omega \times \Omega )\backslash {\Cal R}_\delta $  for any  $\delta > 0$  [14].   Thus it follows from 
Theorem 4.1 and  Lemma 3.1 that   $|K(z,w)| \leq C_1F(z,w)^{-n-1}$   for all   $z, w \in \Omega $.   Combining this fact with 
Lemmas 3.9 and 3.2, and with the Schur test,  we see that the operator  $T$  is bounded on   $L^2(\Omega )$.  Let  
$f$  be a  bounded continuous function on  $\Omega $.   It follows from Lemma 7.1 that
$$
|(f(z) - f(w))K(z,w)|  \leq  \text{diff}(f)\{d(z,w) + 2\}|K(z,w)|
$$
\noindent
 for all   $z, w \in \Omega $.   Hence   $\|[M_f,P]\| \leq  \text{diff}(f)\|T\|$.    $\square $
\medskip
Recall that a continuous function  $f$  on  $\Omega $  is said to have {\it vanishing oscillation} if 
$$
\lim _{z\rightarrow \partial \Omega }\sup \{|f(z) - f(w)| : d(z,w) \leq 1\}  =  0.
$$
\noindent
We denote by  $\text{VO}_{\text{bdd}}$   the collection of continuous functions of vanishing oscillation on  $\Omega $  
that are also  {\it bounded}.   
\medskip
\noindent
{\bf Proposition 7.3.}   {\it For each}  $f \in \text{VO}_{\text{bdd}}$,  {\it the commutator}  $[M_f,P]$  {\it is compact}.
\medskip
\noindent
{\it Proof}.  It suffices to consider  $f \in \text{VO}_{\text{bdd}}$   with  $\|f\|_\infty \leq 1$.   For each  $R > 0$,  we will decompose  
$f$  in the form  $f = g_R + h_R$,  where  $g_R$  has a compact support and  $h_R$   satisfies the conditions  
$\text{diff}(h_R) \leq R^{-1}$   and   $\|h_R\|_\infty \leq 1$.     Since  $g_R$  has a compact support,   $[M_{g_R},P]$  is compact.
On the other hand,   Lemma 7.2  tells us that   $\|[M_{h_R},P]\| \rightarrow 0$  as  $R \rightarrow \infty $.    
Thus such a general decomposition implies the compactness of  $[M_f,P]$.
\medskip
To decompose  $f$,  let  $R > 0$  be given.  Since  $f \in \text{VO}_{\text{bdd}}$,  there is 
a  $t > 0$  such that
$$
|f(z) - f(w)|  \leq  (2R)^{-1}   \quad  \text{if} \ \  z \in H_t \ \  \text{and}  \ \  d(z,w) \leq 1,
\tag 7.4
$$
\noindent
where    $H_t = \{\zeta \in \Omega : -r(\zeta ) < t\}$,   and we may assume  $H_t \neq \Omega $.    Define
$$
\varphi _R(x)  =
\left\{
\matrix
(2R)^{-1}x  &\text{if}  &0 \leq x \leq 2R  \\
\  \  \\
1            &\text{if}  &x > 2R
\endmatrix
\right.  .
$$
\noindent
Then  $\varphi _R$  satisfies the Lipschitz condition  $|\varphi _R(x) - \varphi _R(y)| \leq (2R)^{-1}|x - y|$  for  
$x, y \in [0,\infty )$.  For a non-empty set  $E \subset \Omega $  and  $z \in \Omega $,  we denote 
$d(z,E) = \inf \{d(z,\zeta ) : \zeta \in E\}$  as usual.  
By the triangle inequality for  $d$,  $|d(z,E) - d(w,E)| \leq d(z,w)$  for all   $z, w \in \Omega $.  Hence
$$
|\varphi _R(d(z,E)) - \varphi _R(d(w,E))|  \leq  (2R)^{-1}|d(z,E) - d(w,E)| \leq (2R)^{-1}d(z,w)
\tag 7.5
$$
\noindent
for all   $z, w \in \Omega $.    We now define
$$
g_R(z) = f(z)(1- \varphi _R(d(z,\Omega _t)))   \quad \text{and}  \quad  
h_R(z) = f(z)\varphi _R(d(z,\Omega _t)),
$$
\noindent
where    $\Omega _t = \{\zeta \in \Omega : -r(\zeta ) \geq t\}$.   Since  $\|f\|_\infty \leq 1$   and  $\|\varphi _R \|_\infty = 1$,  we have
$$
|h_R(z) - h_R(w)|  \leq  |f(z) - f(w)|  +  |\varphi _R(d(z,\Omega _t)) - \varphi _R(d(w,\Omega _t))|.
$$
\noindent
If  $h_R(z) - h_R(w) \neq 0$,  then either  $z \in H_t$  or  $w \in H_t$.
Thus if  $d(z,w) \leq 1$  and  $h_R(z) - h_R(w) \neq 0$,  then it follows from (7.4) and (7.5) that  
$|h_R(z) - h_R(w)| \leq 1/R$.  That is,  $\text{diff}(h_R) \leq 1/R$  as promised.   On the other hand,   if  $g_R(z) \neq 0$,  then  
$d(z,\Omega _t) < 2R$.   By Lemma 2.1,  this means that   $-r(z) \geq c(R)t$,  where   $c(R) > 0$  is a constant determined by  $R$.
Hence the support of  $g_R$  is a compact set contained in   $\Omega $.   This completes the proof.   $\square $
\medskip
\noindent
{\bf  Lemma 7.4.}    {\it Let}   $f_1, \dots ,f_k \dots $  {\it be a sequence of continuous functions on}   $\Omega $
{\it satisfying the following four conditions}:

(1)   {\it There is a}  $0 < C < \infty $  {\it such that}   $\|f_k\|_\infty \leq  C$  {\it for every}  $k \in $  {\bf N}.

(2)   {\it For every}  $k \in $  {\bf N},  {\it there exist}   $a_k > b_k > 0$     {\it such that}   $f_k = 0$   {\it on}   
$\Omega _{a_k}\cup H_{b_k}$.

(3)   $\lim _{k\rightarrow \infty }a_k = 0$.

(4)    $\lim _{k\rightarrow \infty }\text{diff}(f_k) = 0$.

\noindent
{\it Then there is an infinite subset}   $I$  {\it of}  {\bf N}  {\it such that}   $f_J \in {\text{VO}}_{\text{bdd}}$   
{\it for every}   $J \subset I$,   {\it where} 
$$
f_J =  \sum _{k\in J}f_k.
$$
\medskip
\noindent
{\it Proof}.     By condition  (3)  and  Lemma 2.1,  we can inductively pick 
a sequence of natural numbers    $k(1) < k(2) < \cdots < k(j) < \cdots $   such that   $a_{k(j+1)} < b_{k(j)}$   and
$$
d(z,w)   \geq  2   \quad \text{if}  \ \  -r(z) \leq a_{k(j+1)}  \ \  \text{and}  \ \  -r(w) \geq b_{k(j)}
\tag 7.6
$$
\noindent
for every   $j \in $  {\bf N}.      Let   $I = \{k(1),k(2),\dots , k(j), \dots \}$.
\medskip
For each   $k \in $  {\bf N},  define   ${\Cal R}_k = \{z \in \Omega :  b_k \leq -r(z) \leq a_k\}$.   Then (2) says that     $f_k = 0$    
on  $\Omega \backslash {\Cal R}_k$.   It follows from (7.6)  that 
$$
\text{if} \  z \in {\Cal R}_{k(j)}  \  \text{and}  \  w \in {\Cal R}_{k(j')}  \ \text{for} \  j \neq j'  \ \text{in}   \  {\bold N},  \ \text{then}  \  d(z,w) \geq 2.  
\tag 7.7
$$
This immediately implies that if   $J \subset I$,  then  $f_J$  is continuous on  $\Omega $.   Moreover,  since  
${\Cal R}_{k(j)} \cap {\Cal R}_{k(j')} = \emptyset $  whenever  $j \neq j'$,  it follows from (1) and (2)  that   $\|f_J\|_\infty \leq C$  for every  
$J \subset I$.  That is, such an  $f_J$  is bounded on  $\Omega $.
\medskip
Let  $j_0 \in $  {\bf N},   and let  $z, w \in \Omega $    satisfy the conditions   $-r(z) \leq a_{k(j_0)}$   and  $d(z,w) \leq 1$.   Then 
it follows from  (7.7) that there is at most one  $j \in $  {\bf N}   such that   $f_{k(j)}(z) - f_{k(j)}(w)  \neq  0$.   Furthermore,  by (7.6),  if such a  
$j$  exist,  then it must satisfy the condition   $j \geq j_0$.  Thus  for  
$z, w \in \Omega $  satisfying the conditions   $-r(z) \leq a_{k(j_0)}$   and  $d(z,w) \leq 1$,  we have
$$
|f_J(z) - f_J(w)|  \leq  \sup \{\text{diff}(f_{k(j)}) : j \geq j_0\}
$$
\noindent
for every   $J \subset I$.  Applying conditions (3) and (4),  this means that for  every  $J \subset I$, 
$f_J$   has vanishing oscillation.   $\square $

\medskip
\noindent
{\bf Definition 7.5.}    (a)  For each $t > 0$,  the symbol  $\Lambda (t)$   denotes the collection of 
{\it continuous}  functions   $g$   on $\Omega $  satisfying the following three conditions:

(1)    $0 \leq g(z) \leq 1$  for every   $z \in \Omega $.

(2)   $g(z) = 1$   when    $z \in \Omega _t = \{\zeta \in \Omega : -r(\zeta ) \geq  t\}$.

(3)   There is a  $t' = t'(g) \in (0,t)$  such that   $g(z) = 0$   whenever   $-r(z) \leq  t'$.

\noindent
(b)   Let    $t > 0$  and   $\delta > 0$.   Then   $\Lambda (t;\delta )$   denotes the collection of functions  
$g \in \Lambda (t)$  satisfying the additional condition   diff$(g) \leq \delta $.

\medskip
\noindent
{\bf Lemma 7.6.}   {\it For  all}    $t > 0$  {\it and}  $\delta > 0$,   {\it we have}    $\Lambda (t;\delta )  \neq \emptyset $.
\medskip
\noindent
{\it Proof}.      This is similar to the proof of Proposition 7.3.   
Let  $\psi : [0,\infty ) \rightarrow [0,1]$  be a Lipschitz function  with Lipschitz constant $\delta $.   
Furthermore,  suppose that  $\psi (0) = 1$   and that  $\psi = 0$  on  $[R,\infty )$  for a sufficiently large  $R$.
Let  $t > 0$  be such that   $\Omega _t \neq \emptyset $  (otherwise, (2) is trivial).  
By Lemma 2.1,   the function  $f(z) =  \psi (d(z,\Omega _t))$  is in   $\Lambda (t;\delta )$.  $\square $ 
\medskip
\noindent
{\bf Lemma 7.7.}  {\it Given any pair of}  $f \in L^\infty (\Omega )$  {\it and}  $h \in L_a^2(\Omega )$,  {\it we have}
$$
\lim _{t\downarrow 0}\sup \{\|T_{fg}h - T_fh\| :  g \in \Lambda (t)\}   =  0.   
\tag 7.8
$$
\medskip
\noindent
{\it Proof}.     Denote  $H_t = \{z \in \Omega : -r(z) < t\}$  as before.  By  Definition 7.5(a),     we have
$$
\|T_{fg}h - T_fh\|^2  \leq \|fgh - fh\|^2   \leq  \|f\|_\infty ^2\int _{H_t}|h(z)|^2dv(z)
$$
\noindent
for all  $g \in \Lambda (t)$,  $f\in L^\infty (\Omega )$  and  $h \in L_a^2(\Omega )$.  This obviously implies (7.8).   $\square $
\medskip
For a bounded operator  $A$  on a Hilbert space  ${\Cal H}$,      denote 
$$
\|A\|_{\Cal Q}  = \inf \{\|A + K\| :  K \ \text{is any compact operator on} \ {\Cal H}\},
$$
\noindent
which is the essential norm of    $A$.
\medskip
\noindent
{\bf Lemma 7.8.}  [16,Lemma 2.1]   {\it Let}   $\{B_i\}$  {\it be a sequence of compact operators on a Hilbert space}  ${\Cal H}$
{\it satisfying the following conditions}:

\noindent
(a)  {\it Both sequences}  $\{B_i\}$  {\it and}  $\{B_i^\ast \}$  {\it converge to} $0$ {\it in the strong operator topology}.

\noindent
(b)  {\it The limit}  $\lim _{i\rightarrow \infty }\|B_i\|$  {\it exists}.

\noindent
{\it Then there exist natural numbers}  $i(1) < i(2) <  \cdots < i(m) < \cdots $     {\it such that the sum}
$$
\sum _{m=1}^\infty B_{i(m)}  =  \lim _{N\rightarrow \infty }\sum _{m=1}^N B_{i(m)}  
$$
\noindent
{\it exists in the strong operator topology and we have}  
$$
\bigg\|\sum _{m=1}^\infty B_{i(m)}\bigg\|_{\Cal Q} =  \lim _{i\rightarrow \infty }\|B_i\|.
$$
\medskip
\noindent
{\bf Definition 7.9.}    For $t > 0$  and   $\delta > 0$,   the symbol  $\Phi (t;\delta )$  denotes the collection of 
{\it continuous}  functions   $f$   on  $\Omega $  satisfying the following three conditions:

(1)    $0 \leq f(z) \leq 1$  for every   $z \in \Omega $.

(2)   $f(z) = 0$   whenever   $-r(z) \geq  t$.

(3)  diff$(f) \leq \delta $.

\medskip
In analogy with [28,Proposition 3.7],   every operator 
in  EssCom$(\{T_g : g \in {\text{VO}}_{\text{bdd}}\})$  satisfies the following ``$\epsilon $-$\delta $"  condition:

\medskip
\noindent
{\bf Proposition 7.10.}   {\it Let}  $X$  {\it be an operator in the essential commutant of}  $\{T_g : g \in {\text{VO}}_{\text{bdd}}\}$.   
{\it Then for every}   $\epsilon > 0$,  {\it there is a}  $\delta = \delta (X,\epsilon ) > 0$  {\it such that}   
$$
\lim _{t\downarrow 0}\sup \{\|[X,T_f]\| : f \in \Phi (t;\delta )\}   \leq  \epsilon .
$$
\medskip
Using 7.4-7.9 above,  the proof of Proposition 7.10 is a repeat of the proof of Proposition 3.7 in [28], 
modified in the obvious way.  For that reason we will omit the proof of Proposition 7.10 here.   
\medskip
\noindent
{\bf  Lemma 7.11.}    {\it Let}   $h_1, \dots ,h_k \dots $  {\it be a sequence of  continuous functions on}  $\Omega $,  {\it and denote}  
$U_k = \{z \in \Omega : h_k(z) \neq 0\}$,  $k \in $  {\bf N}.
{\it Suppose that this sequence has the property that there is an}   $a > 1$  {\it such that}   
$\inf \{d(z,w) :  z \in U_j, w \in U_k\} \geq a$   {\it for every pair of}   $j \neq k$  {\it in}  {\bf N}.
{\it Then the function}    $h =  \sum _{k=1}^\infty h_k$   {\it has the property that}   diff$(h) \leq \sup _{k\in {\bold N}}\text{diff}(h_k)$.  
\medskip
\noindent
{\it Proof}.      Observe that, under  the assumption,  for any pair of  $z, w \in \Omega $   satisfying the condition  
$d(z,w) \leq 1$,   the cardinality of the set  $\{k \in {\bold N}  : h_k(z) - h_k(w) \neq 0\}$  is at most  $1$.      $\square $

\bigskip
\centerline{\bf 8.  Approximate partition of unity}
\medskip
In this section the boundary   $\partial \Omega $  of the domain plays a prominent role.   
It will be beneficial to make a simplification of notation:  for     $\zeta \in \partial \Omega $  and    $t > 0$, let us write
$$
Q(\zeta ,t)  =  \{\xi \in \partial \Omega :  |\zeta - \xi |^2 + |\langle \zeta - \xi ,(\bar \partial r)(\zeta )\rangle | < t\}.
$$
\noindent
In other words, in terms of the notation in Section 2,  we have  $Q(\zeta ,t) = Q_0(\zeta ,t)$.  Similarly,  we will write  
$d\sigma $  for  $d\sigma _0$. That is,   $d\sigma $  is the surface measure on   $\partial \Omega $.
\medskip
\noindent
{\bf Lemma 8.1.}  {\it  There is a constant}  $1 \leq C_{8.1} < \infty $  {\it  such that for any triple of}   
$\zeta ,\xi \in \partial \Omega $  {\it  and}  $t > 0$,   {\it  if}  $Q(\zeta ,t)\cap Q(\xi ,t) \neq \emptyset $,  
{\it  then}   $Q(\xi ,t) \subset  Q(\zeta ,C_{8.1}t)$.
\medskip
\noindent
{\it Proof}.   By the assumption on the defining function  $r$,  there is an  $L$  such that   
$|(\bar \partial r)(x) - (\bar \partial r)(y)| \leq L|x - y|$   for all 
$x,y \in \partial \Omega $.  Suppose that there is a   $w \in Q(\zeta ,t)\cap Q(\xi ,t)$.   Then 
$$
|\zeta - w|^2 + |\langle \zeta - w ,(\bar \partial r)(\zeta )\rangle | <  t   \quad  \text{and} \quad   
|\xi     - w|^2 + |\langle \xi     - w ,(\bar \partial r)(\xi      )\rangle | <  t.
$$
\noindent
From this it  is elementary to obtain $|\zeta - \xi |^2 \leq 2(|\zeta - w|^2 + |\xi - w|^2) < 2t$.   Further,
$$
\align
|\langle \zeta - \xi ,(\bar \partial r)(\zeta )\rangle |  &\leq  |\langle \zeta - w,(\bar \partial r)(\zeta )\rangle | + |\langle w - \xi ,(\bar \partial r)(\zeta )\rangle |   \\
&\leq   |\langle \zeta - w,(\bar \partial r)(\zeta )\rangle |  + |\langle w - \xi ,(\bar \partial r)(\xi )\rangle |  + |w - \xi |\cdot L|\xi - \zeta |  \\
&< 2t + \sqrt{t}\cdot L\sqrt{2t}  =  (2 + \sqrt{2}L)t.
\endalign
$$
\noindent
Therefore if we set   $C_1 = 4 + \sqrt{2}L$,  then  the condition  $Q(\zeta ,t)\cap Q(\xi ,t) \neq \emptyset $ implies  $\xi \in Q(\zeta ,C_1t)$.   
Suppose that  $z \in Q(\xi ,t)$, i.e.,  $|\xi - z|^2 + |\langle \xi - z,(\bar \partial r)(\xi )\rangle | < t$.   Then
$$
|\langle \xi - z,(\bar \partial r)(z)\rangle |  \leq  |\langle \xi - z,(\bar \partial r)(\xi )\rangle | + |\xi -z| \cdot L|\xi -z| < (1+L)t.
$$
\noindent
That is, if   $z \in Q(\xi ,t)$,  then  $\xi \in Q(z,(2+L)t) \subset Q(z,C_1t)$.  Thus the condition   $Q(\zeta ,t)\cap Q(\xi ,t) \neq \emptyset $ implies
 $Q(\zeta ,C_1t)\cap Q(z ,C_1t) \neq \emptyset $   for every   $z \in Q(\xi ,t)$.   By the first conclusion,  we have  $z \in Q(\zeta ,C_1^2t)$.  
 Thus the lemma holds for  $C_{8.1} = C_1^2 =  (4 + \sqrt{2}L)^2$.    $\square $
 \medskip
\noindent
{\bf Corollary 8.2.} {\it  Consider any}  $\zeta \in \partial \Omega $  {\it  and}  $t > 0$.  {\it  If}   $x, y \in \partial \Omega $    
{\it  are such that}   $x \in Q(\zeta ,t)$  {\it  and}  $y \notin Q(\zeta ,C_{8.1}t)$,   {\it  then}   $y \notin Q(x,t)$. 
\medskip
\noindent
{\it Proof}.   If  $x \in Q(\zeta ,t)$,  then  $Q(x,t)\cap Q(\zeta ,t) \neq \emptyset $.  By Lemma 8.1,  we have
$Q(x,t) \subset Q(\zeta ,C_{8.1}t)$.  Therefore if  $y \notin Q(\zeta ,C_{8.1}t)$,   then   $y \notin Q(x,t)$.    $\square $
\medskip
\noindent
{\bf Lemma 8.3.}  {\it  There is a constant}  $0 < C_{8.3} < \infty $  {\it  such that the following bound holds}:  Let   $t > 0$,  {\it  and let}  
$E$  {\it  be a subset of}  $\partial \Omega $   {\it  that has the property}    $Q(x,t)\cap Q(y,t) = \emptyset $  {\it  for all}  
$x \neq y$   {\it  in}    $E$.   {\it  Then for any}  $R \geq 1$  {\it  and any}   $\zeta \in \partial \Omega $,
$$
\text{card}\{x \in E : Q(x,Rt)\cap Q(\zeta ,Rt) \neq \emptyset \}  \leq  C_{8.3}R^n.
$$
\medskip
\noindent
{\it Proof}.    It  $t > T_0$   (see Proposition 2.8),  
then the property of  $E$  implies   card$(E) \leq 1$.  Suppose that   $0 < t \leq T_0$.
If  $Q(x,Rt)\cap Q(\zeta ,Rt) \neq \emptyset $,  then  $Q(x,Rt) \subset Q(\zeta ,C_{8.1}Rt)$   by Lemma 8.1.    
Let  $E_0 = \{x \in E :  Q(x,Rt)\cap Q(\zeta ,Rt) \neq \emptyset \}$.   Since     
$Q(x,t)\cap Q(y,t) = \emptyset $  for all  $x \neq y$  in  $E$,  we have
$$
\text{card}(E_0)\inf _{x\in E_0}\sigma (Q(x,t))  \leq \sum _{x\in E_0}\sigma (Q(x,t))  
=  \sigma \bigg( \bigcup _{x\in E_0}Q(x,t)\bigg)   \leq  \sigma (Q(\zeta ,C_{8.1}Rt)).
$$
\noindent
Applying Proposition 2.8 to the case  $\rho = 0$,   we have
$$
\text{card}(E_0)c_{2.8}t^n   \leq  C_{2.8}(C_{8.1}Rt)^n.
$$
\noindent
Cancelling out  $t^n$  and simplifying,  we see that the lemma holds for the constant  
$C_{8.3} = (C_{2.8}/c_{2.8})C_{8.1}^n$.    $\square $
\medskip
The first order Taylor expansion for  $r$  reads
$$
r(z+u)  =  r(z) + 2\text{Re}\langle u,(\bar \partial r)(z)\rangle + \int _0^12\text{Re}\langle u,(\bar \partial r)(z+xu) - (\bar \partial r)(z)\rangle dx.
$$
\noindent
Thus $r(z+t(\bar \partial r)(z))  =  r(z) + 2t|(\bar \partial r)(z)|^2 +  O(t^2)$  when  $|t|$  is small.   
Recall that  $\bar \partial r$  does not vanish on  $\partial \Omega $.   Hence
when  $z$  is close to  $\partial \Omega $  and  $t$  is positive and small,  we have
$$
r(z+t(\bar \partial r)(z))  \geq  r(z) + t|(\bar \partial r)(z)|^2.
$$
\noindent
Thus for each   $z \in \Omega $  near  $\partial \Omega $,  there is a  $t_z > 0$,  $t_z \approx |r(z)|$,  such that   $r(z+t_z(\bar \partial r)(z)) = 0$.
Let us restate this fact more precisely:    There  exist a  $J \in $  {\bf N}  and  a  $0 < C_p < \infty $  such that for every  
$z \in H_{2^{-J}} = \{\zeta \in \Omega : -r(\zeta ) < 2^{-J}\}$,  there is a  $p(z) \in \partial \Omega $  such that
$$
|z - p(z)|   \leq  C_p|r(z)|.
\tag 8.1
$$
\noindent
In other words,  there is a map  $p : H_{2^{-J}} \rightarrow \partial \Omega $  such that the above bound holds for every  
$z \in H_{2^{-J}}$.     Note that our choice above does not promise any kind of continuity for the map   $p$,  but that 
does not matter for our purpose.  
\medskip
This   $p$  and  the defining function  $r$  together allow us to decompose  $H_{2^{-J}}$  in a manner that is analogous to the radial-spherical 
decomposition for the unit ball in [28].    More specifically,   $p$   plays the role of  ``spherical coordinates",  while   $-r$   is the analogue of
``radial coordinate".    Because we only need  a large-scale, or ``coarse", decomposition,  (8.1) is all that we need to 
know about  $p$.
\medskip
\noindent
{\bf Lemma 8.4.}   {\it  There is a constant}  $0 < C_{8.4} < \infty $  {\it  such that}
$$
\align
|z' - w'|^2 &+ |\langle z' - w',(\bar \partial r)(z')\rangle |  \\
&\leq  3\{|z - w|^2 + |\langle z - w,(\bar \partial r)(z)\rangle |\}  +  C_{8.4}\{|z - z'| + |w - w'|\} 
\endalign
$$ 
\noindent
{\it  for all}  $z, w, z', w'  \in \overline{\Omega } = \Omega \cup \partial \Omega $. 
\medskip
\noindent
{\it Proof}.  It is elementary that  $|z'  - w'|^2 \leq 3|z - w|^2 + 3|z - z'|^2 + 3|w - w'|^2$.  
Since  $\Omega $  is bounded,  there is a  $C_1$  such that  $|\zeta - \xi | \leq C_1$  for all  
$\zeta , \xi \in \overline{\Omega }$.   Hence
$$
|z' - w'|^2 \leq 3|z - w|^2 + 3C_1\{|z - z'| + |w - w'|\}.
\tag 8.2
$$
\noindent
Similarly,  since  $\bar \partial r$  is bounded and satisfies a Lipschitz condition on  $\overline{\Omega }$,  we have
$$
\align
|\langle z' - w',&(\bar \partial r)(z')\rangle |  \leq   (|z - z'| + |w - w'|)|(\bar \partial r)(z')|  +  |\langle z - w,(\bar \partial r)(z')\rangle |  \\
&\leq   C_2(|z - z'| + |w - w'|)   +   |\langle z - w,(\bar \partial r)(z)\rangle |  +  |z - w||(\bar \partial r)(z') - (\bar \partial r)(z)|  \\
&\leq   C_2(|z - z'| + |w - w'|)   +   |\langle z - w,(\bar \partial r)(z)\rangle |  +  C_3|z - z'|.
\tag 8.3
\endalign
$$
\noindent
Obviously, the lemma follows from (8.2) and (8.3).   $\square $
\medskip
We begin the decomposition with natural numbers   $m >  J$   and  $j \geq 1$.   Define
$$
d_{m,j}  =  m2^{-jm},   \quad  a_{m,j}  =  C_{8.1}m2^{-jm}  \quad  \text{and}  \quad   b_{m,j}  =  C_{8.1}^2m2^{-jm}, 
\tag 8.4
$$
\noindent
where  $C_{8.1}$  is the constant in Lemma 8.1.  That is,   $a_{m,j} = C_{8.1}d_{m,j}$  and  $b_{m,j} = C_{8.1}^2d_{m,j}$.
Let   $E_{m,j}$  be a subset of  $\partial \Omega $  that is {\it maximal}  with respect to the property
$$
Q(u,d_{m,j})\cap Q(v,d_{m,j})  =  \emptyset  \quad \text{for all}  \ \  u \neq v \in E_{m,j}.
\tag 8.5
$$
\noindent
By the maximality of  $E_{m,j}$  and  Lemma 8.1,  we have
$$
\bigcup _{u\in E_{m,j}}Q(u,a_{m,j})   =  \partial \Omega .
\tag 8.6
$$
\noindent
Fix a natural number  $N_0$   such that   $N_0 \geq C_{8.3}(C_{8.1}^2)^n$,  where  $C_{8.3}$  is the constant in Lemma 8.3.   
Since  $b_{m,j} = C_{8.1}^2d_{m,j}$,  it follows from (8.5)  and  Lemma 8.3 that
$$
\text{card}\{v \in E_{m,j} : Q(v,b_{m,j})\cap Q(u,b_{m,j})  \neq \emptyset \}   \leq N_0
\tag 8.7
$$
\noindent
for every   $u \in E_{m,j}$.  Now, given any   $m >  J$,   $j \geq 1$   and  $u \in E_{m,j}$,  we  define the sets
$$
\align
A_{m,j,u}   &=  \{z \in \Omega :  p(z) \in Q(u,a_{m,j})  \ \  \text{and} \ \  2^{-(j+1)m} > - r(z)  \geq 2^{-(j+2)m}\}    \quad  \text{and}  \\
B_{m,j,u}   &=  \{z \in \Omega :  p(z) \in Q(u,b_{m,j})  \ \  \text{and} \ \  2^{-jm} > - r(z) > 2^{-(j+3)m}\}.
\endalign
$$
\noindent
It follows from (8.6) that
$$
\bigcup _{j=1}^\infty \bigcup _{u\in E_{m,j}}A_{m,j,u}  =   H_{2^{-2m}} = \{z \in \Omega : -r(z) < 2^{-2m}\}.
\tag 8.8
$$
\noindent
Note that even though we have (8.8),  we do not know that every  $A_{m,j,u}$  
is non-empty from its definition.   Nevertheless,  we have
\medskip
\noindent
{\bf Lemma 8.5.}   {\it There is a constant}  $J < M_{8.5}  < \infty $    {\it such that if}  $m \geq  M_{8.5}$,   {\it then}  
$A_{m,j,u} \neq \emptyset $   {\it for all}   $j \geq 1$   {\it and}  $u \in E_{m,j}$.
\medskip
\noindent
{\it Proof}   By the Taylor expansion for  $r$,  there are constants   $J < M_1 < \infty $   and  $0 < C_1 < \infty $   such that  
if  $m \geq M_1$,  then for every pair of  $j \geq 1$   and   $u \in E_{m,j}$   there is a  $u'$  such that   
$-r(u') = 2^{-(j+(3/2))m}$   and  $|u - u'| \leq C_1(-r(u'))$.  By Lemma 8.4,  we have
$$
|u - u'|^2 + |\langle u - u',(\bar \partial r)(u)\rangle |  \leq  C_{8.4}|u - u'|  \leq  C_{8.4}C_1(-r(u'))  =  C_{8.4}C_12^{-(j+(3/2))m}.
$$
\noindent
Applying Lemma 8.4 again and recalling (8.1),  we have
$$
\align
|u - p(u')|^2 +& |\langle u - p(u'),(\bar \partial r)(u)\rangle |  \leq   3\{|u - u'|^2 + |\langle u - u',(\bar \partial r)(u)\rangle |\} + C_{8.4}|u' - p(u')|  \\
&\leq 3C_{8.4}C_12^{-(j+(3/2))m} + C_{8.4}C_p2^{-(j+(3/2))m}.
\endalign
$$
\noindent
Let  $M_{8.5} \geq M_1$  be  such that   $M_{8.5} \geq 3C_{8.4}C_1 + C_{8.4}C_p$.   If  $m \geq M_{8.5}$,  
then  $u' \in A_{m,j,u}$.   $\square $
\medskip
\noindent
{\bf Lemma 8.6.}   {\it There is a constant}  $M_{8.5} + 100 \leq  M_{8.6} < \infty $    {\it such that for}  $m \geq  M_{8.6}$,       $j \geq 1$,   
{\it and}   $u \in E_{m,j}$,   {\it if}  $z \in A_{m,j,u}$   {\it and}  $w \in \Omega \backslash B_{m,j,u}$,    {\it then}   $d(z,w) \geq(1/13)m$.  
\medskip
\noindent
{\it Proof}.    Set   $M_1 = \max \{M_{8.5}+100,10C_{8.4}\}$,  where  $C_{8.4}$  and  $M_{8.5}$  are the constants in Lemmas 
8.4  and 8.5   respectively.  Consider any  $m \geq  M_1$,   $j \geq 1$   and  $u \in E_{m,j}$.
For a pair of  $z \in A_{m,j,u}$   and  $w \in \Omega \backslash B_{m,j,u}$,  there are three possibilities,   
depending on the value of  $r(w)$.
\medskip
(1)   Suppose that  $-r(w) \geq 2^{-jm}$.   Then  $r(z)/r(w) \leq 2^{-(j+1)m}/2^{-jm} = 2^{-m}$.   Combining this with 
Lemma 2.1,  we have   $c_{2.1}2^{-4d(w,z)} \leq r(z)/r(w) \leq 2^{-m}$.   Hence  
$$
d(z,w) \geq (1/4)m + (1/4)\{\log c_{2.1}/\log 2\}.
$$
\noindent
Let  $M_2 \geq M_1$  be such that   $(1/2)M_2 \geq  |\log c_{2.1}/\log 2|$.    Thus if  $m \geq M_2$, then  for all   $j \geq 1$,    $u \in E_{m,j}$,
$z \in A_{m,j,u}$   and  $w \in \Omega \backslash B_{m,j,u}$,  we have  
$$
d(z,w) \geq (1/8)m
\tag 8.9
$$
\noindent
under the condition  $-r(w) \geq 2^{-jm}$. 
\medskip
(2)   Suppose that  $-r(w) \leq 2^{-(j+3)m}$.   Then  $r(w)/r(z) \leq 2^{-(j+3)m}/2^{-(j+2)m} = 2^{-m}$.   From Lemma 2.1  we now deduce
$c_{2.1}2^{-4d(z,w)} \leq r(w) /r(z) \leq 2^{-m}$.   Thus  (8.9)  again holds under the condition  $-r(w) \leq 2^{-(j+3)m}$   when  $m \geq M_2$.
\medskip
(3)   Suppose that  $2^{-(j+3)m} < -r(w) < 2^{-jm}$.   Then by the definition of  $B_{m,j,u}$ we have  
$p(w) \notin Q(u,b_{m,j})$.  In contrast,  since    $z \in A_{m,j,u}$,  we have  $p(z) \in Q(u,a_{m,j})$.   Since  $b_{m,j} = C_{8.1}a_{m,j}$,  
by Corollary 8.2 we have   $p(w) \notin Q(p(z),a_{m,j})$.    Recall that  $a_{m,j} = C_{8.1}m2^{-jm}$  and  that  $C_{8.1} \geq 1$.
Thus it follows from Lemma 8.4 and (8.1)  that
$$
\align
m2^{-mj}  \leq a_{m,j}  &\leq  |p(z) - p(w)|^2 + |\langle p(z) - p(w),(\bar \partial r)(p(z))\rangle |  \\
&\leq  3\{|z - w|^2 + |\langle z - w,(\bar \partial r)(z)\rangle |\}  +  C_{8.4}C_p\{|r(z)| + |r(w)|\}   \\
&\leq  3\{|z - w|^2 + |\langle z - w,(\bar \partial r)(z)\rangle |\}  +  C_{8.4}C_p\{2^{-(j+1)m} + 2^{-jm}\}.
\endalign
$$
\noindent
Now we pick an  $M_3 \geq M_2$  such that   $M_3 \geq 4C_{8.4}C_p$,  i.e.,  $(1/2)M_3 \geq 2C_{8.4}C_p$.
When  $m \geq M_3$,  elementary manipulations turn the above into the inequality
$$
(1/6)m2^{-mj}  \leq  |z - w|^2 + |\langle z - w,(\bar \partial r)(z)\rangle |.
$$
\noindent
Combining this with Lemma 2.2,  we obtain  
$$
(1/6)m2^{-mj} \leq  C_{2.2}\{d(z,w) +d^2(z,w)\}2^{12d(z,w)}(-r(z)).
$$
\noindent
Since   $-r(z) \leq 2^{-(j+1)m}$,  this implies 
$$
(1/6)m2^{m} \leq  C_{2.2}\{d(z,w) +d^2(z,w)\}2^{12d(z,w)}.
$$
\noindent
From this inequality it is elementary to deduce that there is an  $M_{8.6} \geq M_3$  such that 
if  $m \geq M_{8.6}$,  then  $d(z,w) \geq (1/13)m$.   Combining this with (8.9),  the proof is complete.  $\square $
\medskip
By Lemma 8.5,  for every triple of  $m \geq M_{8.5}$,  $j \geq 1$   and  $u \in E_{m,j}$,  we can pick a
$$
z_{m,j,u} \in A_{m,j,u}.
\tag 8.10
$$
\noindent
This pick will be fixed for the rest of the paper.
\medskip
\noindent
{\bf Lemma 8.7.}   {\it There is a constant}  $M_{8.6} < M_{8.7} < \infty $
{\it such that if}  $m \geq  M_{8.7}$,       {\it then there is an}    $0 < R_m < \infty $   {\it  which has the property that}
$$
B_{m,j,u}  \subset D(z_{m,j,u},R_m)
\tag 8.11
$$
{\it for all}   $j \geq 1$    {\it and}   $u \in E_{m,j}$.  
\medskip
\noindent
{\it Proof}.     Suppose that  $m \geq M_{8.5}$.  Given any  $j \geq 1$  and   $u \in E_{m,j}$,  
we have   $2^{-(j+2)m} \leq - r(z_{m,j,u}) < 2^{-(j+1)m}$  by (8.10).   Now let   $w \in B_{m,j,u}$.   Then
$2^{-(j+3)m} < - r(w) < 2^{-jm}$, which means   $-2^{-2m}r(z_{m,j,u}) \leq - r(w) \leq  -2^{2m}r(z_{m,j,u})$.   
In other words,  we have
$$
2^{k-1}(-r(z_{m,j,u}))  \leq - r(w)  \leq  2^k(-r(z_{m,j,u})) \quad  \text{for some} \ k \in {\bold Z}  \ \text{with} \  |k| \leq 2m.
\tag 8.12
$$
\noindent    
 We have    $p(w) \in Q(u,b_{m,j})$. 
Since    $p(z_{m,j,u}) \in  Q(u,a_{m,j}) \subset  Q(u,b_{m,j})$,   Lemma 8.1 gives us
$Q(u,b_{m,j})$  $\subset $  $Q(p(z_{m,j,u}),C_{8.1}b_{m,j})$.     
Hence    $p(w) \in Q(p(z_{m,j,u}),C_{8.1}b_{m,j})$.   That is,
$$
|p(z_{m,j,u}) - p(w)|^2 + |\langle p(z_{m,j,u}) - p(w),(\bar \partial r)(p(z_{m,j,u}))\rangle| < C_{8.1}b_{m,j}.
$$
\noindent
Applying Lemma 8.4 and (8.1),  we obtain
$$
\align
|z_{m,j,u} - w|^2 + &|\langle z_{m,j,u} -  w,(\bar \partial r)(z_{m,j,u})\rangle| < 3C_{8.1}b_{m,j} + C_{8.4}C_p(|r(z_{m,j,u})| + |r(w)|)  \\
&\leq 3C_{8.1}^3m2^{-jm} + 2C_{8.4}C_p2^{-jm}   <  (3C_{8.1}^3m + 2C_{8.4}C_p)2^{2m}(-r(z_{m,j,u})).
\endalign
$$
\noindent
Let  $M_{8.7} > M_{8.6}$  be such that   $(3C_{8.1}^3M_{8.7} + 2C_{8.4}C_p)2^{-M_{8.7}} \leq 1$.  When  $m \geq M_{8.7}$,  
the above inequality gives us
$$
\align
|z_{m,j,u} - w|^2 + |\langle z_{m,j,u} -  w,(\bar \partial r)(z_{m,j,u})\rangle| &<  2^{3m}(-r(z_{m,j,u}))   \\
&\leq  2^{k + 3m + |k|}(-r(z_{m,j,u})).   
\tag 8.13
\endalign
$$
\noindent
Combining (8.12) and (8.13) with Lemma 3.7,  we obtain   $d(z_{m,j,u},w) < C_{3.7}(1+|k|+3m+|k|) \leq C_{3.7}(1+7m)$.  
Thus when  $m \geq M_{8.7}$,  (8.11) holds for  $R_m = C_{3.7}(1+7m)$.    $\square $
\medskip
In the above we picked constants such that  $M_{8.7} > M_{8.6} \geq  M_{8.5} + 100$  and  
$M_{8.5} > J$.   Thus if    $m \geq M_{8.7}$,  then   $m/13 > 7$.   Now, for every  $m \geq M_{8.7}$,   we define the function
$$
\tilde f_m(x) = 
\left\{
\matrix
1 - \{(m/13) - 4\}^{-1}x  &\text{for}  &0 \leq x \leq (m/13) - 4    \\
\ \  \\
0  &\text{for}  &(m/13) - 4  < x < \infty 
\endmatrix
\right.  .
\tag 8.14
$$
\noindent
Obviously,  $\tilde f_m$  satisfies the Lipschitz condition   $|\tilde f_m(x) - \tilde f_m(y)| \leq \{(m/13) - 4\}^{-1}|x-y|$,       
$x, y \in [0,\infty )$.   For every triple of  $m \geq M_{8.7}$,    $j \in $ {\bf N}  and  $u \in E_{m,j}$,  we  define
$$
f_{m,j,u}(z)   =  \tilde f_m (d(z,A_{m,j,u}))  \quad \text{for}  \ \  z \in \Omega .
$$
\medskip
\noindent
{\bf Lemma 8.8.}  {\it For every triple of}  $m \geq M_{8.7}$,    $j \in $ {\bf N}  {\it and}  $u \in E_{m,j}$,  {\it the function}    
$f_{m,j,u}$  {\it defined above has the following five properties}:

\noindent
(a)   {\it The inequality}   $0 \leq f_{m,j,u} \leq 1$  {\it  holds on}  {\bf B}.

\noindent
(b)    $f_{m,j,u} = 1$  {\it on the set}  $A_{m,j,u}$.

\noindent
(c)   $|f_{m,j,u}(z) - f_{m,j,u}(w)| \leq \{(m/13) - 4\}^{-1}d(z,w)$   {\it for all}   $z, w \in \Omega $.

\noindent
(d)    {\it  If}   $f_{m,j,u}(z) \neq 0$  {\it and}   $w \in \Omega \backslash B_{m,j,u}$,  {\it then}   $d(z,w) \geq 4$.  

\noindent
(e)   {\it We have}  diff$(f_{m,j,u}) \leq \{(m/13) - 4\}^{-1}$.
\medskip
\noindent
{\it Proof}.   (a) and (b) follow directly from the definitions of  $\tilde f_m$   and  $f_{m,j,u}$.   
Then note that
$$
\align
|f_{m,j,u}(z) - f_{m,j,u}(w)|  &= |\tilde f_m (d(z,A_{m,j,u})) - \tilde f_m (d(w,A_{m,j,u}))|     \\
&\leq  {1\over (m/13) - 4}|d(z,A_{m,j,u}) - d(w,A_{m,j,u})|  \leq  {d(z,w)\over (m/13) - 4},
\endalign
$$ 
\noindent
which proves  (c).   For (d),  observe that if  $f_{m,j,u}(z) \neq 0$,  then 
$d(z,A_{m,j,u}) < (m/13) - 4$.   This means that there is a  $z' \in A_{m,j,u}$  such that  
$d(z,z') \leq (m/13) - 4$.   If   $w \in \Omega \backslash B_{m,j,u}$,  then Lemma 8.6 tells us that  
$d(z',w) \geq m/13$.   By the triangle inequality,
$$
d(z,w)  \geq d(z',w) - d(z,z')  \geq  (m/13) - \{ (m/13) - 4\}  = 4.
$$
\noindent
Hence  (d) holds.  Finally, note that  (e)  is an immediate consequence of  (c).     $\square $
\medskip
By (8.7) and a standard  maximality argument,  each  $E_{m,j}$  admits a partition  
$$
E_{m,j}  = E_{m,j}^{(1)}\cup \dots \cup E_{m,j}^{(N_0)}
\tag 8.15
$$
\noindent
such that for every   $\nu \in \{1,\dots ,N_0\}$,  we have   $Q(u,b_{m,j})\cap Q(v,b_{m,j})  =  \emptyset$   for all
$u \neq v$  in $E_{m,j}^{(\nu )}$.    Therefore  for each   $\nu \in \{1,\dots ,N_0\}$,  the conditions  $u, v \in E_{m,j}^{(\nu )}$
and  $u \neq v$  imply   $B_{m,j,u}\cap B_{m,j,v}$  $=$  $\emptyset $.
\medskip
\noindent
{\bf Definition 8.9.}        Let   $m \geq M_{8.7}$  be given.  (a)  For each pair of $\kappa \in \{1,2,3\}$  and  $\nu \in \{1,\dots ,N_0\}$,  
where  $N_0$  is the integer that appears in (8.7) and (8.15),  
let    $I^{(\nu ,\kappa )}_m$  denote the collection of all triples     $m,3j+\kappa ,u$   satisfying the conditions   
$j \in {\bold Z}_+$    and   $u \in E_{m,3j+\kappa }^{(\nu )}$.

\noindent
(b)    For   $\kappa \in \{1,2,3\}$,   $\nu \in \{1,\dots ,N_0\}$  and  $q \in $  {\bf N},
let    $I^{(\nu ,\kappa )}_{m,q}$  denote the collection of all triples     $m,3j+\kappa ,u$   satisfying the conditions   
$0 \leq j \leq q$    and   $u \in E_{m,3j+\kappa }^{(\nu )}$.

\noindent
(c)   Denote  $I_m = \cup _{\kappa = 1}^3\cup _{\nu = 1}^{N_0}I^{(\nu ,\kappa )}_m$.

\medskip
The elements in  $I_m$, equivalently the subscripts in  $A_{m,j,u}$,  $B_{m,j,u}$   and  $f_{m,j,u}$,  are obviously 
quite cumbersome to write as triples.  For this we have the following remedy:
\medskip
\noindent
{\bf Notation 8.10.}    (1)  We will use the symbol  $\omega $ to represent the triple  $m,j,u$.  

\noindent
(2)   For any subset   $I$   of   $I_m$,  denote   $f_I = \sum _{\omega \in I}f_\omega $   and   
$F_I = \sum _{\omega \in I}f_\omega ^2$.
\medskip
\noindent
{\bf Lemma 8.11.}    {\it Let}  $m \geq M_{8.7}$,  $\kappa \in \{1,2,3\}$  {\it and}  $\nu \in \{1,\dots ,N_0\}$.   {\it Then for any}   
$\omega \neq \omega '$  {\it in}   $I_m^{(\nu ,\kappa )}$,   {\it we have}  $B_\omega \cap B_{\omega '} = \emptyset $.
\medskip
\noindent
{\it Proof}.  If  $\omega = (m,3j+\kappa ,u)$  and   $\omega ' = (m,3j+\kappa ,v)$   for a pair of  $u \neq v$  in  $E_{m,3j+\kappa }^{(\nu )}$,  
then by the property of the partition (8.15)  we already know that  $B_\omega \cap B_{\omega '} = \emptyset $.   The other possibility is that  
$\omega = (m,3j+\kappa ,u)$  and   $\omega ' = (m,3j'+\kappa ,v)$    with  $u \in  E_{m,3j+\kappa }^{(\nu )}$   and  
$v \in  E_{m,3j'+\kappa }^{(\nu )}$,  where  $j \neq j'$.   If  $j \neq j'$,  then  $|(3j+\kappa ) - (3j'+\kappa )| \geq 3$,  which ensures  
$B_\omega \cap B_{\omega '} = \emptyset $ by the values of  $-r$  on  $B_\omega $   and  $B_{\omega '}$.   $\square $

\medskip
\noindent
{\bf Lemma 8.12.}    {\it Let}  $m \geq M_{8.7}$,  $\kappa \in \{1,2,3\}$  {\it and}  $\nu \in \{1,\dots ,N_0\}$.    
{\it Then for every subset}   $I$   {\it of}     $I^{(\nu ,\kappa )}_m$,   {\it we have}   $f_I \in \Phi (2^{-m};((m/13) - 4)^{-1})$.
\medskip
\noindent
{\it Proof}.     Let  $I \subset I^{(\nu ,\kappa )}_m$.   
For each   $\omega \in I$,    $f_\omega $  is continuous  on  $\Omega $  and satisfies the condition  $0 \leq f_\omega \leq 1$. 
Lemma 8.11 tells us that for  $\omega \neq \omega '$  in  $I$,  
we have  $B_\omega \cap B_{\omega '} = \emptyset $.  By Lemma 8.8(d),  if  $z, w \in \Omega $  are such that  
$f_\omega (z) \neq 0$   and   $f_{\omega '}(w) \neq  0$,  then   $d(z,w) \geq 4$.  
It follows that   $f_I$  is continuous on  {\bf B}  and that   $0 \leq f_I \leq 1$.
Furthermore, we can invoke Lemma 7.11 to  obtain   diff$(f_I) \leq \sup _{\omega \in I}\text{diff}(f_\omega ) \leq ((m/13) - 4)^{-1}$,  
where the second  $\leq $  follows from Lemma 8.8(e).  
\medskip
Since  $I \subset I^{(\nu ,\kappa )}_m$,  if  $\omega \in I$,  then    
$B_\omega \subset H_{2^{-\kappa m}} = \{\zeta \in \Omega : -r(\zeta ) < 2^{-\kappa m}\}$.   Since Lemma 8.8(d) says that  
$f_\omega = 0$  on  $\Omega \backslash B_\omega $,  we conclude that  $f_I = 0$  on  
$\{\zeta \in \Omega : -r(\zeta ) \geq  2^{-\kappa m}\}$.  Recalling Definition 7.9,  this completes the verification of the 
membership  $f_I \in \Phi (2^{-m};((m/13) - 4)^{-1})$.   $\square $
\medskip
\noindent
{\bf Lemma 8.13.}   {\it Let}   $m \geq M_{8.7}$,  $\kappa \in \{1,2,3\}$  {\it and}  $\nu \in \{1,\dots ,N_0\}$,  {\it and let}   $I$   
{\it be any subset of}   $I^{(\nu ,\kappa )}_m$.    {\it Then for every bounded operator}    $A$   {\it on}  $L_a^2(\Omega )$,  {\it we have}  
$$
\sum _{\omega \in I}T_{f_\omega }AT_{f_\omega }  \in \text{LOC}(A).
$$ 
\medskip
\noindent
{\it Proof}.   Given any   $I \subset I^{(\nu ,\kappa )}_m$,   consider the set   $\Gamma = \{z_\omega : \omega \in I\}$,  where   
$z_\omega $ was picked in  (8.10).   By Lemmas 8.11 and  8.6,   $\Gamma $  is an  $(m/26)$-separated set in  $\Omega $.
Define  $f_{z_\omega } = f_\omega $  for each  $\omega \in I$.  We need to verify that the functions
$\{f_{z_\omega } : z_\omega \in \Gamma \}$  satisfy conditions (1)-(3) in Definition 6.5.  First of all,  (2)  follows from Lemma 8.8(a). 
Lemma 8.7 tells us that  for each   $\omega \in I$,   we have    $B_\omega \subset D(z_{\omega},R_m)$.  By Lemma 8.8(d),   
we have  $f_{z_\omega } = 0$   on   $\Omega \backslash  D(z_{\omega},R_m)$, verifying (1).   
Finally, condition (3)  follows from Lemma 8.8(c).   $\square $

\bigskip
\centerline{{\bf  9.  The essential commutant of}   $\{T_f : f \in {\text{VO}}_{\text{bdd}}\}$}
\medskip
Recall that we write  ${\Cal K}$  for the collection of compact operators on the Bergman space  $L^2_a(\Omega )$.  
Furthermore,  Proposition 6.2 tells us that   ${\Cal K} \subset {\Cal T}$.   Also recall that for each  $f  \in L^\infty (\Omega )$,  
we have the Hankel operator    $H_f$  defined by the formula
$$
H_fh  =  (1 - P)(fh),   \quad   h \in L_a^2(\Omega ).
$$
\medskip
\noindent
{\it Proof of Theorem} 1.1(i).  Obviously,   Proposition 7.3 implies that
$\text{EssCom}\{T_f : f \in {\text{VO}}_{\text{bdd}}\}$  $\supset $   ${\Cal T}$.  
Thus we only need to prove that  $\text{EssCom}\{T_f : f \in {\text{VO}}_{\text{bdd}}\}$  $\subset $   ${\Cal T}$.
\medskip
Let   $X \in $  EssCom$\{T_f : f \in {\text{VO}}_{\text{bdd}}\}$ be given.   To show that  $X$ $\in $ ${\Cal T}$,  
pick any   $\epsilon > 0$.   It suffices to produce a decomposition   $X = Y + Z$   such that  $Y \in {\Cal T}$  and  
$$
\|Z\|  \leq  3N_0\{16(2 + \|X\|) + \|X\| + 2\}\epsilon ,
\tag 9.1
$$
\noindent
where   $N_0$  is  the constant that appears in (8.7) and (8.15).  
\medskip
First, we apply Proposition 7.10,  which provides a  $\delta > 0$  and  a $t^\ast > 0$   such that
$$
\|[X,T_f]\|  \leq  2\epsilon  \quad \text{for every}  \ \  f \in \Phi (t^\ast ;\delta ).
\tag 9.2
$$
\noindent
Then we apply Lemma 7.2,  which tells us that there is a  $\delta ' > 0$  such that
$$
\|H_g\| \leq \epsilon 
\tag 9.3
$$
\noindent
for every bounded continuous function  $g$  on   $\Omega $  with  diff$(g) \leq \delta '$.
With   $\delta $,   $t^\ast $  and  $\delta '$   so fixed,  we pick an integer  $m \geq M_{8.7}$  satisfying the conditions
$$
((m/13) - 4)^{-1} \leq \min \{\epsilon ,\delta ,\delta '\}    \quad  \text{and}  \quad  2^{-m} \leq t^\ast .
\tag 9.4
$$
\noindent
With  $m$  so fixed,  let us consider the function   $F_{I_m}$  given in Notation 8.10(2).   Since
$$
F_{I_m}  =  \sum _{\kappa =1}^3\sum _{\nu =1}^{N_0}F_{I_m^{(\nu ,\kappa )}} 
\tag 9.5
$$
\noindent
and since by Lemma 8.12 each   $F_{I_m^{(\nu ,\kappa )}}$    satisfies the inequality  
$0 \leq F_{I_m^{(\nu ,\kappa )}} \leq 1$  on  $\Omega $,  we have    $0 \leq F_{I_m} \leq 3N_0$   on  $\Omega $.   
By Lemma 8.8(b)   and   (8.8),   we have   $F_{I_m}(z) \geq 1$   whenever   $-r(z) < 2^{-2m}$.    
Thus we have shown that the function
$$
h  =  \chi _{\Omega _{2^{-2m}}}   +  F_{I_m}
\tag 9.6
$$
\noindent
satisfies the inequality    $1 \leq h \leq 3N_0 + 1$   on   $\Omega $,  where 
$\Omega _{2^{-2m}} = \{\zeta \in \Omega : -r(\zeta ) \geq 2^{-2m}\}$.   This guarantees that the positive Toeplitz operator   $T_h$  
is both bounded and invertible on  $L_a^2(\Omega )$.   Moreover,    $\|T_h^{-1}\| \leq 1$.   Since    $T_h \in {\Cal T}$   and
${\Cal T}$   is a $C^\ast $-algebra,   we have   $T_h^{-1} \in {\Cal T}$.
\medskip
By (9.6) and (9.5),  we have the decomposition
$$
X  =  XT_hT_h^{-1}  =  X_0 +  \sum _{\kappa =1}^3\sum _{\nu =1}^{N_0}X_{\nu ,\kappa },
\tag 9.7
$$
\noindent
where
$$
X_0  =  XT_{\chi _{\Omega _{2^{-2m}}}}T_h^{-1}   \quad  \text{and}  \quad  X_{\nu ,\kappa } =  XT_{F_{I_m^{(\nu ,\kappa )}}}T_h^{-1}
$$
\noindent
for  $1 \leq \kappa \leq 3$   and   $1 \leq \nu \leq N_0$.     
Obviously,  the Toeplitz operator  $T_{\chi _{\Omega _{2^{-2m}}}}$  is compact.    
Hence, by Proposition 6.2,   $X_0 \in {\Cal K} \subset {\Cal T}$.   
\medskip
We further decompose each   $X_{\nu ,\kappa }$.  To do that, define the operators
$$
Y_{\nu ,\kappa }   =  \sum _{\omega \in I_m^{(\nu ,\kappa )}}T_{f_\omega }XT_{f_\omega }T_h^{-1}   \quad   \text{and}  \quad
A_{\nu ,\kappa }   =  \sum \Sb \omega ,\omega ' \in I_m^{(\nu ,\kappa )}\\ \omega \neq \omega ' \endSb T_{f_\omega }XT_{f_{\omega '}}T_h^{-1}.
\tag 9.8
$$
\noindent
Obviously,  $Y_{\nu ,\kappa } + A_{\nu ,\kappa } = T_{f_{I_m^{(\nu,\kappa )}}}XT_{f_{I_m^{(\nu,\kappa )}}}T^{-1}_h$  (cf.  Notation 8.10).   
We further define
$$
B_{\nu ,\kappa } =    [X,T_{f_{I_m^{(\nu ,\kappa )}}}]T_{f_{I_m^{(\nu ,\kappa )}}}T_h^{-1} 
+  XH_{f_{I_m^{(\nu ,\kappa )}}}^\ast H_{f_{I_m^{(\nu ,\kappa )}}}T_h^{-1}.
\tag 9.9
$$
\noindent   
It follows from  Lemmas 8.8(d)  and  8.11 that    $F_{I_m^{(\nu ,\kappa )}} = f_{I_m^{(\nu ,\kappa )}}^2$.  
For any real-valued  $f \in L^\infty (\Omega )$,  we have   $T_{f^2} = T_f^2 + H_f^\ast H_f$.
Therefore    
$$
X_{\nu ,\kappa } =  Y_{\nu ,\kappa } + A_{\nu ,\kappa } + B_{\nu ,\kappa }. 
\tag 9.10
$$
\noindent
Since  $T_h^{-1} \in {\Cal T}$,  it follows from Lemma 8.13  and Corollary 6.7  that  $Y_{\nu ,\kappa } \in {\Cal T}$.
\medskip
To estimate  $\|A_{\nu ,\kappa }\|$,  first observe that on  $L^2(\Omega )$,   we have the strong convergence 
$$
\sum \Sb \omega ,\omega ' \in I_{m,q}^{(\nu ,\kappa )}\\ \omega \neq \omega ' \endSb M_{f_\omega }XPM_{f_{\omega '}}
\rightarrow   \sum \Sb \omega ,\omega ' \in I_m^{(\nu ,\kappa )}\\ \omega \neq \omega ' \endSb M_{f_\omega }XPM_{f_{\omega '}}
\quad \text{as} \quad q \rightarrow \infty ,
$$
\noindent
where  $I_{m,q}^{(\nu ,\kappa )}$  was given by Definition 8.9(b).   
Compressing this strong convergence to   $L_a^2(\Omega )$  and using the bound  $\|T_h^{-1}\| \leq 1$,  
we see that there is a  $q \in $   {\bf N}   such that
$$
\|A_{\nu ,\kappa }\|  \leq  2\|Z_{\nu ,\kappa }\|,  \quad \text{where}  \quad  
Z_{\nu ,\kappa }  =  \sum \Sb \omega ,\omega ' \in I_{m,q}^{(\nu ,\kappa )}\\ \omega \neq \omega ' \endSb T_{f_\omega }XT_{f_{\omega '}}.
\tag 9.11
$$
\noindent
Since   $f_\omega f_{\omega '} = 0$   for   $\omega \neq \omega ' $   in      
$I_{m,q}^{(\nu ,\kappa )}$,  by  Lemma 6.8,  there are complex numbers   $\{\gamma _\omega : \omega \in I_{m,q}^{(\nu ,\kappa )}\}$  
of modulus  $1$  and a subset   $I$  of  $I_{m,q}^{(\nu ,\kappa )}$ such that if we define
$$
F = \sum _{\omega \in I}f_\omega ,  \quad   G = \sum _{\omega \in I_{m,q}^{(\nu ,\kappa )}\backslash I}f_\omega ,  \quad 
F' = \sum _{\omega \in I}\gamma _\omega f_\omega    \quad \text{and}  \ \  
G' = \sum _{\omega \in I_{m,q}^{(\nu ,\kappa )}\backslash I}\gamma _\omega f_\omega , 
$$
\noindent
then
$$
\|Z_{\nu ,\kappa }\|   \leq   4(\|T_{F'}XT_G\| + \|T_{G'}XT_F\|).
\tag 9.12
$$
\noindent
Note that   $T_{G'}XT_F = T_{G'}[X,T_F] + T_{G'}T_FX$.  We have  
$F \in \Phi (2^{-m};((m/13) - 4)^{-1})$  by Lemma 8.12.   Hence it follows from (9.4) and (9.2)  that   
$$
\|T_{G'}[X,T_F]\|  \leq  \|[X,T_F]\| \leq 2\epsilon .
\tag 9.13
$$
\noindent
Since  $B_\omega \cap B_{\omega '} = \emptyset $  for all   $\omega  \neq \omega '$   in  $I_m^{(\nu ,\kappa )}$,  we have   
$G'F = 0$  on  $\Omega $,   and consequently    $T_{G'}T_F$  $=$  $-H_{\overline{G'}}^\ast H_F$.    
Since $\text{diff}(F) \leq ((m/13) - 4)^{-1}$,   by (9.4)  and (9.3),  we have
$$
\|T_{G'}T_FX\|  \leq  \|H_F\|\|X\|  \leq    \|X\|\epsilon .
$$
\noindent
Combining this with (9.13),   we see that   $\|T_{G'}XT_F\| \leq (2 + \|X\|)\epsilon $.  The same argument also shows that  
$\|T_{F'}XT_G\| \leq (2 + \|X\|)\epsilon $.   Substituting these  in (9.12) and recalling (9.11),  we obtain
$$
\|A_{\nu ,\kappa }\|  \leq  16(2 + \|X\|)\epsilon .
\tag 9.14
$$
\noindent
Next we estimate  $\|B_{\nu ,\kappa }\|$.
\medskip 
Lemma 8.12 tells us that  $\text{diff}(f_{I_m^{(\nu ,\kappa )}})  \leq ((m/13) - 4)^{-1}$.   Combining this with (9.4) and (9.3),  
and with the fact  $\|T_h^{-1}\| \leq 1$,   we obtain    
$$
\|XH_{f_{I_m^{(\nu ,\kappa )}}}^\ast H_{f_{I_m^{(\nu ,\kappa )}}}T_h^{-1}\|  \leq   \|X\|\|H_{f_{I_m^{(\nu ,\kappa )}}}\|
\leq  \|X\|\epsilon .
$$
\noindent
Again, Lemma 8.12 says  that   $f_{I_m^{(\nu ,\kappa )}} \in \Phi (2^{-m};((m/13) - 4)^{-1})$.
Hence it follows from (9.4) and (9.2)  that   
$$
\|[X,T_{f_{I_m^{(\nu ,\kappa )}}}]T_{f_{I_m^{(\nu ,\kappa )}}}T_h^{-1}\|   \leq  \| [X,T_{f_{I_m^{(\nu ,\kappa )}}}]\|  \leq 2\epsilon .
$$
\noindent
Recalling (9.9),  from the above two inequalities we obtain
$$
\|B_{\nu ,\kappa }\|  \leq  (\|X\| + 2)\epsilon .
\tag 9.15
$$
\noindent
To summarize,  we have shown that for each pair of $1 \leq \kappa \leq 3$   and   $1 \leq \nu \leq N_0$,  
we have the decomposition   (9.10)
where  $Y_{\nu ,\kappa } \in {\Cal T}$   and  where  $A_{\nu ,\kappa }$,  $B_{\nu ,\kappa }$   satisfy estimates  
(9.14)  and  (9.15)  respectively.   Combining (9.10)  with (9.7),  we have   $X = Y + Z$,    where
$$
Y = X_0 + \sum _{\kappa =1}^3\sum _{\nu =1}^{N_0}Y_{\nu ,\kappa }  \quad \text{and} \quad  
Z = \sum _{\kappa =1}^3\sum _{\nu =1}^{N_0}(A_{\nu ,\kappa } + B_{\nu ,\kappa }).
\tag 9.16
$$
\noindent
Now, (9.1) follows from (9.14)  and  (9.15),  and we have shown that    $Y  \in {\Cal T}$.   This completes the 
proof  of part (i) in Theorem 1.1.   $\square $ 

\medskip
\noindent
{\bf Proposition 9.1.}    {\it For}  $X \in {\Cal T}$,    {\it if}  LOC$(X) \subset {\Cal K}$,  {\it then}  $X$  {\it is compact}.
\medskip
\noindent
{\it Proof}.    Let  $X \in {\Cal T}$  and suppose that  LOC$(X) \subset {\Cal K}$.  As we showed above,  for every  
$\epsilon > 0$,  $X$  admits a decomposition  $X = Y + Z$,  where  $Y$  and  $Z$  are given by (9.16),  
with  $X_0$  known to be compact.    
By (9.8) and Lemma 8.13,  the condition    LOC$(X) \subset {\Cal K}$  implies   $Y_{\nu ,\kappa } \in {\Cal K}$.   
Thus   $Y$  is compact.  Since  $Z$  satisfies (9.1),  this shows that   $X$  is compact.    $\square $
\medskip
\noindent
{\bf Proposition 9.2.}    {\it  Let}  $X \in {\Cal T}$.       {\it Suppose that}  $X$  {\it has the property that for every}  $0 < R < \infty $,  
$$
\lim _{z \rightarrow \partial \Omega }\sup \{|\langle Xk_w,k_z\rangle |  :  d(z,w) < R\}  =  0.
\tag 9.17
$$
\noindent
{\it Then}   $X$  {\it is a compact operator}.
\medskip
\noindent
{\it Proof}.     Recall from Proposition 6.6 that     LOC$(X) \subset {\Cal D}(X)$.   
Combining this with  Proposition 9.1,  it suffices to  prove the inclusion  
${\Cal D}_0(X) \subset {\Cal K}$ under the assumption that (9.17) holds for every   $0 < R < \infty $.  
By  Definition 6.3(c),  we need to show that the operator   
$$
T =  \sum _{z\in \Gamma }c_z\langle Xk_{\psi (z)},k_{\varphi (z)}\rangle k_{\varphi (z)}\otimes k_{\psi (z)}
$$  
\noindent
is compact,  where  $\Gamma $ is a separated set in  $\Omega $,   $\{c_z : z \in \Gamma \}$  is a bounded set of coefficients,   
and  $\varphi , \psi  :  \Gamma \rightarrow \Omega $  are maps for which there is a  $0 \leq C < \infty $  
such that    $d(z,\varphi (z)) \leq C$  and    $d(z,\psi (z)) \leq C$  for every   $z \in \Gamma $.
\medskip
By the assumption on  $\varphi $,  $\psi $   and Lemma 2.11,  
there is a partition   $\Gamma = \Gamma _1\cup \cdots \cup \Gamma _k$   such that  
for each  $1 \leq j \leq k$,   the conditions   $z, w \in \Gamma _j$  and  $z \neq w$  imply  
$d(\varphi (z),\varphi (w)) > 2$   and    $d(\psi (z),\psi (w)) > 2$.  Hence for each  $1 \leq j \leq k$,   the sets  
$\{\varphi (z) : z \in \Gamma _j\}$   and    $\{\psi (z) : z \in \Gamma _j\}$  are  $1$-separated.  
This leads to the decomposition   $T = T_1 + \cdots + T_k$,  where  
$$
T_j =  \sum _{z\in \Gamma _j}c_z\langle Xk_{\psi (z)},k_{\varphi (z)}\rangle k_{\varphi (z)}\otimes k_{\psi (z)}
$$
\noindent
for every  $1 \leq j \leq k$.   Thus it suffices to show that   $T_j \in {\Cal K}$     for every  $1 \leq j \leq k$.   
Fix such a  $j$  for the moment.    For each     $\delta > 0$,  denote   $\Gamma _{j,\delta } = \{z \in \Gamma _j : -r(z) \leq \delta \}$.   
Using an obvious finite-rank approximation  and applying Lemma 5.1,    for each  $\delta > 0$,   we have
$$
\|T_j\|_{\Cal Q}  \leq  \bigg\|\sum _{z\in \Gamma _{j,\delta }}c_z\langle Xk_{\psi (z)},k_{\varphi (z)}\rangle k_{\varphi (z)}\otimes k_{\psi (z)}\bigg\|
\leq C_{5.1}^2c\sup _{z\in \Gamma _{j,\delta }}|\langle Xk_{\psi (z)},k_{\varphi (z)}\rangle |,
$$
\noindent
where  $c = \sup _{z\in \Gamma }|c_z|$.   Since  $d(z,\varphi (z)) \leq C$  and    
$d(z,\psi (z)) \leq C$  for every   $z \in \Gamma $,  it follows from (9.17) that the right-hand side tends to  $0$  as  $\delta \downarrow 0$.   
Thus  $\|T_j\|_{\Cal Q} = 0$,   i.e.,  $T_j$  is a compact operator.    This completes the proof.    $\square $
\medskip
As an immediate consequence of Proposition 9.2,  we have
\medskip
\noindent
{\bf Corollary 9.3.}    {\it  Let}  $X \in {\Cal T}$.       {\it Then}  $X$  {\it is compact if and only if}  
$$
\lim _{z \rightarrow \partial \Omega }\|Xk_z\|   =  0.
$$

\bigskip
\centerline{\bf  10.  The essential commutant of the Toeplitz algebra} 
\medskip
We now turn to the proof of part (ii) in Theorem 1.1.
\medskip
\noindent
{\bf Proposition 10.1.}   {\it  If}  $f \in {\text{VO}}_{\text{bdd}}$,   {\it  then}
$$
\lim _{z\rightarrow \partial \Omega }\|(f - f(z))k_z\|  =  0.
$$
\medskip
\noindent
{\it Proof}.    Let   $f \in {\text{VO}}_{\text{bdd}}$   and consider  a large   $R > 0$.   We have  
$$
\align 
\|(f - f(z))k_z\|^2  &=  \int _{D(z,R)}|f(w) - f(z)|^2|k_z(w)|^2dv(w)       \\
&\quad \quad +     \int _{\Omega \backslash D(z,R)}|f(w) - f(z)|^2|k_z(w)|^2dv(w)  \\
&\leq  \sup _{d(z,w) < R}|f(w) - f(z)|^2  +  C_1\|f\|_\infty ^2\int _{\Omega \backslash D(z,R)}{|r(z)|^{n+1}\over F(z,w)^{2n+2}}dv(w).
\endalign
$$
\noindent
Applying Lemma 3.8,  there are constants   $0 < C_2 < \infty $  and   $s > 0$  such that
$$
\|(f - f(z))k_z\|^2 \leq  \sup _{d(z,w) < R}|f(w) - f(z)|^2 + C_2\|f\|_\infty ^22^{-sR}.
\tag 10.1
$$
\noindent
Now we use the fact that  $f$  has vanishing oscillation:  
Using the cutoff functions provided by Lemma 7.6,  for any   $\delta > 0$,  we can write  $f = f_1 + f_2$,  where   $f_1$  has a compact 
support in  $\Omega $  and    diff$(f_2) \leq \delta $.   
Combining this fact with Lemma 7.1,  we see that 
$$
\lim _{z\rightarrow \partial \Omega }\sup _{d(z,w) < R}|f(w) - f(z)|  =  0
$$
\noindent
once an  $R > 0$  is given.  This and (10.1) together imply that   $\|(f - f(z))k_z\|  \rightarrow 0$  as  $z \rightarrow \partial \Omega $.   
This completes the proof.    $\square $
\medskip
\noindent
{\bf Proposition 10.2.}   {\it Suppose that}  $\{z_j\}$   {\it and}  $\{w_j\}$   {\it are sequences in}   $\Omega $   
{\it satisfying the following two conditions}:

(1)  $\lim _{j\rightarrow \infty }r(z_j)  =  0$.

(2)   {\it There is a constant}    $0 < C < \infty $  {\it such that}     $d(z_j,w_j)  \leq  C$   {\it for every}   $j \in $  {\bf N}.    

\noindent
{\it Then for every}    $A \in \text{EssCom}({\Cal T})$   {\it we have}
$$
\lim _{j\rightarrow \infty }\|[A,k_{z_j}\otimes k_{w_j}]\|  =   0.
\tag 10.2
$$ 
\medskip
\noindent
{\it Proof}.   For the given   $\{z_j\}$,  $\{w_j\}$  and  $A$,   suppose that (10.2) did not hold.    
Then,  replacing   $\{z_j\}$,  $\{w_j\}$  by subsequences  if necessary,  
we may assume that there is a   $c > 0$  such that
$$
\lim _{j\rightarrow \infty }\|[A,k_{z_j}\otimes k_{w_j}]\|  =   c.
\tag 10.3
$$ 
\noindent
We will show that this leads to a contradiction.
\medskip
By condition (1) and Lemma 2.1,  there is a sequence  $j_1 < j_2 <  \cdots  < j_\nu < \cdots $   of natural numbers such that 
$-r(z_{j_{\nu +1}}) < -r(z_{j_\nu } )$   for every   $\nu \in $  {\bf N} and such that the set    
$\{z_{j_\nu }  :  \nu  \in {\bold N}\}$     is $1$-separated.   For each   $\nu  \in $  {\bf N},  we now define the operator
$$
B_\nu   =  [A,k_{z_{j_\nu }}\otimes k_{w_{j_\nu }}],
$$
\noindent
whose rank is at most  $2$.  By conditions (1), (2)  and Lemma 2.1,  we also have
that   $r(w_j) \rightarrow 0$  as  $j \rightarrow \infty $.   Thus both sequences of vectors   $\{k_{z_j}\}$  and  
$\{k_{w_j}\}$   converge to  $0$  weakly in   $L_a^2(\Omega )$.   Consequently we have the convergence  
$$
\lim _{\nu \rightarrow \infty }B_\nu  =   0   \quad \text{and} \quad  
\lim _{\nu \rightarrow \infty }B_\nu ^\ast  =   0
$$
\noindent
in the strong operator topology.   Thus by (10.3) and  Lemma 7.8,   there is a subsequence  
$\nu (1) < \nu (2) <  \cdots < \nu (m) < \cdots $  of natural numbers such that  the sum  
$$
B  =  \sum _{m=1}^\infty  B_{\nu (m)}
$$
\noindent
converges strongly with  $\|B\|_{\Cal Q} = c > 0$.  Thus  $B$  is not compact.
Now define the operator  
$$
Y   =  \sum _{m=1}^\infty  k_{z_{j_{\nu (m)}}}\otimes k_{w_{j_{\nu (m)}}}.
$$
\noindent
Since the set    $\{z_{j_\nu }  :  \nu  \in {\bold N}\}$   is $1$-separated and since condition (2) holds, by Proposition 6.4  we have   
$Y \in {\Cal T}$.    Since   $A \in  \text{EssCom}({\Cal T})$,    the commutator   $[A,Y]$  is compact.   On the other hand, we clearly 
have  $[A,Y] = B$,  which is not compact because  $\|B\|_{\Cal Q}  > 0$.    This gives us the contradiction promised earlier.   $\square $
\medskip
\noindent
{\bf Lemma 10.3.}   [27,Lemma 5.1]   {\it  Let}  $T$  {\it be a bounded},  {\it self-adjoint operator on a Hilbert space}   ${\Cal H}$.  
{\it Then for each unit vector}   $x \in {\Cal H}$   {\it we have}   $\|[T,x\otimes x]\|  =   \|(T - \langle Tx,x\rangle )x\|$.
\medskip
\noindent
{\bf Lemma 10.4.}    [27,Lemma 5.2]  {\it  Let}  $T$  {\it be a bounded},  {\it self-adjoint operator on a Hilbert space}   ${\Cal H}$.  
{\it Then for every pair of unit vectors}   $x, y \in {\Cal H}$   {\it we have}
$$
|\langle Tx,x\rangle - \langle Ty,y\rangle |  \leq  \|[T,x\otimes y]\|  +   \|[T,x\otimes x]\| +  \|[T,y\otimes y]\|.
$$
\medskip
For a bounded operator    $A$  on  $L^2_a(\Omega )$,  we define the function  
$$
\tilde A(z)  =  \langle Ak_z,k_z\rangle ,   \quad  z \in \Omega .
$$
\noindent
Recall that   $\tilde A$  is commonly called the  {\it Berezin transform}  of the operator    $A$.  
\medskip
\noindent
{\bf Proposition 10.5.}    {\it  If}    $A \in \text{EssCom}({\Cal T})$,   
{\it then its Berezin transform}   $\tilde A$  {\it is in}    $\text{VO}_{\text{bdd}}$.
\medskip
\noindent
{\it Proof}.   It suffices to consider a self-adjoint   $A \in \text{EssCom}({\Cal T})$.  
Obviously,   $\tilde A$  is bounded, and Proposition 4.6 tells us that it is continuous on   $\Omega $.
If it were true that  $\tilde A \notin \text{VO}$,  
then there would be a  $c > 0$   and sequences    $\{z_j\}$,  $\{w_j\}$   in  $\Omega $  with
$$
\lim _{j\rightarrow \infty }r(z_j)  =  0
\tag 10.4
$$
\noindent
such that for every   $j \in $  {\bf N},   we have   $d(z_j,w_j) \leq 1$   and
$$
|\langle Ak_{z_j},k_{z_j}\rangle  -  \langle Ak_{w_j},k_{w_j}\rangle |  
=  |\tilde A(z_j) -  \tilde A(w_j)|  \geq c.
\tag 10.5
$$
\noindent
But on the other hand, it follows from Lemma 10.4  that
$$
|\langle Ak_{z_j},k_{z_j}\rangle  -  \langle Ak_{w_j},k_{w_j}\rangle |    
\leq \|[A,k_{z_j}\otimes k_{w_j}]\| + \|[A,k_{z_j}\otimes k_{z_j}]\| + \|[A,k_{w_j}\otimes k_{w_j}]\|.
\tag 10.6
$$
\noindent
By (10.4)  and the condition  $d(z_j,w_j) \leq 1$,  $j\in $ {\bf N},  we can apply Proposition 10.2 to obtain
$$
\lim _{j\rightarrow \infty }\|[A,k_{z_j}\otimes k_{w_j}]\| = 0  \quad \text{and} \quad  
\lim _{j\rightarrow \infty }\|[A,k_{z_j}\otimes k_{z_j}]\| = 0.
\tag 10.7
$$
\noindent
By Lemma 2.1,  conditions  (10.4)  and   $d(z_j,w_j) \leq 1$,  $j\in $ {\bf N},  also imply   
$\lim _{j\rightarrow \infty }r(w_j) = 0$.   Thus Proposition 10.2 also provides that
$$
\lim _{j\rightarrow \infty }\|[A,k_{w_j}\otimes k_{w_j}]\| = 0.
\tag 10.8
$$
\noindent
Obviously,   (10.6), (10.7) and (10.8) together contradict (10.5).    $\square $
\medskip
\noindent
{\bf Lemma 10.6.}    {\it  If}    $A \in \text{EssCom}({\Cal T})$,   {\it then} 
$$
\lim _{z\rightarrow \partial \Omega }\|(A -  T_{\tilde A})k_z\|  =  0.
$$
\medskip
\noindent
{\it Proof}.    Again, it suffices to consider a self-adjoint  $A \in \text{EssCom}({\Cal T})$.     
Then  it follows from Lemma 10.3 and  Proposition 10.2  that
$$
\lim _{z\rightarrow \partial \Omega }\|(A -  \tilde A(z))k_z\|  =
\lim _{z\rightarrow \partial \Omega }\|[A,k_z\otimes k_z]\|  =  0.
$$
\noindent
Therefore it suffices to show that
$$
\lim _{z\rightarrow \partial \Omega }\|(T_{\tilde A} - \tilde A(z))k_z\|  =  0.
$$
\noindent
Since  $\|(T_{\tilde A} - \tilde A(z))k_z\| \leq \|(\tilde A - \tilde A(z))k_z\|$,  this follows from  Propositions 10.5 and 10.1.     
$\square $
\medskip
Finally, we are ready to determine the essential commutant of   ${\Cal T}$.
\medskip
\noindent
{\it Proof of Theorem}  1.1(ii).  Again,  it follows from  Proposition 7.3   that
EssCom$({\Cal T})$ $\supset $  $\{T_f : f \in {\text{VO}}_{\text{bdd}}\} + {\Cal K}$.   
\medskip
For the reverse inclusion,  consider any  $A \in $  EssCom$({\Cal T})$.    We need to show that  $A$  $\in $  
$\{T_f : f \in {\text{VO}}_{\text{bdd}}\} + {\Cal K}$.   We know that  $\tilde A \in {\text{VO}}_{\text{bdd}}$   from Proposition 10.5.   
Hence it suffices to show that  $A - T_{\tilde A}$  is compact.   For this we apply Lemma 10.6,  which gives us
$$
\lim _{z\rightarrow \partial \Omega }\|(A -  T_{\tilde A})k_z\|  =  0.
\tag 10.9
$$
\noindent
The membership  $A \in $  EssCom$({\Cal T})$   implies, of course, that  $A \in $  EssCom$\{T_f : f \in {\text{VO}}_{\text{bdd}}\}$.  Hence  
Theorem 1.1(i)  tells us that    $A \in {\Cal T}$.  Consequently,     $A - T_{\tilde A} \in {\Cal T}$.   By Corollary 9.3,   the membership 
 $A - T_{\tilde A} \in {\Cal T}$  and   (10.9)  together  imply that    $A - T_{\tilde A}$  is compact.    $\square $

\bigskip
\centerline{\bf  11.  Berezin transform near the boundary} 
\medskip
The purpose of this section is to show that condition (9.17) is implied by the vanishing of Berezin 
transform near  $\partial \Omega $.   This along with Proposition 9.2 will give us the proof of Theorem 1.2.  
To begin, we need to fix some necessary constants:
\medskip
\noindent
{\bf Lemma 11.1.}    (1)   {\it  There is a}  $0 < c_0  < 1$    {\it  such that}   
$z + {\Cal P}((\bar \partial r)(z);2c_0\sqrt{-r(z)},-2c_0r(z)) \subset D(z,1)$  {\it  for every}  
$z \in \Omega $  {\it  satisfying the condition}   $-r(z) < \theta $.

\noindent
(2)   {\it  There is a}  $b_0 > 0$  {\it  such that}   $D(z,3b_0) \subset z + {\Cal P}((\bar \partial r)(z);c_0\sqrt{-r(z)},-c_0r(z))$
{\it  for every}  $z \in \Omega $  {\it  satisfying the condition}    $-r(z) < \theta $.

\noindent
(3)   {\it  There is an}  $a_0 > 0$   {\it  such that}   $z + {\Cal P}((\bar \partial r)(z);a_0\sqrt{-r(z)},-a_0r(z)) \subset D(z,b_0)$
{\it  for every}  $z \in \Omega $  {\it  satisfying the condition}    $-r(z) < \theta $.

\medskip
\noindent
{\it Proof}.   By Proposition 2.4,  there is a  $0 < c < 1$  such that   $z + {\Cal P}((\bar \partial r)(z);c\sqrt{-r(z)},-cr(z)) \subset D(z,1)$  
for every  $z \in \Omega $  satisfying the condition   $-r(z) < \theta $.   Then $c_0 = c/2$  will do for (1).
\medskip
To prove (2),  take any   $0 < b < 1/2$  such that   $C_{2.5}b < c_0$.  By Proposition 2.5,  we have
$$
D(z,b)  \subset  z + {\Cal P}((\bar \partial r)(z);c_0\sqrt{-r(z)},-c_0r(z)) 
$$
\noindent
whenever   $-r(z) < \theta $.
Thus   (2) holds for the constant  $b_0 = b/3$.
\medskip
Finally, note that  (3) is a direct consequence of Proposition 2.4.   $\square $
\medskip
Once the above constants are fixed,  we can introduce the following  ``polyballs":
\medskip
\noindent
{\bf Definition 11.2.}   (1)  Let  
$$
\align
{\Cal P}  &=  \{(u_1,u_2,\dots ,u_u) \in {\bold C}^n :   |u_1| < a_0  \ \text{and} \   (|u_2|^2 + \cdots + |u_n|^2)^{1/2} < a_0\},   \\
{\Cal Q}  &=  \{(u_1,u_2\dots ,u_u) \in {\bold C}^n :   |u_1| \leq c_0  \ \text{and} \   (|u_2|^2 + \cdots + |u_n|^2)^{1/2} \leq c_0\} 
\quad  \text{and}    \\
{\Cal R}  &=  \{(u_1,u_2\dots ,u_u) \in {\bold C}^n :   |u_1|  < 2c_0  \ \text{and} \   (|u_2|^2 + \cdots + |u_n|^2)^{1/2} < 2c_0\}. 
\endalign
$$
\noindent
(2)   For each    $z \in \Omega $  satisfying the condition   $-r(z) < \theta $,  let  ${\Cal S}_z$  be the linear transformation  
on    ${\bold C}^n$   given by the formula  
$$
{\Cal S}_z(u_1,u_2,\dots u_n)  =  (-r(z)u_1,\sqrt{-r(z)}u_2,\dots ,\sqrt{-r(z)}u_n),   \quad  (u_1,u_2,\dots u_n) \in {\bold C}^n.
$$
\noindent
(3)   For each    $z \in \Omega $  satisfying the condition   $-r(z) < \theta $,  let  ${\Cal U}_z$   be a unitary transformation on  
${\bold C}^n$   such that    ${\Cal U}_z\{(0,u_2,\dots ,u_n) : u_2, \dots ,u_n \in {\bold C}\}$ $=$ 
$\{u \in {\bold C}^n : \langle u,(\bar \partial r)(z)\rangle = 0\}$.

\noindent
(4)   For each    $z \in \Omega $  satisfying the condition   $-r(z) < \theta $,  denote  ${\Cal V}_z = {\Cal U}_z{\Cal S}_z$.
\medskip
\noindent
{\bf Proposition 11.3.}  {\it  Suppose that}   $U$  {\it  is a connected open set in}  ${\bold C}^n$  
{\it  that is symmetric with respect to conjugation}.   {\it  That is},    $(w_1,\dots ,w_n) \in U$  
{\it  if and only if}   $(\bar w_1,\dots ,\bar w_n) \in U$.   {\it  Let}  $F$   {\it  be an analytic function on the domain}  
$U\times U$   {\it  in}   ${\bold C}^n\times {\bold C}^n$.   {\it  If}   $F(\bar z,z) = 0$  {\it  for every}  
$z \in U$,  {\it  then}   $F$  {\it  is identically zero on}  $U\times U$.
\medskip
\noindent
{\it Proof}.   For each   $j \in \{1,\dots ,n\}$,  let  $e_j$   denote the vector in  ${\bold C}^n$  whose   $j$-th component 
is  $1$   and whose other components are  $0$.    We then define
$$
\align  
(d_jF)(w,z)   &=  {1\over 2}\bigg({\partial \over \partial x} + i{\partial \over \partial y}\bigg)
F(w + \overline{(x+iy)e_j},z+(x+iy)e_j)\bigg|_{x=0=y}   \quad \text{and}  \\
(\partial _jF)(w,z)  &= {1\over 2}\bigg({\partial \over \partial x} - i{\partial \over \partial y}\bigg)
F(w+(x+iy)e_j,z)\bigg|_{x=0=y}
\endalign
$$
\noindent
for    $j \in \{1,\dots ,n\}$  and  $w, z \in U$.   It is straightforward to verify that for every multi-index  
$\alpha \in {\bold Z}_+^n$,   we have   $d^\alpha F = \partial ^\alpha F$.   
Since    $F(\bar z,z) = 0$  for every  $z \in U$,   an easy induction on  $|\alpha |$  yields  
$(d^\alpha F)(\bar z,z) = 0$  for every  $z \in U$  and every  $\alpha \in {\bold Z}_+^n$.  Thus if we fix any 
$z \in U$,  then   $(\partial ^\alpha F)(\bar z,z) = 0$  for every  $\alpha \in {\bold Z}_+^n$.   By the standard power-series expansion,  
this means that the analytic function  $f_z(\zeta ) =  F(\zeta ,z)$,  $\zeta \in U$,    vanishes on a small open ball 
containing  $\bar z$.   Since    $U$  is connected,  we conclude that   $f_z = 0$  on    $U$.  
Since this is true for every   $z \in U$,  it follows that  $F$  is identically zero on   $U\times U$.   $\square $
\medskip
\noindent
{\bf Proposition 11.4.}   {\it Let}  $A$  {\it be a bounded operator on the Bergman space}  $L^2_a(\Omega )$.  {\it If}
$$
\lim _{z\rightarrow \partial \Omega }\langle Ak_z,k_z\rangle = 0,
\tag 11.1
$$
\noindent
{\it then for every given}  $0 < R < \infty $  {\it we have}
$$
\lim _{z\rightarrow \partial \Omega }\sup \{|\langle Ak_w,k_z\rangle | :  w \in D(z,R)\}  =   0.
\tag 11.2
$$
\noindent 
\medskip
\noindent
{\it Proof}.   Given (11.1),  suppose that (11.2) failed for some   $0 < R < \infty $.   We will show that this results in a 
contradiction.  First of all,  the failure of (11.2)  for this particular  $R$  means that  there is an  $\epsilon > 0$   and sequences  
$\{z_j\}$,  $\{w_j\}$   in  $\Omega $   such that   
$$
\lim _{j\rightarrow \infty }r(z_j)   =  0
\tag 11.3
$$
\noindent
and at the same time,   $d(z_j,w_j) < R$  and  
$$
|\langle Ak_{w_j},k_{z_j}\rangle |  \geq  \epsilon 
\tag 11.4
$$
\noindent
for every  $j \geq 1$.   Since    $d(z_j,w_j) < R$,   for every  $j \geq 1$  we have a   $C^1$  map  
$g_j : [0,1] \rightarrow \Omega $  such that   $g_j(0) = z_j$,  $g_j(1) = w_j$   and
$$
\int _0^1\sqrt{\langle {\Cal B}(g_j(t))g_j'(t),g_j'(t)\rangle }dt    \leq  R+1.
\tag 11.5
$$
\noindent
By (11.3), (11.5)  and Lemma 2.1,  discarding a finite number of  $j$'s if necessary,  we may assume that  
$-r(g_j(t)) < \theta $  for all   $j$   and   $t \in [0,1]$.   Thus  Lemma 11.1 can be applied on all these paths.  
With the   $b_0$  provided by Lemma 11.1,  we pick an  $m \in $  {\bf N}   such that   $(R+1)/m < b_0$.   
Thus for  every  $j \geq 1$,  there is a partition  
$$
0 = x_j(0) < x_j(1) < \cdots < x_j(m) = 1
$$
\noindent
of the interval   $[0,1]$   such that
$$
\int _{x_j(i-1)}^{x_j(i)}\sqrt{\langle {\Cal B}(g_j(t))g_j'(t),g_j'(t)\rangle }dt  \leq  {R+1\over m} < b_0
\tag 11.6
$$
\noindent
for every  $1 \leq i \leq m$.   Now,  for every pair of  $j \geq 1$  and  $0 \leq i \leq m$,  we define
$$
z_j^{(i)}  =   g_j(x_j(i)).
$$
\noindent
In particular,  we have  $z_j^{(0)} = z_j$   and   $z_j^{(m)} = w_j$  for all   $j$.
\medskip
Recall that we write   $K_z(\zeta ) = K(\zeta ,z)$,   which is
the  (unnormalized) reproducing kernel for  $L^2_a(\Omega )$.   Let us denote
$$
\Phi (w,z)  =  \langle AK_w,K_z\rangle ,
$$
\noindent
$w, z \in \Omega $.   For every pair of  $j \geq 1$  and  $0 \leq i \leq m$,  we define the function
$$
F_j^{(i)}(\zeta ,\xi )  =  |r(z_j^{(0)})|^{n+1}\Phi \left(z_j^{(i)} + {\Cal V}_{z_j^{(i)}}\overline{\zeta },z_j^{(0)}+{\Cal V}_{z_j^{(0)}}\xi \right),
\tag 11.7
$$
\noindent
$ \zeta ,\xi \in {\Cal R}$.     A review of Definitions 11.2 and 2.3  gives us the identity      
$$
z_j^{(i)} + {\Cal V}_{z_j^{(i)}}{\Cal R}
=  z_j^{(i)} + {\Cal P}((\bar \partial r)(z_j^{(i)});2c_0\sqrt{-r(z_j^{(i)})},-2c_0r(z_j^{(i)})).
\tag 11.8
$$
\noindent
Therefore Lemma 11.1 ensures that each    $F_j^{(i)}$  is well defined,  and it is obviously an analytic function of  
${\Cal R}\times {\Cal R}$.     By  (11.8), Lemma 11.1(1) and  (11.5),  if   $w = z_j^{(i)} + {\Cal V}_{z_j^{(i)}}\overline{\zeta }$  for some  
$\zeta \in {\Cal R}$,  then   $d(w,z_j^{(0)}) \leq R+2$.   Thus by (4.1) and Lemma 2.1,   there is a  $C_1 = C_1(R)$   such that   
$$
|F_j^{(i)}(\zeta ,\xi )|  \leq  C_1\|A\|
$$
\noindent
for all   $\zeta , \xi \in {\Cal R}$,  $j \geq 1$   and  $0 \leq i \leq m$.   Hence  for each   $0 \leq i \leq m$,   $\{F_j^{(i)} : j \geq 1\}$  is a normal  
family of analytic functions on   ${\Cal R}\times {\Cal R}$.   Consequently there is a sequence 
$$
j_1 < j_2 <  \cdots < j_\nu < \cdots   
$$
\noindent
in  {\bf N}  such that for every   $0 \leq i \leq m$,   the sequence   $\{F_{j_\nu }^{(i)}\}_{\nu \in {\bold N}}$  is uniformly convergent on each   
compact subset of      ${\Cal R}\times {\Cal R}$.     For every   $0 \leq i \leq m$,   define the function
$$
F^{(i)}  =  \lim _{\nu \rightarrow \infty }F_{j_\nu }^{(i)}  
\tag 11.9
$$
\noindent
on   ${\Cal R}\times {\Cal R}$.   Next we show that every    $F^{(i)}$  is identically zero on   ${\Cal R}\times {\Cal R}$.      
\medskip
We will accomplish this by an induction on  $i$.   First,  let us show that   $F^{(0)}$  is the  zero function.  
For  $j  \geq 1$  and    $\zeta \in {\Cal R}$,  we have
$$
F_j^{(0)}(\overline{\zeta },\zeta )  =  |r(z_j^{(0)})|^{n+1}\Phi \left(z_j^{(0)} + {\Cal V}_{z_j^{(0)}}\zeta ,z_j^{(0)}+{\Cal V}_{z_j^{(0)}}\zeta \right).
$$
\noindent
As we explained above,  (4.1)  and  Lemma 2.1 together guarantee that
$$
|F_j^{(0)}(\overline{\zeta },\zeta )|  
\leq  C_1\bigg|\bigg\langle Ak_{z_j^{(0)} + {\Cal V}_{z_j^{(0)}}\zeta },k_{z_j^{(0)} + {\Cal V}_{z_j^{(0)}}\zeta }\bigg\rangle \bigg|.
$$
\noindent
By (11.3)  and   Lemmas 11.1(1)  and  2.1,      for each  $\zeta \in {\Cal R}$  we have  
$r\big (z_j^{(0)} + {\Cal V}_{z_j^{(0)}}\zeta \big)  \rightarrow 0$  as  $j \rightarrow \infty $.
Thus, combining the above inequality with    (11.1) and (11.9),  we find that  $F^{(0)}(\overline{\zeta },\zeta )$  $=$  $0$   
for every  $\zeta  \in {\Cal R}$.    By Proposition 11.3,  $F^{(0)}$  is identically zero on   ${\Cal R}\times {\Cal R}$.
\medskip
Now suppose that  $0 \leq i < m$  and that we have shown that   $F^{(i)}$  is identically zero on   ${\Cal R}\times {\Cal R}$.   
We need to show that    $F^{(i+1)}$  is also identically zero on   ${\Cal R}\times {\Cal R}$.   By (11.6),  we have  
$d(z_j^{(i)},z_j^{(i+1)}) < b_0$.  A review of Definition 11.2 and  Lemma 11.1  gives us
$$
z_j^{(i+1)} + {\Cal V}_{z_j^{(i+1)}}{\Cal P}  \subset  D(z_j^{(i+1)},b_0)  \subset D(z_j^{(i)},3b_0) \subset z_j^{(i)} + {\Cal V}_{z_j^{(i)}}{\Cal Q}.
\tag 11.10
$$
\noindent
Let  $\xi \in {\Cal R}$  be given.     By (11.7) and (11.10),  for any  $\zeta \in {\Cal P}$,  there is an  $\eta _j(\zeta ) \in {\Cal Q}$  such that
$$
F_j^{(i+1)}(\zeta ,\xi )  =  F_j^{(i)}(\eta _j(\zeta ),\xi ).
\tag 11.11
$$
\noindent
Since  ${\Cal Q}$  is a compact set in  ${\Cal R}$  and since  $F^{(i)} = 0$,  by  (11.9)  we have
$$
\lim _{\nu \rightarrow \infty }\sup \{|F_{j_\nu }^{(i)}(\eta ,\xi )| : \eta \in {\Cal Q}\}  =  0.
$$
\noindent
Combining this with (11.11) and (11.9),  we find that   $F^{(i+1)}(\zeta ,\xi ) = 0$   for every   $\zeta \in {\Cal P}$.   Since  
${\Cal P}$  is a non-empty open subset of  ${\Cal R}$,  this implies that   $F^{(i+1)}(\zeta ,\xi ) = 0$   for every   $\zeta \in {\Cal R}$.  
Since this is true for every  $\xi \in {\Cal R}$,  we conclude that   $F^{(i+1)}$  is identically zero on    ${\Cal R}\times {\Cal R}$.  
This completes the induction on  $i$.
\medskip
In particular,  the above tells us that    $F^{(m)} = 0$  on   ${\Cal R}\times {\Cal R}$,  and consequently  
$$
\lim _{\nu \rightarrow \infty }F_{j_\nu }^{(m)}(0,0)  =  F^{(m)}(0,0)  = 0.  
\tag 11.12
$$
\noindent
Recalling (11.7),  we have
$$
F_{j_\nu }^{(m)}(0,0) = |r(z_{j_\nu }^{(0)})|^{n+1}\Phi \big(z_j^{(m)},z_j^{(0)}\big)  =   |r(z_{j_\nu })|^{n+1}\langle AK_{w_{j_\nu }},K_{z_{j_\nu}}\rangle  .
$$
\noindent
Since   $d(w_{j_\nu },z_{j_\nu }) < R$,  from  (4.1)  and Lemma  2.1   we obtain  
$$
|\langle Ak_{w_{j_\nu }},k_{z_{j_\nu}}\rangle  |  \leq  C_2|F_{j_\nu }^{(m)}(0,0)|.
$$
\noindent
This and (11.12) together contradict  (11.4).  This completes the proof.   $\square $
\medskip
\noindent
{\it Proof of Theorem} 1.2.  This follows immediately from Propositions  11.4 and 9.2.   $\square $

\bigskip
\centerline{\bf   References}
\medskip
\noindent
1.  S.  Axler  and  D.  Zheng,   Compact operators via the Berezin transform, 
Indiana Univ. Math. J. {\bf 47}  (1998),   387-400.

\noindent
2.  W.  Bauer  and  J. Isralowitz,   Compactness characterization of operators in the Toeplitz algebra of the Fock space $F^p_\alpha $, 
J. Funct. Anal. {\bf 263} (2012),   1323-1355.

\noindent
3.  D.  B\' ekoll\' e, C.  Berger,  L.  Coburn  and K.  Zhu, BMO in the Bergman metric on bounded symmetric domains,
J. Funct. Anal.  {\bf 93} (1990),  310-350.

\noindent
4.   C. Berger and  L. Coburn,  On Voiculescu's double commutant theorem, Proc. Amer. Math. Soc. {\bf 124}
(1996),  3453-3457.

\noindent
5.  C.  Berger,  L.  Coburn  and K.  Zhu,   Function theory on Cartan domains and the Berezin-Toeplitz symbol calculus, 
Amer. J. Math.  {\bf 110} (1988),  921-953.

\noindent
6.   K. Davidson, On operators commuting with Toeplitz operators modulo the compact
operators,   J. Funct. Anal.  {\bf 24}  (1977), 291-302.

\noindent
 7.   M.  Didas,  J.  Eschmeier  and  K.  Everard,  
On the essential commutant of analytic Toeplitz operators associated with spherical isometries, 
J. Funct. Anal. {\bf 261} (2011),   1361-1383.

\noindent
8.   X. Ding and S. Sun,  Essential commutant of analytic Toeplitz operators,   Chinese Sci. Bull.
{\bf 42} (1997), 548-552.

\noindent
9.  J. Eschmeier  and K. Everard,  Toeplitz projections and essential commutants,      J. Funct. Anal.
{\bf 269}  (2015),  1115-1135.

\noindent
10.  C. Fefferman,   The Bergman kernel and biholomorphic mappings of pseudoconvex domains,  
Invent. Math.  {\bf 26}   (1974), 1-65.

\noindent
11.  K. Guo and S. Sun, The essential commutant of the analytic Toeplitz algebra and some
problems related to it (Chinese),  Acta Math. Sinica (Chin. Ser.)   {\bf 39}  (1996), 300-313.

\noindent
12.  J. Isralowitz,  M.  Mitkovski  and  B.  Wick, Localization and compactness in Bergman and Fock spaces, 
Indiana Univ. Math. J. {\bf 64}  (2015),   1553-1573. 

\noindent
13.  B. Johnson and S. Parrott,  Operators commuting with a von Neumann algebra modulo the set of compact operators, 
J. Funct. Anal.  {\bf 11}  (1972), 39-61.

\noindent
14.  N. Kerzman,   The Bergman kernel function. Differentiability at the boundary, Math. Ann. {\bf 195} (1972), 149-158.

\noindent
15.  S.  Krantz,   Function theory of several complex variables. Pure and Applied Mathematics. A Wiley-Interscience Publication, 
John Wiley \& Sons, Inc., New York, 1982.

\noindent
16.  P.  Muhly and J.  Xia,  On automorphisms of the Toeplitz algebra,  Amer. J. Math.  {\bf 122} (2000),  1121-1138.

\noindent
17.  J. Munkres,  Analysis on manifolds, Addison-Wesley Publishing Company, 
Advanced Book Program, Redwood City, CA, 1991.

\noindent
18.  M. Peloso, Hankel operators on weighted Bergman spaces on strongly pseudoconvex domains, 
Illinois J.  Math. {\bf 38}  (1994), 223-249.
  
\noindent
19.  S.  Popa,  The commutant modulo the set of compact operators of a von Neumann algebra,   
J. Funct. Anal. {\bf 71} (1987),  393-408.

\noindent
20.  R. M.  Range,   Holomorphic functions and integral representations in several complex variables. 
Graduate Texts in Mathematics  {\bf 108},   Springer-Verlag, New York, 1986.
  
\noindent
21.  D. Su\'arez, The essential norm of operators in the Toeplitz algebra on  $A^p({\bold B}^n)$, 
Indiana Univ. Math. J.  {\bf 56}  (2007), 2185-2232.
  
\noindent
22.  H. Upmeier, Toeplitz operators and index theory in several complex variables. Operator Theory: Advances
and Applications  {\bf 81},   Birkh\"auser Verlag, Basel, 1996.
  
\noindent
23.   D. Voiculescu, A non-commutative Weyl-von Neumann theorem, 
Rev. Roumaine Math. Pures Appl.  {\bf 21} (1976),   97-113.

\noindent
24.  J. Xia,  On the essential commutant of   ${\Cal T}$(QC),   Trans. Amer. Math. Soc.  {\bf 360} (2008),
1089-1102.

\noindent
25.   J. Xia,  Singular integral operators and essential commutativity on the sphere,  
Canad. J. Math.  {\bf 62}  (2010),  889-913.

\noindent
26.  J. Xia,  Localization and the Toeplitz algebra on the Bergman space,   J. Funct. Anal.
{\bf  269}  (2015), 781-814.

\noindent
27.   J. Xia,  On the essential commutant of the Toeplitz algebra on the Bergman space,   
 J. Funct. Anal.  {\bf 272}  (2017), 5191-5217.
 
\noindent
28.  J. Xia,    A double commutant relation in the Calkin algebra on the Bergman space, 
J. Funct. Anal. {\bf 274} (2018), 1631-1656.

\noindent
29.   J.  Xia  and  D. Zheng,  Localization and Berezin transform on the Fock space,  
J. Funct. Anal. {\bf 264}  (2013),   97-117.

\medskip
\noindent
Department of Mathematics, State University of New York at Buffalo, Buffalo, NY 14260

\medskip
\noindent
E-mail: yiwangfdu\@gmail.com

\medskip
\noindent
E-mail: jxia\@acsu.buffalo.edu

\end